
       \input amstex
\documentstyle{amsppt}
\baselineskip =24 pt plus 2 pt
\lineskip=15 pt minus 1 pt
\magnification =\magstep1
\lineskiplimit =14pt
\bigskip

\def \pslr{$PSL(2,\Bbb R)$}
\def\pr{\tilde{\pi}_r}
\def\epst{\tilde{\eps}}
\def\hr{H^2(\Bbb H,\nu_r)}
\def\tr{\text{tr}}
\def\lh{L^2(\Bbb H,\nu_r)}
\def\lbh{B(\lh)}
\def\eps{\epsilon}
\def\ovt{\overline{\otimes}}
\def\calg{\Cal L(\G)}
\def\lbg{\Cal L(\G)\overline{\otimes} B(L^2(F,\nu_r))}
\def\lfg{l^2(\G)\overline{\otimes} L^2(F,\nu_r)}
\def\lbf{B(L^2(F,\nu_r))}
\def\lf{L^2(F,\nu_r)}
\def \pt {\tilde{\pi}^0 _r}

\def\b{B(H_r)}
\def\nk{||k||_{\lambda ,r}}
\def\G{\Bbb {\Gamma}}
\def\h{\Bbb H}

\def\i{\int_ G}
\def\c{c_r}
\def\l{L^2 (G)}
\def\s{\sup\limits_{\zeta \in \Bbb H}}
\def\a{\hat A(z,\zeta)}
\def\ra{\rightarrow}

\def\dz{|d(z,\zeta)|^r}
\def\dn{d\nu _0 (\zeta)}
\def\k{k(z,\overline \zeta)}
\def\bh{\hat B(H_r)}
\def\g{\gamma}
\def\hg{\Bbb H /\Bbb {\Gamma}}
\def \e{\epsilon}
\def\n{||\ ||_{\lambda ,r}}
\def\ns{||\ ||_{\lambda ,s}}
\def\na{||A||_{\infty ,r}}
\def\an{\langle e^{\frac{s-r}{2}}_z ,
 e^{\frac{s-r}{2}}_z \rangle}
\def\z{\zeta}
\def\st{*_s}
\def\fm{|\langle\pi (g)\zeta,\eta \rangle _H|^2}
\def\d{\text{d}}
\def\Im {\text{Im}}
\def\ovl{\overline}
\def\car{\Cal A_r}
\def\cas{\Cal A_s}
\def\cat{\Cal A_t}
\def\lar{L^2(\car)}
\def\las{L^2(\cas)}
\def\lat{L^2(\cat)}
\def\vv{\vert\vert}
\def\tbrf{T^r_{B^{-1/2}_rf}}
\def\jsr{j_{s,r}}
\def\brs{B^{-1/2}_sB^{1/2}_r}
\def\lr{l_{s,r}}
\def\cbs{\Cal B_s}
\def\cbr{\Cal B_r}
\def\lbr{L^2(\cbr,\tau)}
\def\lbs{L^2(\cbs,\tau)}
\def\sinab{_{s \in (a,b)}}
\def\dds{\frac {\text{d}}
 {{\text d}s}}
\def\st{\ast_t}
\def\ss{\ast_s}
\def\cbt{\Cal B_t}
\def\lbt{L^2(\cbt,\tau)}
\def\arg{\text{\ arg \ }}
\def \g{\gamma}
\def\la{\lambda}
\def\us{^{-1,s}}
\def\ut{^{-1,t}}
\def\ps{\tilde{\psi}}
\def\lt{L^2 (M,\tau)}
\def\t{\theta}
\def\G{\Gamma}
\def\hcat{\hat\Cal A_t}
\centerline{On the $\G$-equivariant form of
 the Berezin's quantization of the upper
half plane}
\centerline{ by Florin G. R\u adulescu
\footnote{ Department
of Mathematics, University of Iowa, Iowa City, IA 52242}
\footnote { Member of the
 Institute of Mathematics of the Romanian Academy}}
\bigskip

\centerline{Introduction}

Let $\G$ be a fuchsian subgroup of \pslr.
In this paper we consider the $\G$-equivariant form of
 the Berezin's
quantization of the upper half plane which will
 correspond to a deformation
quantization of the (singular) space $\Bbb H/\G$.
 Our main result is that the
von Neumann algebra associated to the $\G-$ equivariant
 form
 of the quantization
is stable isomorphic with the von Neumann algebra
 associated to $\G$. Moreover
the  dimension of each algebra in the deformation
 quantization, as a
left module over the group  von Neumann algebra
 $\Cal L(\G)$, is a
linear function of the deformation parameter
(the "Planck constant").

Recall that the Muray-von Neumann construction, for  the dimension
of projective, left modules over type II von Neumann algebras
 with trivial
center, allows all positive real numbers as  possible
  value for the
dimension. Consequently
 the above isomorphism is meaningful for all values of the
deformation parameter.

This will be particularly interesting when $\G$ is the
group $\G=PSL(2,\Bbb Z)$.
We use  the terminology (introduced in [KD], [FR], see also [DV])
of von Neumann algebras $\Cal L(F_t), t>1$
 corresponding
to free groups with a (possible) fractional
 "number $t$ of generators"
(even if the group itself may not make sense).
In this case the von Neumann algebras associated to
 the equivariant form of the
Berezin quantization will be free groups von Neumann algebras
where the  "number of generators" is a bijective function of the
deformation parameter.

The difference between the Berezin quantization of
 the upper half plane
and its $\G$-equivariant form is easy to establish.
 In the classical case the
von Neumann algebras associated to the deformation are simply
isomorphic to $B(H)$, the algebra of all bounded
operators on a Hilbert space. In the equivariant
case these algebras are type $II_1$ factors. This is
 a consequence
of the formulae for  the traces  in this algebras:
 classically the trace
is an integral over $\Bbb H$ of the restriction to the diagonal
of the reproducing kernel,
 while in the $\G$-equivariant case
the trace is the
 integral over a fundamental domain of $\G$ in $\Bbb H$.
This last fact is in a particular a generalization of
 the computation
in [GHJ] of the type $II_1$ factor trace of a product
of two Toeplitz operators having automorphic forms as symbols.

There exists a remarkable analogy, at least at the formal level,
between  Rieffel's ([Ri])
 construction
for the irrational rotation  C$^\ast-$algebra and
the construction in this paper. In Rieffel's construction,
 a
deformation quantization for the torus $\Bbb T^2$ is
realized  (in a more involved manner),  starting
  from the lattice
$\Bbb Z^2$ in $\Bbb R^2$. The von Neumann algebras
 in the deformation
for $\Bbb T^2$
are all isomorphic to the hyperfinite $II_1$ factor,
 with the exception of the
rational values for the parameter.

 For any deformation quantization, with suitable
properties, there exists
an associated  $2-$cocycle in the Connes's cyclic cohomology
 of a certain "smooth"
subalgebra, for each parameter value. In this paper, we construct
 in a natural way,
a dense family of subalgebras on which the 2-cocycles live.
 This algebras are endowed with a norm
that  is a continuous analog of the $\sigma(l^1,l^{\infty})$
 norm on finite
matrices.

Surprisingly, the formulae defining the 2-cocycle are very similar
to the formulae in the paper of Connes and Moscovici ([CM]) where
cyclic cocycles are constructed from Alexander-Spanier cycles.
We will also show that the two cocycles appearing in the deformation
quantization for $\Bbb H/\G$ are related in a natural to a canonical
element (see [Gr], [Gh]) in
the second bounded cohomology group $H^2_{\text{bound}}(\G,\Bbb Z)$.

Our construction, when $\G=PSL(2,\Bbb Z)$,
 could be thinked of as yet another
definition for  the von Neumann algebras corresponding to
 free groups with
"fractional number of generators". A possible advantage
of such an approach could be the fact that  the
"smooth structure" (see [Co])
may be used to
 define a  family of $C^{\ast}$-algebras that could be
 a candidate for  the $C^{\ast}$ -algebras
corresponding to free groups with "fractional number of generators".
It would be very interesting if one could establish a
direct relation between Voiculescu's  random matrix model for free
group factors and the "continuous matrix" model coming from
Berezin's equivariant quantization.

 We    prove that, for deformation
quantization of algebras   with the property that
 the associated cyclic 2-cocycles are
bounded with respect to the uniform  norms on the algebras,  there
exist a time-dependent, linear differential equation, whose
 associated
evolution operator  induces an isomorphism between the algebras
associated to distinct values of the deformation parameter.
This  depends on a rather standard technique to prove vanishing
of cohomology groups of von Neumann algebras by fixed point
 theorems.

The cocycles we are constructing would be bounded if
 a certain bounded
function on $\Bbb H^2$ ( defined by $z,w\ra\arg(z-\ovl w)$) would
 be a Schurr
multiplier (see [Pi],[CS]) on the Hilbert spaces of analytic
 functions
 $H^2(\Bbb H,y^{r-2}\text{d}x\text{d}y)$.

The proof of the main result of this paper is based on
 an observation,
related to the Plancherel formula for the universal
cover
$\widetilde{PSL(2,\Bbb R)}$ of \pslr\ established
 by L. Pukanzsky in [Pu] (see also [Sa]).
The fact we are using is  that the projective, unitary representations
of \pslr\ in the continuous series that
extends the analytic discrete series of \pslr\,
 also have square integrable coefficients over \pslr\,
 as the representations
in the
discrete series do.

 The  computation of  the isomorphism class of the algebras
in the deformation is then completed
  by  a  generalization,
to  projective, unitary representations,
of a method found in [GHJ] (and referred there  to [AS],[Co]).
 This method computes the Murray-von Neumann dimension
 of the Hilbert
 space of a representation of
a Lie group as left module (via restriction) over the
 group algebra
of a lattice subgroup.

\proclaim{Acknowledgment} Part of this paper was
elaborated while the author benefited of the generous
hospitality of the Department of Mathematics at the University
of Toronto. The author gratefully acknowledges
the very enlightening discussions he had during the
elaboration of this paper with A. Connes, G. A. Elliott,
 P. de la Harpe, V.F.R. Jones,
J. Kaminker, S. Klimek, V. Nistor, G. Pisier, S. Popa,
 M. A. Rieffel,
 G. Skandalis,
D. Voiculescu, A. J. W. Wasserman.
\endproclaim

\def\bbb{\Bbb}
\centerline{0. Definitions and outline of the proofs}

Recall that a von Neumann algebra is a selfadjoint subalgebra of
$B(H)$, which is unital and closed in the weak operator topology.
 A type $II_1$ factor is a von Neumann algebra $M$ with trivial
center and
such that there exists a
weakly continuous, linear functional (called trace)
 $\tau:M\ra\Bbb C$  with
$\tau(xy)=\tau(yx)$ and so that $0$ is  the only  positive element
in  the kernel of $\tau$. We  normalize $\tau$ by
$\tau(1)=1$.  Let $\Cal L(\G)$ be the weak closure of the
group algebra
$\Bbb C(\G)$,  represented in $B(l^2(\G))$,
 by left convolution operators. If $\G$ has  nontrivial,
infinite conjugacy classes, then $\Cal L(\G)$ is a type
$II_1$ factor, the trace being simply the evaluation at the
neutral element in $\G$.

 Such  algebras are usually associated with
 a discrete group $\G$ with infinite (nontrivial)
conjugacy classes. Let $\Cal L(\G)$ be the weak closure of the
group algebra
$\Bbb C(\G)$,  represented in $B(l^2(\G))$,
 by left convolution operators.

The following construction goes back to the
 original paper ([MvN]) of Murray and von Neumann. Let $M$ be a type
$II_1$ factor with trace
$\tau$ and  let $t$ be a positive
 real number. Denote by $\tau$ also, when no confusion
arises,
 the tensor product trace $\tau \ovl\otimes
\text{tr}_{B(H)}$ on $M\ovl\otimes B(H)$.
 Let $p$ be any selfadjoint idempotent
in $M\ovl\otimes B(H)$. Then the isomorphism class of the
type $II_1$ factor $p(M\ovl\otimes B(H))p$
is independent on the choice of $p$ as long as $p$ has trace $t$.
This type $II_1$ factor is usually denoted by
$M_t$. Clearly $\Cal F(M)=\lbrace t\vert\ M_t\cong M\rbrace$
is a multiplicative subgroup of $\Bbb R_{+}/\{0\}$, referred to,
by Muray and von Neumann as the fundamental group of $M$.

In the same paper referred above,
  given a weakly
continuous representation
of $M$ into some $B(K)$,
 the authors define   a positive real number $\text{dim}_MK$
which measures the dimension of $K$ as a left Hilbert module over
$M$.  The  dimension, in type $II_1$ factors, may take
 any positive real value. The original
terminology for $\text{dim}_MK$
 was the coupling constant of $M$ in $K$.

This number has all the formal features of   a dimension theory,
 that is $\text{dim}_M(K_1\oplus K_2)=
\text{dim }_M(K_1)+\text{dim }_M(K_2)$. The dimension number
 is normalized
 (when $M$
 is the von Neumann algebra $\Cal L(\G)$ of a group $\G$)
 by the condition
$\text{dim}_{\Cal L(\G)}{l^2(\G)}=1$. More generally, for  arbitrary
 $M$, let
$L^2(M,\tau)$ be  the Hilbert space corresponding
 to Gelfand-Naimark-Segal
construction for the trace $\tau$ on $M$, that is $L^2(M,\tau)$
is the Hilbert space
completion of $M$ as a vector space with respect to the scalar
product $\langle a,b\rangle_\tau=\tau(b^{\ast}a), a, b\in M$.  Then
we have
$\text{dim}_M(L^2(M,\tau))=1$.

To obtain a Hilbert space of an arbitrary dimension $t$,
 $0<t\leq1$ one takes
a projection $e'$ in the commutant
 $$M'=\lbrace x \in B(L^2(M,\tau))\vert \lbrack x,M\rbrack=0\rbrace,$$
of trace $\tau_{M'}(e')=t$. Then $e'L^2(M,\tau)$ is a left
Hilbert module over $M$ of dimension $t$. Note that $M'$ is
also a type $II_1$ factor. This last statement is more transparent
 in the case
 $M=\Cal L(\G)$. In this case $M'$ is  isomorphic and
 antisomorphic to $M$. In fact, $\Cal L(\G)'$ is
the von Neumann algebra $\Cal R(\G)$ generated
by right convolutors on $l^2(\G)$.
If $t>1$, to obtain  a module over
$M$ of dimension $t$ one has to replace
 $M$ by $M\otimes M(n,\Bbb C)$ where
$n$ is any integer bigger than $t$.

 From the above construction it is easy deduced that, if  $K$ is a
left Hilbert module over $M$ (that is we have
a unital embedding of $M$ into $B(K)$), then the commutant
$M'=\{x \in B(K)\vert \lbrack x,M\rbrack=0\rbrace$
 is antisomorphic to the algebra $M_t$
  with $t=\text{dim}_MK$.

Since for $M= \Cal L(\G)$, $M$ is
always antisomorphic to itself,  in this case
we may replace antisomorphic simply by isomorphic
everywhere in the above
statements.

Finally, one may   construct  a left Hilbert module over the
 von Neuman algebra
of a  discrete group $\G$ which is a lattice subgroup in a
semisimple Lie group $G$, for every representation $\pi$ of
$G$ which belongs to the discrete series representations of $G$.
Let
 $H_\pi$ be
 the Hilbert
 space of the representation and let
$\d_\pi$ be the coefficient with which this representation
 enters in the
Plancherel formula for $G$. Then the following result was
proved in [GHJ] (see also [Co],[AS]):

Let $M$ be the von Neumann  algebra generated by the image of $\G$
by $\pi$ in $B(H_\pi).$ If $\G$ is a lattice then $M$
 is isomorphic to
$\Cal L(\G)$ and
$$\text{dim}_MH_\pi=(\text{covol \ }\G)d_\pi.$$

We now recall the definition of the Berezin quantization from
[Be1], [Be2]. Although the setting could be more general, we will
restrict ourselves in this paper to the case
when the phase space is $\Bbb H$ with the geometry given
by the noneuclidian metric $(\Im\ z)^{-2}\d z\d\ovl z$.

The meaning of quantization was recalled  in the  introduction to
[Be1] and it was "a construction, starting from the classical
mechanics of a system, of a quantum system
which had the classical system as its limit as $h \ra 0$,
where  $h$ is Planck's constant" (we quote from the above mentioned
paper).

In the case when the phase space is $\Bbb H$, Berezin's
construction of an algorithm for the quantum system is the
construction of a family of  multiplications $\ast_h$ indexed
by $h$ on a suitable vector subspace of functions on $\Bbb H$,
which are associative, whenever this
comparison makes sense.

The property that the quantum system has the classical
 system as a limit,
when $h\ra 0$, means that for suitable functions $f,g$ on $\bbb H$
$$\lim\limits_{h \ra 0}\frac {1}{h}
\lbrack f\ast_hg-g\ast_h f\rbrack=
\lbrace f,g\rbrace,$$
where the Poisson bracket $\lbrace\cdot,\cdot\rbrace$ is computed
according to the formula:
$$\lbrace f,g\rbrace=(\Im z)^{-2}
\lbrack(\frac {\d}{\d z}f)(\frac {\d}{\d \ovl z}g)-
(\frac {\d}{\d z}g)(\frac {\d}{\d \ovl z}f)\rbrack.$$
The Berezin
algorithm for the multiplication operation
$\ast_h$ is realized by identifying any
suitable functions  on $\bbb H$ with
bounded linear operators
 on the Hilbert space of analytic functions
$$H_{1/h}=H^2(\bbb H,(\Im\ z)^{-2+1/h}\d z\d\ovl z).$$
The multiplication $\ast_h$ then corresponds, via this
 identification,
to the   composition operation for linear operators
on $H_{1/h}$.
There exists two possibilities to realize the correspondence
between functions on $\bbb H$ and linear operators on $H_{1/h}$.
In both methods the functions on $\Bbb H$ are identified with a special
type of symbol for  linear operators.

 Let $P_{1/h}$ be the
orthogonal projection from
$L^2(\bbb H,(\Im\ z)^{-2+1/h}\d z\d\ovl z)$ onto \break
$H^2(\bbb H,(\Im\ z)^{-2+1/h}\d z\d\ovl z)$ and let
$M^{1/h}$ be the multiplication operator on\break
$L^2(\bbb H,(\Im\ z)^{-2+1/h}\d z\d\ovl z)$ with the function $f$.
Then the  (not necessary bounded)
 Toeplitz operator with symbol $f$ is
$T^{1/h}_f=P_{1/h}M^{1/h}P_{1/h}.$

 A function $f$ on $\bbb H$
is the covariant symbol of a linear operator $A$ on
$H_{1/h}$ if $A$ is the Toeplitz operator $T^{1/h}_f$ on
$H_{1/h}$ with symbol $f$.  We will use
in this situation the notation $\overset\text{$\circ$}\to{A} =f$ or
$\overset\text{$\circ$}\to{A}(z,\ovl z)=f(z,\ovl z)$.

A function $f$ on $\bbb H$
is the contravariant symbol of a linear operator $A$ on
$H_{1/h}$ if $f$ is the restriction to the diagonal $z=w$ of
a function $\tilde f$ on $\bbb H^2$ which is analytic in the first
variable and antianalytic in the second variable.
 The relation between
$A$ and $\tilde f$ is explained bellow:
 Let $e^{1/h}_z, z\in \bbb H$ be the
evaluation vectors in $H_{1/h}$,
 that is $\langle f,e^{1/h}_z\rangle=
f(z), z \in \bbb H$, for all $f$ in $H_{1/h}$.
Then $\tilde f$ is given by
$$\tilde f(z,\ovl w)=\frac{\langle Ae^{1/h}_w,e^{1/h}_z\rangle}
{\langle e^{1/h}_w,e^{1/h}_z\rangle},z,w \in \bbb H.$$
 We will use the notation
$\tilde f=\hat A$ or $\tilde f(z,\ovl w)=
\hat A(z,\ovl w)$.

 The main theorem in Berezin's papers [Be1], [Be2] is that by using
the correspondence between functions on $\bbb H$ and linear
operators on $H_{1/h}$, given by any of this two symbols, one
gets a quantization which has the required classical limit.

Moreover there exist a natural duality relation between the two type
of symbols which is realized by using the pairing given
by the operatorial trace on $B(H_{1/h})$ : For suitable
bounded linear operators $A,B$ on  $H_{1/h}$ (see e.g. [Co] for a
rigorous treatment) one has
$$\leqno (0.1)\ \ \ \ \  \text{tr}_{B(H_{1/h})}(AB)=
\int\limits_{\bbb H}\hat A(z,\ovl z)
\overset\text{$\circ$}\to{B} (z,\overline z)
(\Im z)^{-2}\d z \d \ovl z .$$

Finally the multiplication operation $\ast_h$
 is  covariant with respect
to the action of group \pslr on $\Bbb H$. Hence
there exists a projective, unitary representation
$\pi_{1/h}$ on the Hilbert space $H_{1/h}$ with the
following property:
If $f=\hat A$ is the contravariant symbol of an operator
$A$ on $B(H_{1/h})$ then, for any group
element  $g$ in \pslr\, the function on $\bbb H$
defined by
$z\ra f(g^{-1}z)$
is the symbol of the operator $\pi_{1/h}(g)A\pi_{1/h}(g^{-1})$.

Let $\G$ be a discrete subgroup of \pslr\ of finite covolume.
The covariance property for the Berezin multiplication shows that
this operation is  an inner operation
 on a suitable space of  functions
on $\bbb H$, that are $\G-$ equivariant. From a borelian viewpoint,
$\G-$equivariant, measurable functions on $\bbb H$ are identified
  with functions on $F$ where $F$ is any fundamental domain in
$\bbb H$ for the action of $\G$ on $\bbb H$.

We will show in the third paragraph of this paper that there exists
a suitable  vector space $\Cal V_{1/h}$
 of $\G-$ equivariant functions on
$\Bbb H$, dense in $L^2(F)$
and so that $\Cal V_{1/h}$ is an involutive
  algebra with respect to the product
given by $\ast_{1/h}$.
The main theorem of this paper is the following

\proclaim{Theorem}Let $\G$ be a discrete subgroup
 of \pslr\ of finite covolume. Let $F$ be  a fundamental domain in
$\bbb H$ for the action of $\G$ on $\bbb H$.
For every $h>0$, there exists a
dense vector subspace $\Cal V_{1/h}$ of $L^2(F)$ which is
closed under conjugation and under
 the product $\ast_h$. Here
$\ast_h$  is the
Berezin's product, if we identify functions on $F$ with
$\G$-equivariant functions on $\bbb H$.
 Let $\tau$ be the functional on this
vector space defined by $\tau(f)=
\int\limits_Ff(z) (\text{Im \ }z)^{-2}\d z \d \ovl z .$ Then $\tau$
 is a trace,
that is $\tau(f\ast_hg)=\tau(g\ast_h f),$ for all suitable $f,g$.

Let  $\Cal A_{1/h}$ be the von Neumann algebra
obtained by taking the weak closure of the vector space $\Cal V_{1/h}$
 in the
Gelfand-Naimark-Segal representation associated with
 the trace $\tau$ on $\Cal V_{1/h}$. Then $\Cal A_{1/h}$  is
isomorphic to the type $II_1$ factor $(\Cal L(\G))_t$ with
$t=$ (covol $\G)(\frac{1/h-1}{\pi})$.
\endproclaim

This result is particularly interesting when $\G=PSL(2,\bbb Z)$.
In this case the algebras in the deformation are isomorphic to
the type $II_1$ factors $\Cal L(F_t)$
 corresponding to free groups with
real number of generators. This factors where introduced in
([KD] and independently in [FR]) based on the random matrix
 techniques developed
by Voiculescu (see also [DV] for an entropy theoretic  viewpoint
definition of this factors).

The interesting feature that appears is that
the (real) "number of generators"  $t$ in  $\Cal L(F_t)$, for the
algebras in the deformation, is a bijective function
 on the deformation
parameter (the Planck' constant). Recall that it is still
an open problem (hinted in [MvN] and first time explicitly mentioned
by R. Kadison) if the isomorphism class of $\Cal L(F_N)$
 depends on
$N$.

The explanation of this behavior when $\G=PSL(2,\Bbb Z)$ is that
in this case $\Cal L(PSL(2,\Bbb Z))=\Cal L(F_{\frac {7}{6}})$
(by [KD]) and  that the following formula holds:
$\Cal L(F_t)_r=\Cal L(F_{(t-1)r^{-2}+1}), $ for all $t>1$,
$r>0$, ([KD],[FR]).

The proof of the main theorem  follows  from the following facts.
By the covariance property recalled above is easy deduced that
$\Cal A_{1/h}$ is isomorphic to the commutant $\{\pi_{1/h}(\G)\}'$
of the image of $\G$ through $\pi_{1/h}$ in $B(H_{1/h})$. If
the group  cocycle corresponding to the projective representation
$\pi_{1/h}$ vanishes by restriction to
 $\G$ then $\pi_{1/h}\vert_{\G}$
may be perturbed (by scalars of modulus 1) to a representation
of $\G$.

We are then
in a situation which is very similar to the theorem we recalled
at the beginning of this paragraph.
 If $1/h=r$ is an integer then $\pi_{1/h}$ is a representation
in the discrete series of \pslr\ with coefficient $(r-1)/\pi$.

  Hence,
if $1/h$ is an integer, then
 the dimension of $H_{1/h}$ as a left Hilbert module over
$\Cal L(\G)$ (via $\pi_{1/h}\vert_{\G}$) is $(r-1)/\pi$(covol $\G$).
Consequently, by what we recalled at the
beginning of this paragraph, the
algebra $\Cal A_{1/h}$, which is the commutant of $\pi_{1/h}(\G)$
in $B(H_{1/h})$, will be isomorphic to
$\Cal L(\G)_t$ with $t=(r-1)/\pi$(covol $\G$).

If $1/h$ is now  a  positive real number, not necessary an integer,
then the projective  representation $\pi_{1/h}$ is no longer a
representation, so the above argument does no longer apply. There
exists still a striking similarity with the previous situation
which may be read off from the Plancherel formula for
the universal cover $\widetilde{PSL(2,\bbb R)}$ of \pslr.

 The
projective representations $\pi_{1/h}$ lift to actual unitary
representations of $\widetilde{PSL(2,\bbb R)}$ and they now belong,
as it was observed in Pukanszky article ([Pu], see also [Sa]),
to the continuous series of representations
 of $\widetilde{PSL(2,\bbb R)}$.
The coefficient with which the representation
$\pi_{1/h}$ intervene in the
continuous series  is given by the same algebraic formula
 as in the integer case for \pslr,
that is $(1/h-1)/\pi$.

The above mentioned property for the representations
$\pi_{1/h}$ may be better understood if we look directly to the
computations involved in determining the coefficient of
a representation in the Plancherel formula for the
discrete series.

If $1/h$ is an integer then
 the representation $\pi_{1/h}$ has square
summable coefficients which are verifying the generalized
Schurr orthogonality relations (see [Go],[HC]).
 This relations are (with $d g$
Haar measure on \pslr):
$$\int\limits_{PSL(2,\bbb R)}
\vert\langle\pi_{1/h}(g)\zeta,\eta\rangle\vert^2\d g=
\frac{1/h-1}{\pi}
\vv\zeta\vv ^2  \vv \eta\vv^2, \zeta,\eta \in H_{1/h}.$$
If $1/h$ is not an integer then one can still check this
 relations holds true for the projective representation
$\pi_{1/h}$. Note that this doesn't depend on the possible choice
of scalars of modulus 1 which would appear if we consider
$\pi_{1/h}$ be induced from a representation of the universal cover.

This fact means that each of the projective, unitary
   representations
$\pi_{1/h}$ for \pslr\  is contained in a "skewed" form
of the left regular representation of \pslr. One can deduce from here
an analogue of the theorem from the book [GHJ] holds true for the
representations $\pi_{1/h}$  and hence that
$\Cal A_{1/h}=\{\pi_{1/h}(\G)\}'$  is isomorphic to
$\Cal L(\G)_t$ with $t=$ (covol $\G)(\frac{1/h-1}{\pi})$.

Note that in the preceding setting, if $f$ is a $\G$-equivariant
function then the Toeplitz operator with symbol $f$ in
$B(H_{1/h})$ commutes with $\pi_{1/h}(\G)$ and hence it belongs
(for suitable functions $f$) to $\Cal A_{1/h}$.

The duality  relation (0.1) between the covariant and contravariant
symbol for operators in $\Cal A_{1/h}$ still exists if one
replaces the operatorial trace with the trace on the type
$II_1$ factor $\Cal A_{1/h}$. The relation takes now the
following form:
$$\leqno (0.2)\ \ \ \ \ \tau(AB)=
\int\limits_{F}\hat A(z,\ovl z)
\overset\text{$\circ$}\to{B} (z,\overline z)
(\Im z)^{-2}\d z \d \ovl z .$$
This is now a generalization of the formula 3.3.e in [GHJ], computing
the
 trace of a product of two Toeplitz
 operators having symbols automorphic forms.
In fact the covariant symbols for operators in $\Cal A_{1/h}$
may be regarded as a generalization of automorphic forms (in fact any
pair of automorphic forms gives rise to such a symbol which could
be eventually the symbol of an unbounded operator).

With the terminology we have just introduced our main result
also shows that
\proclaim{Proposition} Let $\pi_r$ be the projective representations
of \pslr on the Hilbert space $H_r=
H^2(\bbb H,(\text{Im\ }z)^{-2+r}\d z\d\ovl z)$ which are
 given by the same
formula as the unitary representations in the discrete
 analytic series
for \pslr\ when $r$ is an integer. Assume that $\G$ is a lattice
subgroup of \pslr\ so that the second group cohomology cocycle
for the projective representation $\pi_r$ of \pslr\
vanishes by restriction to $H^2(\G,\bbb T).$

Then von Neumann algebra $M$ generated by $\pi_r(\G)$ in
$B(H_r)$ is isomorphic to $\Cal L(\G)$ and
$\dim _MH_r=\frac{r-1}{\pi}(\text{covol\ }\G)$.
\endproclaim

In the remaining part of the paper we will be concerned with
certain cohomology classes that are associated with a deformation of
algebras. Formally to any deformation quantization one could
associate a 2-Hochschild cohomology cocycle $ a,b\ra a\ast'_rb$
 (which lives on a
suitable subalgebra on which derivations are possible) defined by
$$a\ast'_r b=
\frac{\d}{\d r}(a\ast_rb).$$

That this formally verifies  the properties of a 2-Hochschild cocycle
can be seen easy by taking the derivative in the
deformation parameter of the relation expressing
the associativity of the multiplication.

This 2-cocycle should measure, in a certain sense,
 the obstruction for the
algebras in the deformation to be isomorphic. In particular
if this element vanishes in the second cohomology group,
 then one could
hope to find an eventually unbounded operator
$X_r$ for all $r$ so that
$$\leqno (0.3)\ \ \ \ \  X_r(a\ast_rb)-(X_ra)\ast_r b-a\ast_r(X_rb)=
a\ast'_r b,$$
for suitable $a,b$.

If this operator could be made selfadjoint and if the evolution
 equation
$$\dot y(r)=X_r y(r),$$
would have a solution for a dense set of initial values then let
the associated evolution operator be  $U(s,t)$. Recall
that $U(s,t)$   is defined by the
condition that
$U(s,t)y$ is the solution of the differential eqaution at point
$s$ with initial $y$  condition at $t$. Then $U(s,t)$  would be an
isomorphism of algebras.
Indeed we would have formally that
if $\dot y(r)=X_r y(r)$ and $\dot z(r)=X_r z(r),$ then
$$\frac{\d}{\d r} {(y(r)\ast_rz(r))}=
\dot y(r)\ast_r z(r)+y(r)\ast'_r z(r)+
y(r)\ast_r\dot z(r).$$

The identity (0.3) would then show that the last expression can
be further reduced to $X_r(y(r)\ast_rz(r))$. The unicity of the
solutions of the differential equation shows that
 the evolution operator
must then be an isomorphism of algebras.

Of course to make this formal argument work properly,
 a lot of conditions about
domains of unbounded operators should be checked and
 this seems practically
impossible.

A possible approach to overcome this difficulties in particular
 cases would
be to use the quadratic forms instead of looking at
bounded operators.
This amounts to looking at
$$\frac{\d}{\d r}\langle (a \ast_r b,c\rangle_{\tau}=
\frac{\d}{\d r}\tau(a\ast_r\ast_rc^\ast),$$
for convenient $a,$ $b,$ $c$ rather
 then looking at $\frac{\d}{\d r}(a\ast_rb)$.

If we renormalize   the previous trilinear functional by discarding
the terms which come from the "skewing effect"  due to the fact
that the trace of a product of two elements depends on the
deformation parameter parameter, we will end up with a
cyclic cohomology two cocycle $\psi_r$   associated
to the deformation.

We introduce a new  norm $\vv\cdot\vv{\lambda,r}$ for each $r$ on
 the algebras in the deformation and denote the set of all
elements in $B(H_r)$ that are finite with respect to this norm by
$\hat{B(H_r)}$. This norm is the analogue of the usual
$\sigma(l^1,l^{\infty})$
norm on finite matrices and is defined for general
operators on $H_r$ by considering their kernels  to be "continuous
matrices" (see paragraph 2 for the precise definition of the norm
$\vv\cdot\vv{\lambda,r}$.
The remarkable property of the norms $\vv\cdot\vv{\lambda,r}$
is
\proclaim{Proposition}
 The algebras $\hat{B(H_r)}$ consisting
of all
elements in $B(H_r)$ that are finite with respect to this norm
  are involutive Banach algebras
not only with respect to the product $\ast_r$ but
also with respect to
all the products $\ast_s$ for $s\geq r$. Moreover
$\hat{B(H_r)}$ is a (dense) vector subspace of $\hat{B(H_s)}$.
\endproclaim

\proclaim{Proposition} We
 consider the subalgebras $\hat{\car}=\car\cap\hat{B(H_r)}$.
Then  the derivatives involved in the definition
of $\psi_r$  make sense for all $a,b,c$ in $\hat{\car}$  and
$\psi_r$ is bounded with respect to this norm.
More precisely
$$\vert\psi_r(a,b,c)\vert\leq\text{const}
\vv a\vv_{\lambda,r}\vv b\vv_{\lar}\vv c\vv_{\lar},a,b,c \in
\hat{\car}.$$
\endproclaim
In the last paragraph of this paper we will show that

\proclaim{Theorem}
 A cyclic
two cocycle $\psi$ on an arbitrary type $II_1$ factor,
 which has the property that
$$\psi(a,b,c)
\leq \vv a\vv_{\infty}\vv b\vv_2 \vv c \vv_2,a,b,c \in M,$$
vanishes in the second cyclic cohomology group
 $H^2_{\lambda}(M,\Bbb C)$ (see  A. Connes' article [Co] for the
definition of this groups).
  Moreover
in this case  the cocycle implementing $\psi$ could be
chosen so as to correspond to an antisymmetric
 bounded operator on both
$L^2(M,\tau)$ and $M$.
\endproclaim

But then such a solution would
lead to a linear differential equation with bounded linear
operators which is known to have a well defined evolution operator.

Thus if we  have a deformation quantization of algebras
in which the  cyclic cohomology cocycle verifies the  more
restrictive boundedness condition mentioned above,
 then the associated evolution
operator  induces an isomorphism between the algebras in the
deformation, at different values of the parameter.

In the case of the cocycles $\psi_r$ that arise in connection with
the equivariant form of the Berezin quantization, it is not clear
if one could prove  that $\psi_r$ is  such
a bounded linear functional  as above.

The best we can do is to write down an explicit formula
for $\psi_r(A,B,C)$ in terms of the symbols of the operators
$A,B,C$.
The  formula for $\psi_r(A,B,C)$ is the same as the formula
$\tau(ABC)$. The difference between the two formulas is that
for $\psi_r(A,B,C)$ one has to juxtapose to the integrand
defining $\tau(ABC)$ an Alexander Spanier cocycle
$\theta$ which is a  diagonally $\G-$ equivariant
function on $\bbb H^3$.
If
$$\leqno(0.4) \ \ \phi(z,\ovl{\zeta})=
i\text{arg}((z-\ovl{\zeta})/2i)= \text{ln}
((z-\ovl{\zeta})/2i)-
\ovl{\text{ln}((z-\ovl{\zeta})/2i)}, \ \ z,\zeta\text{ in }
\h$$
then $\theta$ has the expression
$$ \ \theta(z,\eta,\zeta)=
\phi(z,\ovl{\zeta})+\phi(\zeta, \ovl{\eta})
+\phi(\eta,\ovl{z}),\ \ z,\eta,\zeta \text{ in } \h.$$

The cocycle $\psi_r$ would be bounded if one could prove that
the function $\phi$ is a bounded Schurr multiplier on
the space $B(H_r$ for $r$ is in an interval.

If $\G=PSL(2,\bbb Z)$ then is easy to see that $\theta$
is vanishing as a $\G-$ equivariant
cocycle (in the $\G-$equivariant Alexander-Spanier cohomology)
 because
$$\leqno (0.5)\ \ \ \ \  \theta(z,\eta,\zeta)=
\tilde{\phi}(z,\ovl{\zeta})+\tilde{\phi}(\zeta, \ovl{\eta})
+\tilde{\phi}(\eta,\ovl{z}),\ \ z,\eta,\zeta \text{ in } \h,$$
with a  diagonally, $\G-$ equivariant  $\tilde\phi$. The
formula for $\tilde\phi$ is
$$\tilde{\phi}(z,\ovl\zeta)=\text{arg}((z-\ovl{\zeta})/2i)
+\arg(\Delta(z))-\ovl{\arg(\Delta(\zeta))}
, \ \ z,\zeta\text{ in }
\h.$$

Here  $\Delta$ is the unique automorphic form
of order 12. The disadvantage for  $\tilde\phi$ is that
$\tilde\phi$ is not a bounded function although
$\theta$ is.  The fact that one can not find
a bounded $\G-$equivariant $\tilde\phi$ solving the
above equation is related to the non-vanishing
of $H_{\text{bound}}^2(PSL(2,\bbb Z)$ (see [Ghys]).

\def \pslr{$PSL(2,\Bbb R)$}

\centerline{1. Berezin quantization of the upper half plane}
\bigskip

\bigskip
In this paragraph we recall some facts concerning the Berezin's
quantization of the upper half plane ([Be.1], [Be.2], [Upm]).
Berezin realizes the deformation quantization for
  the upper half plane (and
in fact for more general symmetric domains)
by using symbols for bounded
operators acting on Hilbert spaces of analytic functions. As the
 bounded operators
on Hilbert spaces of analytic functions
 are allways given by reproducing
kernels, this symbols will allways exist for bounded operators.

We let $H_r$ be the Hilbert space of square
 integrable analytic functions
on, the upper half plane $\Bbb H$ with respect to
 the measure $\nu _r$,
which has density $(\text{Im\ }\  z)^{r-2}$ with
 respect to the canonical Lebesgue
measure $dzd\overline {z}$ on $\Bbb H$.  $H_r$
 is nonzero for $r>1$ ([Ba]).

The choice of the measure  $\nu _r$ for
 the Hilbert spaces $H_r$ is dictated
by the fact that this Hilbert spaces are
 enacted with projective, unitary
representations $(\pi_r)_{r>1}$ of
 $PSL(2,\Bbb R)$. We recall first from
([Ba], see also [Sa], [Puk]) the
 construction of this representations.
Also recall that $\nu_0$ is a invariant
 measure on $\Bbb H$ under the
action of $SL(2,\Bbb R)$.

\proclaim {Definition 1.1}
([Ba]) Let $H_r$, $r>1$ be the Hilbert space of all
 analytic functions with
$||f||_r^2 <\infty$ where the Hilbert norm is defined by

$$||f||_r^2  =\int_{\Bbb H}
|f(z)|^2 (\text{Im\ }\ z)^{r-2}dzd\overline z =
\int_{\Bbb H} |f(z)|^2 d\nu _r (z).$$
Assume
 that $G=PSL_2 (\Bbb R)$ acts on $\Bbb H$ by M\" oebius transforms:
$$G\times {\Bbb H} \ni (\pmatrix a&b
\\c&d \endpmatrix ,\ z)\rightarrow \frac {az+b} {cz+d}\in \Bbb H$$
and let: $j(g,z)=(cz+d)$. Also choose (see [Ma, pag 113]) a normal
 branch
of $\text{arg} (j(g,z))=\text{arg}(cz+d)$, for
 all $z\in \Bbb H$, $g=\pmatrix a&b\\c&d \endpmatrix \in G$
where the function arg  takes its values in $-\pi <t\leq \pi$.

Using this branch for $(cz+d)^r =\text{exp} (r \text{ln} (cz+d))$
 one defines
$$(\pi_r (g)f)(z)=(cz+d)^{-r}f(g^{-1}z),\ g=
\pmatrix a&b\\c&d \endpmatrix \in G,\ f\in H_r,\ z\in \Bbb H \ .$$
Then([Ba]) $\pi_r : PSL_2 (\Bbb R)\rightarrow B(H_r)$
 is a projective, unitary
representation of $PSL_2 (\Bbb R)$  with
 cocycle $c_r(g_1 ,g_2)\in \{z\vert\  |z|=1,\ z\in \Bbb C\}$
defined by
$$\pi_r(g_1 ,g_2)=c_r(g_1 ,g_2)\pi_r(g_1)\pi_r(g_2),\  g_1,g_2\in G.$$
If $r=2,3,...$ then $\pi_r$ is an actual
 representation of $PSL_2 (\Bbb R)$
and belongs to the discrete series of representations
 for $PSL_2 (\Bbb R)$
(see[La]).
\endproclaim
Recall that any Hilbert space of analytic functions
 has a naturally associated
reproducing kernel ([Aro]). For $H_r$ this has been
 computed allready by
[Ba] and we recall the formulae.
\proclaim {Theorem 1.2} ([Ba]).
 The reproducing kernel $k_r(z,\zeta )$ for $\Bbb H$
is given by the formula
$$k_r(z,\zeta )=
\frac {c_r} {((z-\overline \zeta)/ 2i)^r}
 \text{for \ all} z,\zeta \in \Bbb H.$$
In particular the following functions on $\Bbb H$ defined
 for all $z\in \Bbb H$,
$$e^r_z (\zeta)=
\frac {c_r} {((\zeta-\overline z )/ 2i)^r} ,\ \zeta \in \Bbb H$$
belong to $H_r=H^2(\Bbb H,\nu _r)$ and
$\langle f,e^r_z\rangle _r =
f(z),\ \text{for all}\ f\in H_r,\ z\in   \Bbb H$.
\endproclaim
\proclaim {Corollary} 1.3 ([Ba]). The orthogonal
 projection $P_r$ from $L^2(\Bbb H,\nu_r)$
onto $H_r=H^2(\Bbb H,\nu _r)$ is given by the formula
$$(P_r f)(z)=\langle f,e^r_z\rangle =
c_r\int _{\Bbb H} \frac {f(\zeta)}
{((z-\overline \zeta)/ 2i)^r}d\nu_r (\zeta)$$
for all $f$ in $L^2(\Bbb D,\nu_r)$. In the  terminology of
 vector valueded integral, this is
$$P_r f=\int _{\Bbb H} f(\zeta)e^r_{\zeta}d\nu_r (\zeta).$$
\endproclaim
The fact that the Hilbert spaces $H_r$ have evaluation
 vectors $e^r_z,\ z\in \Bbb H$,
shows that all bounded linear operators $B$ on $H_r=H^2(\Bbb H,\nu _r)$
are given by integral kernels. Berezin defined the contravariant symbol
for an operator $B$ on $H_r$ to be its
 normalized integral kernel $\hat B(z,\overline \zeta )$,
 which is a function on
$\Bbb H^2$.
\proclaim {Definition1.4} ([Be1,2]).
 For $B$ in $B(H_r)$ let the contravariant symbol
 $\hat B=\hat B(z,\overline \zeta ) \ z,\zeta \in \Bbb H$
 be the function on $\Bbb H^2$, analytic in $z$,
 antianalytic in $\zeta$
defined by:
$$\hat B(z,\overline \zeta )=
\frac {\langle Be^r_{\zeta},e^r_z\rangle}
{\langle e^r_{\zeta},e^r_z\rangle} , \ z,\zeta \in \Bbb H^2.$$
Then $\hat B$ completely determines $B$ by
 the formula:
$$(Bf)(z)=\langle B f,e^r_z\rangle =
\langle B(\int f(\zeta )e^r_{\zeta} d\nu_r(\zeta),e^r_z\rangle$$
$$=\int f(\zeta )\langle Be^r_{\zeta},e^r_z\rangle d\nu_r(\zeta)=
c_r\int _{\Bbb H} \frac {\hat B(z,\overline \zeta)f(\zeta)}
 {((z-\overline \zeta)/ 2i)^r}d\nu_r (\zeta) , z\in \Bbb H, \
f\in H^2(\Bbb H,\nu _r).$$

\endproclaim
The following properties for $\hat B$ are obvious
 consequences of the definition.
In particular they show that the above integral is absolutely convergent.
\proclaim {Proposition}1.5 ([Be]).
 Let $B$ be any bounded linear operator acting
on $H_r=H^2(\Bbb H,\nu _r)$ and let $\hat B=
\hat B(z,\overline \zeta ), \ z,\zeta \in \Bbb H$
be its symbol as above. Denote by $||B||_{\infty ,r}$ the uniform norm of
$B$ (as an element of $B(H_r)$). Then:

a) For all $\zeta$ in $\Bbb H$, the function on $\Bbb H$ defined by
$z\rightarrow c_r((z-\overline \zeta)/ 2i)^{-r}
   \hat B(z,\overline \zeta )$
belongs to $H^2(\Bbb H,\nu _r)$ and has $L^2$norm less than
 $c^{1/2}_r ||B||_{\infty ,r}(\text{Im\ }  \zeta)^{-r/2}$.

b) For all $z$ in $\Bbb H$, the function on $\Bbb H$ defined
 by $ \zeta \rightarrow c_r((z-\overline \zeta)/ 2i)^{-r}
 \hat B(z,\overline \zeta )$
belongs to $H^2(\Bbb H,\nu _r)$ and has $L^2$norm less
 than $c^{1/2}_r ||B||_{\infty ,r}
(\text{Im\ }  \zeta)^{-r/2}$.

c) If $B$ is identified with $P_r BP_r$ as an operator acting
 on $L^2(\Bbb H,\nu _r)$
then the formula
$$ (Bf)(z)=
c_r\int _{\Bbb H}\hat B(z,\overline \zeta)f(\zeta) d\nu_r (\zeta)$$
holds for all $f$ in $L^2(\Bbb H,\nu _r)$.

d) $|B(z,\overline z)\leq ||B||_{\infty ,r}$, for all $z$ in $\Bbb H$.

e) The symbol $\hat {B^*}=
B^* (z,\overline \zeta), z,\overline \zeta \in \Bbb H$
of $B^*$ (the adjoint of $B$) is given by the formula:
$$\hat {B^*}(z,\overline \zeta)=
\overline {B(\zeta,\overline z)}, z,\zeta \in \Bbb H.$$
\endproclaim
Proof. a) follows from the fact
 that $ c_r((z-\overline \zeta)/ 2i)^{-r}\hat B(z,\overline \zeta )$
is by definition $\langle Be^r_{\zeta},e^r_z\rangle$.
 Thus for fixed $\zeta$
in $H$, this is $(Be^r_{\zeta})(z)$. The $L^2$ norm of $Be^r_{\zeta}$ is,
by the boundedness of $B$, less than
$$||B||_{\infty , r} ||e^r_{\zeta}||=
c^{1/2}_r ||B||_{\infty , r} (\text{Im }\ \zeta)^{-r/2}$$
as $\overline  {B(z, \zeta)}$ is antianalytic in $\zeta$.
$$ ||e^r_{\zeta}||^2=\langle e^r_{\zeta},e^r_{\zeta}\rangle =
 e^r_{\zeta} (\zeta)=\frac {c_r} {(\text{Im} \ \zeta)^r}.$$
This completes the proof of a) and b) is similar.
 The point e) is obvious while
for point d) we observe that
$$|B(z,\overline z)|=
|\langle B e^r_{z},e^r_{z}\rangle \langle e^r_{z},e^r_{z}\rangle ^{-1}|$$
$$\leq ||B||_{\infty , r} ||e^r_{z}||^2 ||e^r_{z}||^{-2}=
||B||_{\infty , r}.$$
Point c) now follows from the fact that all
 $f$ in  $L^2(\Bbb H, d\nu _r)$ we have
$$c_r\int _{\Bbb H}\hat B(z,\overline \zeta)f(\zeta) d\nu_r (\zeta)=
\langle f,B(z,\cdot)\rangle _r =
\langle P_r\ f,\overline {B(z,\cdot)}\rangle _r $$
for all $z$ in $\Bbb H$, as $\overline {B(z,\cdot)}$ is analytic in $\zeta$.

We now recall the definition of the Berezin
 product of the symbols $\hat A,\ \hat B$, which
are functions on $\Bbb H^2$. The product will depend on the quantization
variable $r$, and it represents in fact the symbol of the composition
(in $B(H_r))$ of the bounded operators on $H_r$ that
 are represented by the symbols
 $\hat A,\ \hat B$.
\proclaim {Definition} ([Be 1]). Let $A,\ B$ be
 two functions on $\Bbb H^2$,
analytic in the first variable, antianalytic in the second and so that
$((z-\overline \zeta)/ 2i)^r A(z,\overline \zeta)$, as
 a function of $\overline \zeta$
on $H$, keeping $z$ in $\Bbb H$ fixed, belongs to $H^2(\Bbb H,\nu _r)$,
and so that
 $((z-\overline \zeta)/ 2i)^r B(z,\overline \zeta)$ as a function
of $z$, on $\Bbb H$, keeping $\zeta$ in $\Bbb H$ fixed,
 belongs to  $H^2(\Bbb H,\nu _r)$.
For $r>1$ the product of the two symbols $A,B$ is defined by the formula:
$$(A\ast _r B)(z,\overline \zeta) =
 c_r ((z-\overline \zeta)/ 2i)^r
\int _{\Bbb H}((z-\overline \eta)/ 2i)^{-r} A(z,\overline \eta)
 ((\eta-\overline \zeta)/ 2i)^{-r}
B(\eta,\overline \zeta) d\nu_r (\eta).$$
for all $z,\zeta$ in $\Bbb H^2$.
\endproclaim
The following formula is an easy consequence of
 the integral representation for
operators acting on $H_r =H^2(\Bbb H , \nu_r)$.
 It shows that the above product
is the symbol of the composition of the operators in $B(H_r)$.
\proclaim {Proposition} ([Be]).
 If $\hat A=\hat A(z,\zeta), \hat B=
\hat B(z,\zeta), z,\zeta \in \Bbb H^2$
are the symbols of two bounded
 linear operators $A$ and $B$, respectively,
acting boundedly on $H_r$,
 then $\hat A\ast _r\hat B=
(\hat A\ast _r\hat B) (z,\overline \zeta) ,z,\zeta \in \Bbb H^2$
is the symbol of the composition $AB$ in $B(H_r)$ of
 the two operators $A$ and $B$.
\endproclaim
Proof. We have to compute
 $\langle AB e^r_{\zeta},
e^r_{z}\rangle _r \langle e^r_{\zeta},e^r_z \rangle^{-1}_r$
for $z,\zeta$ in $\Bbb H$, and show that this is equal to
  $(\hat A\ast _r\hat B) (z,\overline \zeta)$,
if we know that $\hat A (z,\overline \zeta)=
\langle A e^r_{\zeta},e^r_z \rangle _r \langle e^r_{\zeta},
e^r_z \rangle^{-1}_r$
and similary for $B$. We have:
$$\langle AB e^r_{\zeta},e^r_z \rangle _r=
\langle A\int_{\Bbb H}(Be^r_{\zeta})(\eta)e^r _{\eta}
 d\nu _r (\eta),
e^r_z\rangle _r=$$
$$\int_{\Bbb H}
(Be^r_{\zeta})(\eta)\langle A e^r _{\eta},e^r_z\rangle d\nu _r (\eta)=
\int_{\Bbb H}\langle Be^r_{\zeta},
e^r _{\eta}\rangle _r \langle A e^r _{\eta},e^r_z\rangle d\nu _r (\eta)$$
and hence for all $z,\zeta$ in $\Bbb H$
$$
 \langle AB e^r_{\zeta},e^r_z \rangle _r \langle e^r_{\zeta},
e^r_z \rangle^{-1}_r=\langle e^r_{\zeta},e^r_z \rangle^{-1}_r =$$
$$\langle e^r_{\zeta},
e^r_z \rangle^{-1}_r \int_{\Bbb H}
\hat A(z,\overline \eta) \langle e^r_{\eta},e^r_z \rangle_r
 \hat B(\eta ,\overline \zeta) \langle e^r_{\zeta},
e^r_{\eta} \rangle d\nu _r (\eta)=$$
$$ c_r ((z-\overline \zeta)/ 2i)^r
\int _{\Bbb H}\hat A(z,\overline \eta)((\eta-\overline z)/ 2i)^{-r}
\hat B(\eta ,\overline \zeta)((\eta-\overline {\zeta})/ 2i)^{-r}
 d\nu_r (\eta).$$
Finnally, we  recall the fact that togather with the above contravariant
 symbols
for operators acting on $H_r$, Berezin also
introduced another type of symbols called covariant symbols for operators
acting on $H_r$. Contrary to the contravariant
 symbols for a bounded operator,
the covariant symbol
does not allways exists.
\proclaim {Definition} ([Be]). Let $A$ be a bounded operator
 on $H_r$ and
let $f$ be a bounded measurable function $f$ on $\Bbb H$.
 Let $M^r_f$ be
the bounded multiplication operator on $L^2(\Bbb H ,\nu _r)$ with the
function $f$ and let $T^r_f =P_r M^r_f P_r$ in $B(H_r)$ be
 the Toeplitz operator
with symbol $f$. Then $f$ is called the covariant symbol of
 $A$ and one uses
(following [Be]) the notation $f=\AA  $, if $A=T^r_f$.
 Clearly in this
case the uniform norm of $A$ is bounded by the essential
 norm of $f$.
\endproclaim
\bigskip

The relation between the two symbols for a given bounded
 operator $A$ on
$H_r$ is obtained as follows.

\proclaim {Proposition} ([Be]) If $f$ is any bounded measurable function
on $\Bbb H$, and $T^r_f =P_r M^r_f P_r$ is the
 corresponding Toeplitz operator
with symbol $f$, acting on the Hilbert space $H_r$,
 then the (contravariant)
symbol $\hat A=\hat A (z,\overline \zeta)$ for the
 operator $A=T^r_f$ is
$$\hat A (z,\overline \zeta)=
\langle e^r_{\zeta},e^r_z \rangle^{-1}
 \int_{\Bbb H} f(a) \langle e^r_a,e^r_z \rangle
 \langle e^r_{\zeta},e^r_a \rangle d\nu_r (a)=$$
$$c^{-1}_r \int_{\Bbb H} f(a)
\frac {\langle e^r_a,e^r_z \rangle
\langle e^r_{\zeta},e^r_a \rangle}
{\langle e^r_{\zeta},e^r_z \rangle
 \langle e^r_{a},e^r_a \rangle} d\nu_0 (a) =
c_r \int_{\Bbb H} f(a)
\frac {((z-\overline \zeta)/ 2i)^r
((a-\overline a)/ 2i)^r}
 {((z-\overline a)/ 2i)^r ((a-\overline \zeta)/ 2i)^r}
   d\nu_0 (a).$$
In particular $\hat A (z,\overline z)=
 (B_r f)(z,\overline z)$,
where $B_r$ is the operator on
$L^2(\Bbb H ,\nu _0)$ with kernel
 the function on $\Bbb H^2$ defined by
$(z,a)\rightarrow [\frac {\text{Im } \ z\ \text{Im }\ a}
 {|z-\overline a |^2}]^2 =k_r(z,a).$
\endproclaim
Denote $\rho (z,a)$ to be
$\frac {\text{Im } \ z\ \text{Im}\ a} {|z-\overline a |^2}$.
Then it is well known from hyperbolic geometry
 that $\rho$ is the inverse of
the hyperbolic cosinus of the hyperbolic distance
 between $z$ and $a$, for
$z,a$ in $\Bbb H^2$. In particular the kernel $k_r$
 is invariant under
the diagonal action of the group $PSL(2,\Bbb R)$ and
 hence by ([Se], see also [Ku])
the operator $B_r$ is a function $B_r (\Delta)$
  (in the sense of functional calculus
for unbounded selfadjoint operators) of the invariant
 laplacian $\Delta=
(\text{Im}\ z)^2 \frac {\partial}{\partial z}
\frac {\partial}{\partial \overline z}$
 on $\Bbb H$.

The explicit formula for $B_r=B_r(\Delta)$ as function
 of the laplacian is
determined also in [Be] as an infinite product
 of resolvants of $\Delta$.
As we will not need the explicit formula for $B_r$ we
 will only recall the
essential properties of $B_r$ as stated in [Be].

\proclaim {Theorem} ([Be]). The operator $B_r$is
 a bounded, positive operator on $L^2(\Bbb H ,\nu _0)$.
Moreover $||B_r||\leq 1$ and the operators
 $B_r$ pairwise commutte for $r>1$. Also $B_r$
converges strongly to 1 when $r \rightarrow \infty$.
\endproclaim
Finally we will recall the duality relation between the
 covariant and the
contravariant symbols. The rigurous form of this
 statements may be found in
the paper by Coburn [Co].

\proclaim {Theorem} ([Be]). Let $A$ be  a bounded
 linear operator on $H_r$
of contravariant symbol  $\hat A=
\hat A (z,\overline \zeta) \ z,\overline \zeta \in \Bbb H$
and let  $B$ be  another bounded
 linear operator of covariant symbol
 $f\in L^{\infty} (\Bbb H ,\nu _0)
,\ f(z)=
\overset\text{$\circ$}\to{B} (z,\overline z) ,\ z\in \Bbb H$.
 If
 the operator $AB$ is in the ideal of  trace
class operators, then
$$tr_{B(H_r)}(AB)=
\int _{\Bbb H}\hat A (z,\overline z)
 \overset\text{$\circ$}\to{B} (z,\overline z) d \nu _0 (z).$$
\endproclaim

The deformation quantization of the upper halfplane is now
 realized by
Berezin, by observing that for any two functions $f,\ g$ on
 the upper halfplane,
which are so that there exists functions
 $\hat A=\hat A (z, \zeta),\  \hat B=
\hat B (z,\overline \zeta), \ z, \zeta \in \Bbb H$,
analytic
in $z$ and antianalytic in $\zeta$ with
$$f(z)=\hat A (z, \overline z);\ g(z)=
\hat B (z, \overline z), \ z\in \Bbb H,$$
one may define
$$(f*_r g)(z)=(\hat A *_r \hat B)(z, \overline z),\ z\in \Bbb H.$$
This product is well defined for example as long as
 $\hat A, \hat B$ are symbols
of bounded operators acting on the spaces $H_r$.

The main result in [Be] is that under suitable conditions
 on the functions
$f,g$ the limits $\lim\limits_{r\rightarrow \infty} f*_r g$
and $\lim\limits_{r\rightarrow \infty} r(f*_r g - f*_r g)$
exists and are equal to respectively $fg$ and
$$(\text{Im}\ z)^2 (\frac{\partial f}
{\partial z} \frac{\partial g}{\partial \overline z}
- \frac{\partial g}{\partial z} \frac{\partial f}
{\partial \overline z}).$$
The last limit is the $SL(2,\Bbb R)$ invariant form
 of the Poisson bracket
on $\Bbb H$.

Finnally, all the above formulae hold true if one
 replaces $\Bbb H$ by
$\Bbb D$ the unit disk. The form of the formulae
 will be the same except
that we will have to replace the factor $(z-\overline \zeta)/{2i}$
 by
$(1-z\overline \zeta )$ (so that $\text{Im}\ \zeta$ will
 be replaced by $1-|\zeta |^2$).

In the case of the unit disk the group which replaces $SL(2,\Bbb R)$ is
$SU(1,1)$, which is the group of all matrices
$$\{\pmatrix a&b\\{\overline b}&{\overline a}
 \endpmatrix \vert \  a,b \in \Bbb C,
|a|^2 - |b|^2 =1 \}$$
The modular factor $(cz+d)$ is now
 replaced by $\overline b z+\overline a,\ z\in \Bbb D$.

\centerline{2. Smooth algebras associated to the Berezin quantization}
\bigskip

\bigskip
In this paragraph, by analogy with finite matrices, we construct a (weakly
dense) subalgebra $\hat B(H_r)$ of $B(H_r)$, which is a Banach algebra with
respect to a certain norm (the analog of the $(l^1,l^{\infty})$ norm
 on finite
matrices). We will show that the algebras $\hat B(H_r)$
 are well behaved,
globally, with respect to the Berezin product $*_s$ for
 all $s\geq r$.
In particular the symbols corresponding to the operators
 in $\hat B(H_r)$
will form an algebra under all the operation $*_s$.
 Thus, operations, such
as differentiation of $A*_s B$  will have a  sense for
 $A,B$ in $\hat B(H_r)$.

The norm which defines the algebras $\hat B(H_r)$ is
 the analog (modulo a
weight) of the supremum (after lines) of the absolute sum of elements
in all the rows of a given matrix.

We first start by a criteria (which is essentially contained in Aronszjan
memorium ([Aro])) on the contravariant
symbol  of  a bounded operator, for the operator  to be positive.

\proclaim {Lemma 2.1} Let $A$ in $B(H^2(\Bbb H ,\nu_r ))$ be a positive,
bounded operator on $H_r$ of uniform norm $||A||_{\infty ,r}$ and with
contravariant symbol
 $\hat A=\hat A (z,\overline {\zeta}) \ z, \zeta \in \Bbb H$.
Then there exists a constant
 $M>0$ so that the following matrix inequality
holds for all $N$ in $\Bbb N$ and
 all $z_1 ,z_2 ,...z_N$ in $\Bbb H$:
$$0\leq [\hat A (z_i, \overline {z_j})((z_i - \overline {z_j})/2i)^{-r}
]^N_{i,j=1}\leq [M((z_i - \overline {z_j})/2i)^{-r}]^N_{i,j=1}$$
Moreover, given $A$, the best constant $M$ for which
 the above inequality holds
for all $N,\ z_1 ,z_2 ,...z_N$ is
 the uniform norm $||A||_{\infty ,r}$.

Conversely, let $K$ be a kernel
 $k=k(z,\overline {\zeta})$ on $\Bbb H^2$,
analytic in $z$, antianalytic in $\zeta$ for which there exists a positive
constant $M$ such that the above inequality
 holds with $k$ replacing $\hat A$,
for all $N$ in $\Bbb N$ and all $z_1 ,z_2 ,...z_N$ in $\Bbb H$.

Then $k$ is the contravariant symbol of a positive bounded operator on
$\Bbb H$ of uniform norm less than $M$.
\endproclaim
Proof. We start with the direct part of our statement.
 The inequality we
have to prove is equivalent to showing (with $M=||A||_{\infty ,r}$)
that for all $f$ in
 $L^{2} (\Bbb H ,\nu _r)$ one has
$$0\leq \int \int_{\Bbb H^2} \langle A e^r_{\zeta},
e^r_z \rangle f(\zeta) \overline {f(z)} d\nu _r(z,\zeta)
\leq M
\int \int_{\Bbb H^2} \langle  e^r_{\zeta},
e^r_z \rangle f(\zeta) \overline {f(z)} d\nu _r(z,\zeta).$$
We identify $A$ with $P_r AP_r$ as an
 operator acting boundedly on $L^{2} (\Bbb H ,\nu _r)$.
 By proposition 1.5c  $A$ has as reproducing kernel,
 the function on $\Bbb H^2$,
$(z,\zeta)\rightarrow  \langle A e^r_{\zeta},e^r_z \rangle$.
 Hence the
above inequality is
$$0\leq \langle A f, f\rangle \leq M \langle P_r f, f\rangle
=M\langle P_r f, P_r f\rangle ,
f\ \text{in}\  L^{2} (\Bbb D ,\nu _r).$$
Since $A=P_r AP_r$ and $A$ is positive, of norm $M=
||A||_{\infty ,r}$,
this holds true.
The converse is along this lines. This completes the proof.
\bigskip

We introduce the following notation for a square root of
 the inverse of the
hyperbolic distance between points in the upper half plane.
\proclaim
{Notation 2.2}  For $z,\zeta$ in $\Bbb H$ denote
 by $d(z,\overline {\zeta})$ the quantity:
$$d(z,\overline {\zeta})=
\frac {(\text{Im}\ z)^{1/2} (\text{Im}\ \zeta)^{1/2}}
{[(z-\overline {\zeta})/2i]},\ \text{for\ all\ }
z,\zeta\ \text{in} \ \Bbb H$$
or in the unit disk representation
$$d(z,\overline {\zeta})=
\frac {(1-|z|^2)^{1/2}(1-|\zeta |^2)^{1/2}}
{1-z\overline {\zeta}},\text{for\ all\ }
 z,\zeta \ \text{in} \ \Bbb D.$$
Let $D_h (z,\zeta)$ be the hyperbolic distance
 between points in $\Bbb H$
(respectively $\Bbb D$). Then
$$|d(z,\overline {\zeta})|^2=\rho (z,\zeta)=
cosh^{-1}(D_h(z,\zeta)),\ \text{for\ all\ }
z,\zeta \ \text{in}\ \Bbb H .$$
\endproclaim
In the next lemma we will show that the growth of a symbol
 $\hat A=\hat A (z,\overline {\zeta})$,
when $z,\zeta$ are approaching the boundary of
 $\Bbb H$, is of the same
type (in absolute value) as the growth of $\rho (z,\zeta)^{-r/2}$ when
$z,\zeta$ are approaching the boundary.
 Moreover the estimate depends
only on the uniform norm of $A$.
\proclaim {Corrollary 2.3} Let $A$ in $B(H_r)$ be
 a bounded operator and
let $\hat A=\hat A (z,\overline {\zeta}) \ z, \zeta \in \Bbb H$ be
its (contravariant)
symbol. Then

i) if $A$ is positive then
$$| A (z,\overline {\zeta})||d(z,\zeta)|^r\leq (A(z,\overline {z}))^{1/2}
(A(\zeta,\overline {\zeta}))^{1/2}\leq
 \sup\limits_{w\in \Bbb H} A(w,\overline w),
z,\zeta \ \text{in}\ \Bbb H .$$
ii) If $A$ is arbitrary in $B(H_r)$ then
$$| A (z,\overline {\zeta}) (d(z,\zeta))^r|
\leq 4||A||_{\infty ,r}\
\text{ for all } z,\zeta \in \Bbb H .$$
\endproclaim
Proof. i) If $A$ is a positive element in
  $B(H_r)$, then by the preceding
proposition, for all $z,\zeta \ \text{in}\ \Bbb H$, the matrix
$$\pmatrix {A(z,\overline {z}) (\text{Im}\ z)^{-r}}
&{A (z,\overline {\zeta})(z - \overline {\zeta})/2i)^{-r}}\\{\ }
&{\ }
\\{\overline {A (z,\overline {\zeta})(z - \overline {\zeta})/2i)^{-r}}}&
{A(\zeta ,\overline {\zeta})
(\text{Im}\ \overline \zeta)^{-r}} \endpmatrix$$
is positive. (The expresion for the bottom left corner is justified by
Proposition 1.5.d.

It is obvious that if $(a_{i,j})^2_{i,j=1}$ is a positive matrix, then so
is $(a_{i,j}\ \lambda _i\overline {\lambda _j})^2_{i,j=1}$ for
 all $\lambda _1,\lambda _2$ in $\Bbb C$.
We let $\lambda _1 =(\text{Im}\  z)^{r/2},\ \lambda _2 =
(\text{Im}\ \zeta)^{r/2}$
and apply this statement to the above matrix.
 We obtain that the following matrix,
$$\pmatrix {A(z,\overline {z})}&
{A (z,\overline {\zeta})(d(z,\overline {\zeta}))^{r}}\\{\ }&{\ }
\\\overline{{A(z,\overline {\zeta})(d(z,\overline {\zeta}))^{r}}}
&{A(\zeta ,\overline \zeta)}\endpmatrix$$
is positive for all $z,\zeta$ in $\Bbb H .$ The condition that
 its determinant
be positive proves condition i).

To prove ii) we observe that by i) we already  know
that for all positive elements
$A$ in $B(H_r)$  one has
$$|A (z,\overline {\zeta}) (d(z,\zeta))^r|
\leq \sup\limits_{w\in \Bbb H} A(w,\overline w)
\leq ||A||_{\infty ,r}$$
(we use again Proposition 1.5.c). To get the general statement we use the
fact that any element $A$ in $B(H_r)$ has the expression ([Dix]).
$$A=[(\text{Re}\ A)_+ -(\text{Re}\ A)_- ]+
i[(\text{Im}\ A_+)-(\text{Im}\ A_-)]$$
where
 $(\text{Re}\ A)_{\underline +},(\text{Im}\ A)_{\underline +}$
 are positive of uniform
norm less then the uniform norm $||A||_{\infty ,r}$ of $A$.

The preceding proposition shows allready that the symbol
 $\hat A=\hat A (z,\overline {\zeta})$
of an operator $A$ in $B(H_r)$ will define, by the integral
 formula, a bounded
operator on all the Hilbert spaces $H_s$, with $s\geq r$.
 Moreover the
uniform norm doesn't increase.
\proclaim {Proposition 2.4} Fix $r>1$ and let $s\geq r$.
 Let $j_{s,r}$
be the map which assigns to an operator $A$ on $B(H_r)$ with
 contravariant
symbol $\hat A=
\hat A (z,\overline {\zeta}),\ z,\overline {\zeta}\in  \Bbb H$
the operator acting on
 $H_s$ with the same contravariant symbol, that is
$$(j_{s,r}(A) f)(z)=\int_{\Bbb H}\hat A (z,\overline {\zeta})
\langle  e^s_{\zeta},
 e^s_z \rangle f(\zeta) d\nu _s(\zeta) ,\ f\ \text{in}\ H_s.$$
Then $j_{s,r}$ takes its values in  $B(H_s)$ and
$$||j_{s,r}(A) ||_{\infty ,s}\leq 4 ||A||_{\infty ,r}$$
$$ j_{s,r}(B(H_r)_+)\subseteq B(H_s)_+ .$$
Moreover $j_{s,r}$ has a weakly dense image in  $B(H_s)$.
 We convene to
denote $j_{s,r}(A)$ also by $A$, this being justified by
 the fact that
both operators have the same symbol.
\endproclaim
Proof. As in the preceding proof, to show the estimate
on the norm, it is
sufficient to prove that $j_{s,r}(B(H_r)_+)\subseteq B(H_s)_-$
 and that
$$||j_{s,r}(A) ||_{\infty ,s}
\leq  ||A||_{\infty ,r} \text{ for } A \text{ in } B(H_r)_+ .$$
If $A$ belongs to $B(H_r)_+$, we
 have seen in Proposition 2.1, with
$M=||A||_{\infty ,r}$, that for all
 $N$ in $\Bbb N$, $z_1,z_2,...z_N$ in
$\Bbb H$ one has the following matrix inequality:
$$(2.1)\ \
0\leq [A(z_i ,\overline z_j)((z_i - \overline {z_j})/2i)^{-r}]_{i,j}
\leq M[((z_i - \overline {z_j})/2i)^{-r}]^N_{i,j}.$$
 From general matrix theory we know that if the  matrix  inequality
\break
$0\leq (a_{i,j})^N_{i,j=1} \leq (b_{i,j})^N_{i,j=1}$  holds for
 some
matrices $(a_{i,j})_{i,j=1},\ (b_{i,j})_{i,j=1}$ in $M_N(\Bbb C)$
 then so
does the matrix inequality
$$0\leq (c_{i,j}\cdot a_{i,j})^N_{i,j=1}
 \leq (c_{i,j}\cdot b_{i,j})^N_{i,j=1}$$
for any other positive matrix $(c_{i,j})^N_{i,j=1}$ in  $M_N(\Bbb C)$ .

On the other hand by ([Shapiro, Shielda]), for $s\geq r$ the matrix
$$[\frac {1}{[(z_i - \overline {z_j})/2i]^{s-r}}]^N_{i,j =1}$$
is allways positive for all $N$ in $\Bbb N$, $z_1,z_2,...z_N$ in
$\Bbb H$.
Hence the above remarks and formula (2.1) prove
 that for all $N$ in $\Bbb N$, $z_1,z_2,...z_N$ in
$\Bbb H$, the following matrix inequality holds:
$$0\leq [A(z_i ,\overline z_j)((z_i - \overline {z_j})/2i)^{-s}]^N _{i,j}
\leq  ||A||_{\infty ,r} [((z_i - \overline {z_j})/2i)^{-s}]^N _{i,j}.$$
The converse to Proposition 2.1 shows that $j_{s,r}(A)$
 is a positive
operator in $B(H_s)$ of uniform norm less than $||A||_{\infty ,r}$.

It remains to check that $j_{s,r}(B(H_r)$ is
 weakly dense in $B(H_s)$.
It is sufficient to check that $j_{s,r}(\Cal C_2(H_r)$ is weakly dense in
$B(H_s)$ (we use here the standard notation $\Cal C_2(H)$
 for the Hilbert-
Schmidt operators on $H$). In the canonical identification
 of $\Cal C_2(H_r)$
with $H_r\overline {\otimes} H_r$ (see [Sch],[B.S]) an
 operator $A$ in
$j_{s,r}(B(H_r))$  with symbol $\hat A=\hat A (z,\overline {\zeta})$ will
correspond to the function on $\Bbb H ^2$, defined by
$$(z , \overline {\zeta})\rightarrow
\hat A (z,\overline {\zeta} [(z - \overline {\zeta})/2i]^{-(s-r)}.$$
For fixed $\eta _1,\eta _2$ in $\Bbb H$ and $\epsilon >0$,
 the function
on $\Bbb H ^2$,defined by
$$l_{\eta _1,\eta _2 ,\epsilon}(z,\overline {\zeta}) =
 (z - \overline {\zeta})/2i]^{(s-r)}
[(z - \overline {\eta _1})/2i]^{-s}][(\eta _2 -
 \overline {\zeta})/2i]^{-s}]\ \exp(-\epsilon z)
\exp (-\epsilon \overline {\zeta})$$
belongs to $H_r\overline {\otimes}\
\overline H_r$. As $\epsilon \rightarrow 0$ the
function $(z , \overline {\zeta})
\rightarrow (1-z  \overline {\zeta})^{(s-r)}
\phi_{\eta _1,\eta _2 ,\epsilon}(z,\overline {\zeta})$
converges weakly in
 $H_r\overline {\otimes}\  \overline H_r$ to the function
on  $\Bbb H ^2$
$$ (z , \overline {\zeta})
\rightarrow [(z - \overline {\eta _1})/2i]^{-s}][(\eta _2 -
\overline {\zeta})/2i]^{-s}].$$
But this is the symbol of the 1 dimensional operator in $B(H_s)$ defined
(on $H_s$) by
$$f\rightarrow e^r_{\eta _1}\langle f, e^r_{\eta_2} \rangle .$$
This proves the weak density of $j_{s,r}(\Cal C_2(H_r)$ in $B(H_s)$ and
this completes the proof.

We now define for $A$ in $B(H_r)$, which  by looking at its symbol,
 may
be thinked of as a continuous matrix,
 the analogue of the $(l^{\infty },l^1)$
norm for finite matrices. This norm will
 make the inclusions $j_{s,r}, \ s\geq r$
continuous and it will induce a new,
 involutive, Banach algebra structure
on a subalgebra in  $B(H_r)$.
\proclaim {Definititon 2.5} Let
 $|d(z,\overline {\zeta})|$  be the function
on the hyperbolic distance between $z$ and $\zeta$ in  $\Bbb H $ defined
in 2.2. For $A$ in $B(H_r)$ with contravariant symbol
 $\hat A=\hat A (z,\overline {\zeta}),
\ z,\overline {\zeta}\in  \Bbb H$ define $||A||_{\lambda ,r}$ to be the
maximum of the following integrals:
$$\sup\limits_{z\in \Bbb H}c_r \int_{\Bbb H}|
A(z,\overline {\zeta})d(z,\overline {\zeta})^r|
d\nu _0 ({\zeta})\text{\ and}$$

$$\sup\limits_{\zeta \in \Bbb H}
c_r \int_{\Bbb H}|A(z,\overline {\zeta})d(z,\overline {\zeta})^r|
d\nu _0 (z).$$
Let $\hat B(H_r)$ be the vector space of all $A$ in $B(H_r)$
 for which the
quantity $||A||_{\lambda ,r}$ is finite. Clearly $||\ ||_{\lambda ,r}$
is a norm on $\hat B(H_r)$.
\endproclaim
We will now show that with respect to the
 norm $||\ ||_{\lambda ,r}$,
$\hat B(H_r)$ becomes an involutive Banach
 algebra. Moreover the image
of $\hat B(H_r)$ in $B(H_s)$ for $s\geq r$ is closed under the products
$*_s$ for $s\geq r$ and this will allow us to differentiate the product.

\bigskip

\proclaim {Lemma 2.6} Let $A,B$ be bounded linear operators in $B(H_r)$
with symbols $\hat A=\hat A (z,\overline {\zeta})$ and respectively
$\hat B=\hat B (z,\overline {\zeta}),\ z,\overline {\zeta}\in  \Bbb H$.
If $A,B$ belong to the algebra $\hat B(H_r)$ defined in
  in the preceding statement,
then for all $s\geq r$, the product
$$j_{s,r}(A)*_s j_{s,r}(B)$$
which we  denote  simply by $A*_s B$, belongs to $\hat B(H_r)$
(more precisely to $j_{s,r}(\hat B(H_r))$). Moreover
$$||A*_s B||_{\lambda ,r}\leq
 2^{s-r}(\frac {c_s}{c_r})||A||_{\lambda ,r}||B||_{\lambda ,r}.$$
\endproclaim

Proof. It is sufficient to show that the integral:
$$c_r c_s \int_{\Bbb H}
|\hat A(z,\overline {\eta})||(z - \overline {\eta })/2i|^{-s}
|\hat B(\eta, \overline {\zeta})||(\eta
  - \overline {\zeta})/2i|^{-s}d\nu _s (\eta)
|d(z,\overline {\eta})|^r d\nu _0 (\zeta)$$
(and a similar integral corresponding
 the other term in the definition of
$||\ ||_{\lambda ,r}$) are bounded
 by $2^{s-r}(\frac {c_s}{c_r})||A||_{\lambda ,r}||B||_{\lambda ,r}.$

Regrouping the factors in the integral
 it follows that we have to estimate:
$$c^2_r (\frac {c_s}{c_r}) \int \int_{\Bbb H}|
\hat A(z,\overline {\eta})||d(z,\overline {\eta})|^r
|\hat B(\eta,\overline {\zeta})||d(\eta,\overline {\zeta})|^r
 (M(z,\eta ,\zeta))^{s-r}
 d\nu _0 (\eta ,\zeta)$$
where $M$ is the positive valued function on
 $\Bbb H ^3$ defined by the formula
$$M(z,\eta ,\zeta)=
|\text{Im}\ \eta||(z - \overline {\zeta })/2i|
(z - \overline {\eta })/2i|^{-1}
|(\eta - \overline {\zeta })/2i|^{-1}$$
We will show in the next lemma that $M(z,\eta ,\zeta)\leq 2$ for all
$z,\eta ,\zeta$ in $\Bbb H$. Hence the integral
 is bounded by
$$2^{s-r}(\frac {c_s}{c_r})c^2_r \int \int_{\Bbb H}
|\hat A(z,\overline {\eta})||d(z,\overline {\eta})|^r
|\hat B(\eta,\overline {\zeta})||
d(\eta,\overline {\zeta})|^r d\nu _0 (\eta ,\zeta) \leq$$
$$\leq 2^{s-r}(\frac {c_s}{c_r})c_r
\int_{\Bbb H}|\hat A(z,\overline {\eta})||d(z,\overline {\eta})|^r
(c_r\int_{\Bbb H}|\hat B(\eta,
\overline {\zeta})||d(\eta,\overline {\zeta})|^r d\nu _0
 (\zeta)) d\nu _0 (\eta )\leq$$
$$\leq 2^{s-r}(\frac {c_s}{c_r})
||B||_{\lambda ,r}c_r \int_{\Bbb H}|\hat A(z,\overline {\eta})||
d(z,\overline {\eta})|^r
d\nu _0 (\eta )\leq 2^{s-r}(\frac {c_s}{c_r})
||A||_{\lambda ,r} ||B||_{\lambda ,r}.$$
This completes the proof of the lemma, subject to checking the following
estimate.
\proclaim {Lemma} 2.6 Let $M=M(z,\eta ,\zeta)$ be
 the function $\Bbb H ^3$
defined by the formula
$$M(z,\eta ,\zeta)=\frac {(\text{Im}\ \eta)|
(z - \overline {\zeta })/2i|}{|(z - \overline {\eta })/2i|
|(\eta - \overline {\zeta })/2i|} ,\ z,\eta ,\zeta \in \Bbb H.$$
Then $0\leq M(z,\eta ,\zeta)\leq 2$ for all
 $z,\eta ,\zeta$ in $\Bbb H$.
\endproclaim

Proof. We use the formula
$$\rho (z,\zeta)=
\frac {\text{Im} \ z\ \text{Im}\ \zeta}{|(z-\zeta)/2i|^2}=
\cosh^{-1}(D_h (z,\zeta)),\ z\zeta \in  \Bbb H,$$
(recall that
 we used the notation $D_h (z,\zeta)$ for
the hyperbolic distance between the points
$z,\zeta$ in $\Bbb H$).
Then
$$M(z,\eta ,\zeta)=
\rho (z,\eta) \rho (\eta,\zeta) \rho ^{-1} (z,\zeta),\ z,\eta ,\zeta$$
and
 hence $M$ is invariant under the diagonal action of the group
 $PSL_2(\Bbb R)$
on $\Bbb H^3$.
 Consequently $M(z,\eta ,\zeta)$ is an invariant function
in the hyperbolic geometry of $\Bbb H$. Hence we may replace
 $M$, to compute
its supremumum, by the corresponding form of $M$ in the
 terminology
of the hyperbolic geometry
on the unit disk.

Hence we may assume that $$M(z,\eta ,\zeta)=
\frac {(1-|\eta |^2)|1-z \overline \zeta|}
{|1-z \overline \eta ||1-\eta \overline \zeta|},
\ z,\eta ,\zeta \in \Bbb D.$$
Since $M$ is invariant under the diagonal action
 of the group $SU(1,1)$
and since this group acts transitively on
 $\Bbb D$ it is sufficient to estimate
this when $\eta =0$.

Thus we want un upper estimate for
$$M(z, 0,\zeta)=
 |1-z \overline \zeta| \text{ for } z,\ \zeta \in \Bbb D.$$
The number 2 is clearly
 the lowest upper bound. This completes the proof.
\bigskip

In the following lemma we
 prove that $\hat B(H_r)$ is continuosly embedded
$B(H_r)$. In fact we will
 show a stronger statement; any kernel function
$k=k(z,\zeta)$ on  $\Bbb H^2$ which has the property
 that the norm $||k||_{\lambda ,r}$
is finite defines an integral operator
 on $L^2 (\Bbb H ,\nu _r)$ of uniform
norm less than $||k||_{\lambda ,r}$.

This will be an easy consequence of the operator
 interpolation technique
used by Berezin in [Be 2].

\proclaim {Proposition 2.7} Let $k=k(z,\zeta)$ be
 a function on $\Bbb H^2$
such that the maximum of the following expressions
 is finite:
$$\sup\limits_{z}c_r \int_{\Bbb H}|k(z,\zeta)|\cdot
|d(z, \overline \zeta|^r d\nu _0 (\zeta)$$
$$\sup\limits_{\zeta}c_r \int_{\Bbb H}|k(z,\zeta)|\cdot
 |d(z, \overline \zeta|^r d\nu _0 (\zeta).$$
Denote the maximum of the above two expressions
 by $||k||_{\lambda ,r}$.
 Then the integral
operator $k$ on  $L^2 (\Bbb D ,\nu _r)$ defined
 by the formula:
$$(kf)(z)=c_r \int_{\Bbb H}\frac {k(z,\zeta)}
{|(z - \overline {\zeta })/2i|^r}f(\zeta)d\nu _r (\zeta)
,\ z\in \Bbb H,\ f\in L^2 (\Bbb D ,\nu _r) $$
is continuous and its norm is bounded by $||k||_{\lambda ,r}$.
In particular
for all $B$ in $B(H_r)$ we have that
$$||B||_{\infty ,r}\leq ||B||_{\lambda ,r}.$$
\def \psl{PSL(2,\Bbb R)}

\def\b{B(H_r)}
\def\nk{||k||_{\lambda ,r}}
\def\d{\Bbb D}
\def\h{\Bbb H}

\def\i{\int_{\Bbb H}}
\def\c{c_r}
\def\l{L^2 (\Bbb D ,\nu _r)}
\def\s{\sup\limits_{\zeta \in \Bbb H}}
\def\a{\hat A(z,\zeta)}
\def\ra{\rightarrow}

\def\dz{|d(z,\zeta)|^r}
\def\dn{d\nu _0 (\zeta)}
\def\k{|k(z,\zeta)|}
\def\bh{\hat B(H_r)}
\def \e{\epsilon}
\def\n{||\ ||_{\lambda ,r}}
\def\na{||A||_{\infty ,r}}
\def\an{\langle e^{\frac{s-r}{2}}_z , e^{\frac{s-r}{2}}_z \rangle}
\endproclaim
Proof. Proving that $k$ is bounded on $\l$ is equivalent,
 by using the isometry
$$g\ra g(z) (\text{Im}\ z)^{r/2}$$
from $\l$ onto $L^2 (\Bbb D ,\nu _0)$, proving that the
 operator $k_0$
on $L^2 (\Bbb D ,\nu _0)$ defined by
$$(k( g))(z)=
\i \k \dz g(\zeta)\dn ,\ z\in \h ,\ g\in L^2 (\Bbb H ,\nu _0)$$
is a bounded linear operator. Moreover
 the uniform norms of $k$ and $k_0$
coincide.

Our hypothesis is that
$$\sup\limits_{z \in \Bbb H} \c \i \k \dz \dn \leq \nk$$
$$\sup\limits_{\zeta \in \Bbb H} \c \i \k \dz d\nu _0 (z) \leq \nk.$$
The interpolation technique used in [Be], Th 2.4, pp1131,
shows that
$k_0$ is bounded, of uniform norm less then $\nk$.

We now summarize the properties that we obtained so far
 for the algebras
$\bh$ and $\b$. The remarcable point about this is that
 the algebra $\bh$
is closed under all the products $(*_s)$ for $s\geq r$
 and moreover with
respect to the norm $\n$, $\bh$ is a banachique algebra
 with respect to the product
$(*_s)$  and the norm $\n$.
 (That is $||A*_s B||_{\lambda ,r}\leq \text{const }
 ||A||_{\lambda ,r}||B||_{\lambda ,r}$
for all $A,B$ in$\bh$ with a constant depending only on $s\geq r$).

\proclaim {Corrolary 2.8}
$\bh$ is an involutive Banach algebra with respect to the
 norm $\n$ from
Definition 2.5 (and the product $(*_s)$). Moreover, with
 respect to this
norm on $\bh$, and the uniform norm on $\b$, the inclusion of $\bh$ into
$\b$ is continuous with weakly dense image.

Moreover for all $s\geq r$, $\bh$ is mapped by $j_{s,r}$
 continouosly into
$\hat B(H_s)$ and the image of $\bh$ in $\hat B(H_s)$ is
 closed under the
multiplication $*_s$ and we have
$$||A*_s B||_{\lambda ,r}
\leq (\frac {c_s}{c_r})2^{s-r} ||A||_{\lambda ,r}||B||_{\lambda ,r}$$
for all $A,B$ in $\bh$, $s\geq r$.

Finnally  $\bh$ is weakly dense in $B(H_s)$ for
 all $s\geq r$.
\endproclaim

Proof. The last statement will follow from the next lemma in
 which we show
that  $B(H_{r-2-\epsilon})$ is contained in $\bh$ for all
 $\epsilon >0$.
One uses here also the fact that $j_{s,r}(\b)$ is weakly
 dense in $B(H_s)$ for all
$s\geq l$. The only assertation that needs to be proved
 is that $\bh$
embeds continuously into $\bh$ for $s\geq r$.This follows from the fact
that $d$ takes only subunitary values (see definition 2.1).

It is easy to see that the uniform norm on $\b$ is not equivalent
 to the
norm $\n$, i.e. $\bh$, as expected, is a strictly smaller
 algebra than
$\b$. (This may be verified by looking at the values of
 the two norms $\n$
and $||\ ||_{\infty ,r}$ on one dimensional projections in $\b$.

However the two norms are equivalent on the image of
 $B(H_{r-2-\epsilon})_+$
in $\bh$ for all $\epsilon >0$. We will sketch a
 proof of this (simple)
fact although we are not going to use it.

First we note a corrollary of Lemma 2.
\proclaim {Corrollary 2.8} If $\ 1<r<s-2-\e$, and $\e$
 is strictly positive
then $js,r$ maps $\b$ into $B(H_s)$ and
$$||js,r(A)||_{\lambda ,s}\leq (\text{const})
(\frac {s-r}{2}-1)^{-1}\na$$
for all $A$ in $\b$.
\endproclaim
Proof. If $A$ belongs to $\b$ then we have
 proved that
$$\sup\limits_{z,\zeta \in \Bbb H^2}|A(z,\zeta)|\dz \leq 4\na$$
for all $A$ in $\b$. Hence
$$\c \int  |A(z,\zeta)||d(z,\zeta)|^s
\dn \leq 4||A||_{\infty}\ \frac {c_r}{c_s}
c_s \int  |d(z,\zeta)|^{s-r}\dn=$$
$$ =4\na \i [\frac {(\text{Im} \ z)(\text{Im} \ \zeta)}
{|(z - \overline {\zeta })/2i|}]^{s-r}\dn=$$
$$=\text{const } \na (\text{Im}\ z)^{s-r} \an =
\text{const } \na (\frac {s-r}{2}-1)^{-1}.$$
\proclaim {Corrollary} 2.9 Let $r<s-2-\e$ with
 $\e$ strictly positive.
Then for all $A$ in $B(H_r)_+$, which is identified,
 by $js,r$, with
an element in $\b$, we have
$$ ||A||_{\lambda ,s}
\leq \text{const}(\frac {s-r}{2}-1)^{-1}\na .$$
\endproclaim

Proof. We proved as in the proof of the above corrollary.
 Since $A$ is positive
in $\b$ we have in addition that
$$\sup\limits_{z,\zeta \in \Bbb H}
 |A(z,\overline \zeta)||d(z,\zeta)|^r \leq
\sup\limits_{z \in \Bbb H}
 |A(z,\overline z)|\leq ||A||_{\infty, s}$$
(we use here proposition 1.5.d).
The argument then  follows as above.
\bigskip

Finnaly we mention the following interesting
 behaviour of the weak topology
on the unit ball in $\bh$ which is a corrollary
 of the preceding discussion
(although we are not going to make use of this fact, either.)

\proclaim {Lemma 2.10} The unit ball of $\bh$
 (with respect to the norm
$\n$) is compact with respect to the weak operator topology on $\b$.
\endproclaim

Proof. We already know that $\hat B(H_r)_1 \subseteq  B(H_r)_c$ for
some $c>0$ (by Corrollary 2.9). Let $B_n$ be any sequence $\hat B(H_r)_1$.
by the preceding assertation we may assume that $B_n$
 converges weakly
to $B$ for some $B$ in $B(H_r)_c$.

Let $$f_n(z,\zeta)=|B_n(z,\zeta)|
|d(z,\zeta)|^r, \ z,\zeta\in \h,$$
$$f(z,\zeta)=|B(z,\zeta)||d(z,\zeta)|^r,
 \ z,\zeta\in \h.$$
We know that the positive functions $f_n(z,\zeta)$,\ $f(z,\zeta)$ on
$\Bbb H^2$ are uniformely bounded by $C$ on $\Bbb H^2$.
Moreover $f_n$ converges punctually to $f$.

We want to take to the limit the following inequality that holds
for all the functions $f_n$:
$$\i f_n(z,\zeta)\dn \leq 1$$
for all $z$ in  $\h$.

For fixed $z$ in $\h$, we use the Lebesgue dominated
 convergence theorem
(as all functions are uniformely bounded by $c$) to
 deduce that for any
subset $E$ of $\h$ of finite $\nu_0$-measure we have

$ \int_{E}f(z,\zeta)\dn \leq 1$ for all $z$ in $\b$.

As $f$ is positive we may now use Fatou lemma to deduce that
$$\sup\limits_{z \in \Bbb H} \i f(z,\zeta)\dn \leq 1 .$$
A similar computation holds for $\s \i f(z,\zeta)\dn $.
 Hence $B$ belongs
to $\hat B(H_r)_1$.

\bigskip
\centerline{Part III}
\bigskip

\centerline{ The Berezin quantization for quotient
 space $\Bbb H /\Bbb {\Gamma}$}
\bigskip

\bigskip
\def \psl{PSL(2,\Bbb R)}

\def\b{B(H_r)}
\def\nk{||k||_{\lambda ,r}}
\def\G{\Bbb {\Gamma}}
\def\h{\Bbb H}

\def\i{\int_ G}
\def\c{c_r}
\def\l{L^2 (G)}
\def\s{\sup\limits_{\zeta \in \Bbb H}}
\def\a{\hat A(z,\zeta)}
\def\ra{\rightarrow}

\def\dz{|d(z,\zeta)|^r}
\def\dn{d\nu _0 (\zeta)}
\def\k{k(z,\overline \zeta)}
\def\bh{\hat B(H_r)}
\def\g{\gamma}
\def\hg{\Bbb H /\Bbb {\Gamma}}
\def \e{\epsilon}
\def\n{||\ ||_{\lambda ,r}}
\def\na{||A||_{\infty ,r}}
\def\an{\langle e^{\frac{s-r}{2}}_z , e^{\frac{s-r}{2}}_z \rangle}
\def\z{\zeta}
\def\st{*_s}
\def\fm{|\langle\pi (g)\zeta,\eta \rangle _H|^2}
\def \pt {\tilde{\pi}^0 _r}
\def\lbg{\Cal L(\G)\overline{\otimes} B(L^2(F,\nu_r))}
\def\lfg{l^2(\G)\overline{\otimes} L^2(F,\nu_r)}
\def\lbf{B(L^2(F,\nu_r))}
\def\lf{L^2(F,\nu_r)}
\def\calg{\Cal L(\G)}
\def\ovt{\overline{\otimes}}
\def\tr{\text{tr}}
\def\lh{L^2(\Bbb H,\nu_r)}
\def\lbh{B(\lh)}
\def\hr{H^2(\Bbb H,\nu_r)}
\def\eps{\epsilon}
\def\epst{\tilde{\eps}}
\def\pr{\tilde{\pi}_r}

In this paragraph we are analysing the $\G$-invariant form of the Berezin
quantization. Let$(\ast_s)_{s \in (a,b)}$ be the product for
Berezin symbols which was derived from the composition rule
of linear operators acting on the Hilbert spaces\break
$H^2(\Bbb H, (\text{Im}\ z)^{-2+s}\text{d} z\text{d} \overline z)$.
  Let $\G$ be a fuchsian subgroup of $\psl$ which we allow
to be of finite or infinite covolume. Since the Berezin quantization is
constructed such that $\psl$ acts as a group of symmetries, it follows
that the symbols $k=\k ,\ z,\zeta \in \Bbb H$, that are invariant  under
the diagonal action of the group $\G$ (i.e. $k(\g z,\overline {\g \z})=\k ,\
\g \in \G, z,\z \in \h$) are closed under any of the products $*_s$.

 By analogy with [Be] (see also [KL]) to obtain a deformation
quantization for $\hg$, we  let the algebras in the corresponding deformation
to be the vector space of $\G$-invariant symbols and we let the
 multiplication be defined
by the products $*_s$. For the integrals entering the
formula of $\ast_s$ to be convergent we impose  suitable
   condition on the growth of the
symbols.

This algebras are  then identified with the commutant
 of the (projective) representation
of  $\G$ on $H_r$ that  is obtained by restriction of
  the projective representation
$\pi_r$ of $\psl$  to $\G$.

Using a generalization of a theorem by [AS], [Co], [GHJ], in the form
stated in the monograph [GHJ] we will prove that the
 algebras in the deformation
quantization are type $II_1$ factors or respectively
 properly semifinite algebras (corresponding
 the case when $\G$ has   finite or ,respectively, infinite covolume).

The above mentioned generalization of the theorem by [A.S.], [Co], [GHJ]
to projective reprezentations of $\psl$, is used to determine the
isomorphism class of the algebras in the deformation: they are all stably
isomorphic to the von Neumann  algebra $\Cal L (\Gamma)$ associated
to $\G$. Moreover, the
 dimension of the Hilbert space $H_r$ as a left module over the
corresponding  algebra
in the deformation, tends to zero when the deformation parameter $r=1/h$
tends to infinity.

Finnally, the Berezin  formula computing the trace for operators on $\b$
as the integral of the restriction of the symbol on $\h$, is now replaced
by an integral over a fundamental domain. This generalizes a
 rezult in [GHJ],
(formula 3.3.e, see also the manuscript notes [Jo]).
In particular we will show that the $\G-$invariant Berezin are a  natural
generalization of  the notion of automorphic forms for the group $\G$.

The following definition extends the formal dimension ([Go]) of
 a representation
with square integrable coefficients to more general projective
 representations.
\proclaim {Definition 3.1} Let $G$ be a unimodular locally compact group
with Haar measure $\text{d} g$. Let $\pi :G\ra B(H)$ be a projective,
 unitary representation of
$G$.
Assume that $\pi$ is (topologically) irreducible and also assume that $\pi$
has square integrable coeficients: i.e. there exists at least one nonzero
$\eta$ in $H$ so that
$$\i |\langle\pi (g)\eta,\eta \rangle _H|^2 dg =d^{-1}_{\pi} ||\eta||^2,$$
for some stricty positive number $d_{\pi}$. Then $d_{\pi}$ is independent
on the choice of $\eta$. Moreover,  the following equality holds
true  (the Schurr orthogonality relations):
$$\i |\langle\pi (g)\zeta,\eta \rangle _H|^2 dg =d^{-1}_{\pi}||\zeta ||^2
 ||\eta||^2 \text{  for all  }\zeta ,\eta \text{ in } H. $$
\endproclaim
Before proving the lemma, we recall the following folklore lemma:
\proclaim {Lemma} Let $G$ be an unimodular localy compact group. Let
$\{c(g,h)\}_{g,h \in G}$ be a familly of modulus 1 complex numbers
 defining an element
in $H^2(G,\Bbb T)$, that is
$$c(g,hk)c(h,k)=c(gh,k)c(g,h), \ g,h,k \text{ in } G.$$
 For all $h$ in $G$ define an unitary operator $R_h$ on $L^2(G)$
  by the following
 relation:
$$(R_h f)(g)=c(g,h)^{-1}f(gh),\ f\in L^2(G),\ g,h \text{ in } G.$$
Then $(R_h)_{h\in G}$ is a projective, unitary representation of $G$ on
$L^2(G)$ with
cocycle
$\{c(g,h)\}_{g,h\in H}$. Thus
$$R_{hk}= c(h,k) R_h R_k,\ k,h \text{ in } G.$$
\endproclaim
Proof. We have $(R_{hk}f)(g)= c(g,hk)^{-1}f(ghk)$.
On the other hand
$$( R_h (R_k f))(g)=c(g,h)^{-1}(R_k f)(gh)=c(g,h)^{-1}c(gh,k)^{-1}f(ghk).$$
Hence
$$ (R_{hk}f)(g)= c(g,hk)^{-1}c(g,h)c(gh,k) ( R_h (R_k f))(g)=$$
$$=c(g,h)( R_h (R_k f))(g), \
\text{ for all }f \in L^2(G),g,h\in G.$$
This completes the proof.
\bigskip

Let $\pi$ be a projective, unitary
 representation of $G$.
Let $ c(g,h)$ be the complex number of modulus $1$ defined by the relation
$$\pi (gh)=c(g,h) \pi (g)\pi (h),\ \text{ for }g,h\in G$$
Because of the above lemma, we may  think to the  projective,
 unitary representation $\pi$
as subrepresentation of the "skewed", left regular representation
 of $G$ on $L^2(G)$, with cocycle $c$.
\bigskip

Proof of Proposition 3.1. We follow the lines contained in the proof of the
analogue result for the  discrete series of  representations of $G$,
 as it is found in the
exposition in [Ro], chapter 17.

Keeping $\eta$ fixed in $H$, we define the linear map $T$ from the vector
space
$W=\{\z \in H \vert \i\fm dg<\infty \}$ into  $L^2(G)$ by
$(T\z )(g)= \langle\pi (g)\zeta,\eta \rangle.$
Then $$(T\pi (h)\z )(g)=\langle\pi (g)\pi (h)\zeta,\eta \rangle=$$
$$c(g,h)^{-1}\langle\pi (g,h)\zeta,\eta \rangle =R_h (T\z )(g) \text
{ for\ all\ } g\in G.$$
Thus $T\pi (h)=R_h T$ and $T$ is an (eventually) unbounded, closed operator
with domain $W$. Moreover the operators  $ \pi (h),h\in G$ map $W$ into $W$.

The theorem 15.13 in [Ro] shows that $T$ is a multiple of an isometry
and in particular $W=H$.
In particular we find that
$$\i |\langle\pi (g)\zeta,\eta \rangle |^2 dg<\infty ,
\ \text{ for all }\z \text{ in }H. $$
Similar arguments involving the, cocycle perturbed, left,
 regular representation
of $G$ on  $L^2(G)$, will show that:
$$\i |\langle\pi (g)\zeta,\eta ' \rangle |^2 dg<\infty ,
 \text{ for all }\z ,\eta ' \text{ in }H $$
The rest of the argument is exactly as in the proof
 of Theorem 16.3 in [Ro]. This completes the proof.
\bigskip

For the projective representation $(\pi _r)_{r>1}$ Recall that
the discrete series of representation $(\pi _n)_{n\in \Bbb N}$ of
$\psl$  have formal dimension $\frac {n-1}{\pi}$. We will prove a similar
statement for the projective representation $(\pi _r)_{r>1}$. It will
 then follow
 that the projective, unitary representations
$(\pi _r)_{r>1}$ have similar properties to  square integrability
 with $d_{\pi _r }=\frac {r-1}{\pi}$

The explanation for this fact should be looked into the Plancharel formula
for the universal cover $\widetilde {\psl}$ of $\psl$.
L. Pukanzki [Pu] has shown that the reperesentations
$(\pi _r)_{r>1}$ are a summand  in the continuous series in the
 Plancherel formula for $\widetilde {\psl}$. The coefficient
 of the representation
$\pi_r$ in this formula is
$d_{\pi _r }=\frac {r-1}{\pi}$.
\def \psl{PSL(2,\Bbb R)}

The main result of this paragraph is the following:
\proclaim {Theorem 3.2}
Let $\G$ be a fuchsiangroup in $\psl$ of finite  or infinite
covolume.
 Let $r>1$ and let $\Cal A_r$ and $\hat{\Cal A_r}$  be the
vector space of all symbols $k=\k, \ z,\overline \z \in\h$ analytic
in $z$ and  antianalytic in $\z$, that are $\G$ invariant (i.e.
$\k =k(\gamma z,\overline{\gamma \z}), \ \gamma \in \Gamma,
\ z,\z \in \h$) and which  correspond to bounded operators
in $\b $ or $\bh$, respectively. Then the vector spaces
  $\Cal A_r$ and $\hat{\Cal A_r}=\hat B(H_r)\cap  \Cal(A_r)$
are closed under the product $*_r$.

Let $F$ be a fundamental domain for the action of $\G$ in
$\h$ (see [Leh]).
Then

i) $\Cal A_r$ is a type $II_1$ factor or an infinite semifinite
von Neumann corresponding respectively to the case when $\G$
is of finite or infinite covolume. Moreover $\Cal A_r$ is
isomorphic to the commutant $\{\pi _r(\G)\}'\subseteq \b$
of the image of  $\G$ in $\b$ via the  representation $\pi _r$.

ii) In the finite covolume case the trace on $\Cal A_r$
  is given by the formula
$$\tau (k)=\frac{1}{\text{area}(F)} \int_F k(z,\overline z)
d\nu _0 (z) ,\ k\in \Cal A_r.  $$
  When $\G$ has infinite covolume in $\psl$,
 one  defines
a semifinite, faithfull, normal trace on $\Cal A_r$  by the formula
$$\tau (k)= \int_F k(z, z)d\nu _0 (z).$$

iii) Assume that the group cocycle in the second cohomology
 group of
$\psl$, $H^2(PSL^2(\Bbb R),\Bbb T)$ associated with
 the projective
representation of $\psl$ on $H_r$,
 vanishes by restriction to  $H^2(\G, \Bbb T)$.
Then $\Cal A_r$ is isomorphic to $e(\Cal L(\G)\otimes \b)e$
where $e$ is any projection in $L(\G)\otimes \b$ of trace $\tau(e)$=
covol\ $(\G)d_{\pi_r}=\frac{r-1}{\pi}(\text{covol}(\G))$.
\endproclaim
\proclaim {Remark} Note that the condition in (iii) is always
satisfied when $\G$ is not cocompact (see [Pa2]).
\endproclaim

The proof of the theorem will be splitted into more lemmas.
In a slightly different form the first two lemmas
  may be found in the book of [Ro]. The first lemma
shows that the representations $\pi_r$ have the property that
they move (modulo scalars) $\pi_r(g)$ maps the
vector $e^r_z$ into $e^r_{g^{-1}z}$
for all $z$ in $\h$.

Recall from paragraph 1, that the representation $\pi_r$ were
defined by the formula
$$(\pi_r(g)f)(z)=(j(g,z))^{-r}f(g^{-1}z),$$
for all $g$ in $G$, $f$ in $H^2(\h ,\nu_r)$ and $z$ in $\h$.
They are projective representations  because the factor
$$(j(g,z))^r=\exp(r\ln(j(g,z))=\exp(r\ln(cz+d))$$
for $g=\pmatrix a&b\\c&d \endpmatrix$ in $\psl$, $z$ in $\h$,
involves the  choice of a branch for $\ln (cz+d)$.
\proclaim{Lemma 3.3}  ([Ro]). Let $r>1$ and let $e^r_z\in H^2(\h, d\nu_r)$
be the evaluation vector at $z$, i.e. $\langle f, e^r_z  \rangle_r = f(z)$
for all $f$ in $H_r= H^2(\h, d\nu_r)$. Then there exists $\Theta_{g}\in \Bbb C$
of modulus 1 so that
$$\pi_r(g)e^r_z=\Theta_{g}(j(g^{-1},z))^r e^r_{g^{-1}z}$$ for all $g$
in $G$, $z\in \h$.
\endproclaim
Proof. Indeed for any $f$ in $H_r$ we have
$$\langle f,\pi_{r}(g), e^r_z  \rangle_{H_r}=
\langle \pi_r(g^{-1}) f,\pi_r(g^{-1})\pi_r(g) e^r_z  \rangle_{H_r}$$
by the unitary of $\pi_r$. Since  $\pi_r$ is a projective representation
of $G$, there exists a scalar $\overline {\Theta_{g}}$
of modulus 1 so that $\pi_r(g^{-1})\pi_r(g)= \overline {\Theta_{g}}\  Id_{\b}$.
Hence
$$\langle f,\pi_{r}(g), e^r_z  \rangle=
\Theta_{g}\langle \pi_r(g^{-1}) f, e^r_z  \rangle
=\Theta_{g}( \pi_r(g^{-1}) f) (z)=$$
$$\Theta_{g}(j(g^{-1},z))^{-r}f(g^{-1}z)=\Theta_{g}(j(g^{-1},z))^{r}
\langle f,e^r_{g^{-1}z}\rangle_r .$$
It is now easy to prove that an operator $A$ in $\b$ commutes
with $\pi_r(\G)$ if and only if its symbol has the property
that its symbols is $\G$ invariant. This will complete the
proof of i) by showing that the $\G$ invariant symbols correspond
to elements in the commutant of $\pi_r(\G)$ in $\b$, and hence
they are clearly closed under multiplication.
\proclaim {Lemma}3.4 An operator $A$ in $\b$ of Berezin symbol
$\hat A=\a , z,\z \text{ in  }\h$ commutes  with $\pi_r(\G)$
if and only if $\hat A$ is $\G$ invariant under the diagonal
action, that is
$$\hat A(\g z,\overline {\g \z})=
\a \ \text{ for all }z,\z \text{ in } \h, \g \in \G.$$
\endproclaim
Proof. Assume that $A$ commutes with $\pi_r(\G)$.
Then for all $\g$ in $\G$,
$$ A(\g ^{-1}z,\g ^{-1} \z)=
\frac{\langle A e^r_{\g^{-1}\z},e^r_{\g^{-1}z}\rangle}
{\langle e^r_{\g^{-1}\z},e^r_{\g^{-1}z} \rangle}=$$
$$ \frac {\langle A (\Theta_{\g^{-1}}(j(\g,\z))^{-r})\pi_r(\g^{-1})e^r_{\z},
(\Theta_{\g^{-1}}(j(\g,z))^{-r})\pi_r(\g^{-1})e^r_{z}\rangle}
{\langle  (\Theta_{\g^{-1}}(j(\g,\z))^{-r})\pi_r(\g^{-1})e^r_{\z},
(\Theta_{\g^{-1}}(j(\g,z))^{-r})\pi_r(\g^{-1})e^r_{z}\rangle}$$
where we used the preceding lemma.
 The scalars are canceling
themselves and hence we get
$$\hat A(\g ^{-1}z,\overline {\g ^{-1} \z})=
\frac{\langle A\pi_r(\g^{-1})e^r_{\z},\pi_r(\g^{-1})e^r_{z}\rangle}
 {\langle \pi_r(\g^{-1})e^r_{\z},\pi_r(\g^{-1})e^r_{z}\rangle}.$$
Since $\pi_r(\g^{-1})$ is unitary and $\pi_r(\g^{-1})$ commutes
with $A$, we get $A(z,\overline \z)$. The converse follows this
lines too. This completes the proof of i). (We also use here
the fact that $\bh$ is closed under multiplication and hence
so is $\Cal {\hat A}_r = \bh \cap \Cal A_r).$

Proof of ii). To prove ii) we will first check that the formula
for $\tau$ defines indeed a trace. Note that by point i)
$k(\g ^{-1}z,\overline {\g ^{-1} z})=k(z,\overline z)$ for
all $z$ in $\h,\ \g$ in $\G$ and hence the integral for
$\tau$ doesn't depend on the choice of $F$.

 To check that
$\tau$ is a trace it is sufficient to check that $\tau (k^* *_rk)=
\tau (k *_r k^*)$
for every $\G-$ equivariant kernel $k=k(z,\overline\z)$ on $\h ^2$ which
is so that  the integrals in the formulae for
$\tau (k^* *_rk)$ and
$\tau (k *_r k^*)$  are absolutely convergent.

Assume that $k$ is as above and that $k$ is so that the integrals
involved in $\tau (k^* *_rk)$ or $\tau (k *_r k^*)$ are convergent.

We have: $k^*(z,\overline {\z})=\overline{k(\z,\overline z)}$ for
$z,\z$ in $\h$ and hence
$$(k^* *_rk)(z,\overline {\z})=
c_r[(z - \overline {\zeta })/2i]^r \int_{\h}\frac{k^*(z,\eta)}
{[(z - \overline {\eta })/2i]^r}
\frac{k(\eta,\overline {\z})}{[(\eta - \overline {\zeta })/2i]^r}
d\nu _r (\eta)$$
so that
$$(k^* *_rk)(z,\overline {z})=
c_r \int_{\h}k^*(z,\overline {\eta})k(\eta,\overline {\z})
|d(z,\overline {\eta})|^{2r}d\nu _0 (\eta)$$
and hence
$$\tau (k^* *_rk)=(\text{const})c_r\int_F(\int_{\h}
|k(\eta , \overline {z})|^2|d(z,\overline {\eta})|^{2r}
\text{d}\nu _0 (\eta))\text{d}\nu _0 (z).$$
 Similarly
$$\tau (k *_r k^*)=(\text{const})c_r\int_F(\int_{\h}
|k(\eta , \overline {z})|^2|d(z,\overline {\eta})|^{2r}
\text{d}\nu _0 (\eta))\text{d}\nu _0 (z).$$
The constant in front of the two integrals equals $1$ if $F$
has infinite covolume and equals  (area $F)^{-1}$ if $\G$
has finite covolume (the hyperbolic) area is computed here relative
to $\nu_0$).

By renaming the variables in the integral for $\tau (k^* *_rk)$
we get:
$$\tau (k^* *_rk)=(\text{const})c_r\int_F(\int_{\h}
|k(z , \overline {\eta})|^2|d(\eta, z)|^{2r}
d\nu _0 (z))d\nu _0 (\eta).$$
The second expression for $\tau (k^\ast *_rk)$ is different from
the one for $\tau (k *_r k^*)$ only in the choice of the
domain of integration (as $|d(z,\overline \eta)|=
|d(\overline \eta ,z)|$
for all $z,\eta$). But in both integrals we are integrating
over a fundamental domain for the diagonal action of $\G$ in
$\h^2$, while the integrand is $\G-$invariant under the same
action of $\G$. Since we are integrating  positive functions,
the two integrals must be equal.

We will now show that the functional $\tau$ is positive definite.
 Let
$k$ represent a positive operator in $\Cal A_r \subseteq \b$.
Then by Lemma 2.1., $k(z,\overline z)$ is positive $z\in \h$.
Hence if $\tau (k)=
(\text{const})\int_F k(z,\overline z)d\nu _0 (z)$
is $0$ then it follows that $k(z,\overline z)=0$ for all $z$
in $F$ and thus for all $z$ in $\h$ by the $\G$-invariance.

But $k=k(z,\overline \z)$ is analytic in $z$ and
antianalytic in $\z$.  Hence if  $k$ vanishes on the diagonal
$z=\z, \ z\in \h$, then it must be  identically zero. Thus we have
shown that $\tau$ is definite, i.e. that if $\tau (k)=0$ and
$k$ is the symbol of a positive element in $\Cal A_r$ then $k=0$.

In the case of a group $\G$ of finite covolume, $\tau$ is
well defined on $\Cal A_r$. Indeed let  $||k||_{\infty ,r}$
be the uniform norm of the operator on $H_r$ defined by $k$.
 We then have
$$|\tau (k)|=
|\frac {1}{\text{area}(F)}\int_F k(z,\overline z)d\nu _0 (z)|
\leq \sup\limits_{z \in F}| k(z,\overline z)|\leq ||k||_{\infty ,r}.$$
We   used for the last inequality Lemma 1.5.d.

In the case of a group $\G$ of infinite covolume we will have
to check in addition that there exists sufficiently many positive
elements in $\Cal A_r$ so that  the trace $\tau$ takes finite
value on them, and so that $1$ (the unit of $\Cal A_r$) is a
weak increasing limit of positive  elements  having finite trace.

To obtain such elements we let
 $k_f=k_f(z,\overline z), \ z,\overline z \in \h$
be the symbol of a Toeplitz operator $T^r_f=P_rM_fP_r$ on
$\b$. The symbol $f$ of the operator  $T^r_f$ is assumed to be a
positive function $f$ with compact support in the interior
$\overset\text{$\circ$}\to{F}$ of $F$ and then extended by $\G$
 invariance to the
whole upper half plane $\h$, (see also the next paragraph).

Since $f$ is $\G-$ equivariant, the operator $T^r_f$
commutes with $\pi_r(\G)$. Also $T^r_f$ is clearly a positive
operator. If we choose an increasing net of such functions $f$
that converges point wise  to the constant function $1$ then
it follows that $1_{\Cal A_r}$ is a weak (increasing) limit
of operators of the form $T^r_f$.

The functional $\tau$ evaluated on $T^r_f$ gives:
$$\int_F k_f(z,\overline z)d\nu _0 (z)=
\int_F c_r [(\text{Im} \ z)^r \int \frac
{f(a)}{|1-\overline z a|^{2r}}d\nu _0 (a)]d\nu _0 (z)=$$
$$=\text{const}\int_F f(a)[\int_{\h}|d(z,\overline a)|^{2r}
d\nu _r (a)]d\nu _0 (z).$$

But $\sup\limits_{z}\int_{\h}|d(z,\overline a)|^{2r}
d\nu _0 (a) <\infty$
and hence $\tau (T^r_f)$
 is finite.
This completes the proof of ii).

We now turn to the proof of iii). We will follow the lines for
the proof in the case of actual representations of $\psl$
which is contained in the monography [GHJ]. The only difference
is  that  the representation $\pi _r$is not
 a subrepresentation in the left regular representation
of $G$.Instead, for arbitrary real $r$,
one considers $\pi _r$ is a subrepresentation of the (projective)
unitary representation $\tilde{\pi _r}$ of $G$ into $B(L^2(\h,\nu _r))$ and
$\pi _r|_\G$ is a subrepresentation in the regular representation
of $\G$ into $B(l^2(\G))$.
\proclaim {Proposition} 3.5  Let $\tilde{\pi _r}:G=\psl\rightarrow
B(L^2(\h,\nu _r))$ be the (projective) unitary representation
of $G$ onto $ L^2(\h,\nu _r)$ defined by the same formula
as $\pi _r$:
$$(\tilde{\pi _r}(g)f)(z)=
(j(g,z))^{-r}f(g^{-1}z), \ f\in L^2(\h,\nu _r),
\ g\in G,\ z\in \h.$$
Assume that there exists complex numbers $c(\g)$ of modulus $1$
so that \break
 $\g \rightarrow c(\g) \tilde{\pi _r}(\g)=\tilde{\pi}^0 _r(\g)$
is a unitary representation of $\G$
on $L^2(\h,\nu _r)$. Let $F$ be a fundamental domain for $\G$
in $\h$.
Let $(e_{\g})_{\g \in \G}$ be the canonical orthonormal basis
for $l^2(\G)$ and let
 $V_r:l^2(\G)\otimes  L^2(F,\nu _r)\rightarrow L^2(\h,\nu _r)$
be defined by the formula
$$ V_r(\sum\limits_{\g}e_{\g}\otimes f_{\g})=
\sum\limits_{\g}\tilde{\pi}^0 _r(\g)
(\g^{-1})(f_{\g})$$
for all elements $\sum\limits_{\g}e_{\g}\otimes f_{\g}$ in
$l^2(\G)\otimes  L^2(F,\nu _r)$. Note that the functions
$f_{\g}$ are identified with elements in $L^2(\h,\nu _r)$ by
defining them to be zero outside $F$.

Let $R_{\g}:\G \rightarrow B(l^2(\G))$ be the right regular
representation of $\G$ (i.e. $R_{\g}e_{h}=e_{h}\g^{-1},$
 for all $
h,\g \in \G$).

Then $V_r$ is an unitary and
$$V^*_r\tilde{\pi _r}(\g) V_r=R_{\g}\overline\otimes
 Id_{B( L^2(F,d\nu _r)},
\ \text{ for all }\g\ \text{ in }  \G.$$
\endproclaim
Proof. We construct first an inverse $U_r$ for $V_r$ which
is defined on $L^2(\Bbb H,\nu _r)$ with values
 into $\ l^2(\G)\otimes  L^2(F,\nu _r)$
by
$$(3.1) \ \ U_r f=\sum\limits_{\g}e_{\g}\otimes f_{\g}, \
f_{\g}=\tilde{\pi}^0 _r(\g)(f\chi_{\g^{-1}F}),
\ \text{ for all }\ f\ \text{ in }  L^2(\Bbb H,\nu _r).$$
We denote by $\chi_{A}$ the characteristic function of $A$
in $\h$. Clearly $ \tilde{\pi}^0 _r(\g)(f\chi_{\g^{-1}F})$
has its support in $F$ for all $\g$ in $\G$ as
$$ \tilde{\pi}^0 _r(\g)(f\chi_{\g^{-1}F})(z)=
c(\g)(j(\g,z))^{-r}(f\chi_{\g^{-1}F})(\g^{-1}z)$$
is nonzero only for $z$ in $F$. Hence $U_r$ is well defined.
We check first that $V_rU_rf=f$ for all $f$ in $L^2(F,\nu _r)$.
Indeed if $f_{\g}$ are as formula (3.1) then
$$V_r(\sum\limits_{\g}e_{\g}\otimes f_{\g})=\sum\limits_{\g}
\tilde{\pi}^0 _r(\g^{-1})(f_{\g})=\sum\limits_{\g}
\tilde{\pi}^0 _r(\g)[\tilde{\pi}^0 _r(\g^{-1})(f\chi_{\g^{-1}F})]=$$
$$=\sum\limits_{\g}f\chi_{\g^{-1}F}=f$$
Clearly $U_r$ is surjective as $F$ is a fundamental domain
for $\G$ in $\h$ so that $\h$ is covered by translates of
$F$ by $\G$.
Moreover, $U_r$ is unitary because, given $f$ in $L^2(F,\nu _r)$
and letting $f_{\g}$ be defined by formula (3.1), then
$$||U_r f||^2=||\sum\limits_{\g}e_{\g}\otimes f_{\g}||^2_
{l^2(\G)\otimes  L^2(F,\nu _r)}=\sum\limits_{\g}
||\tilde{\pi}^0 _r(\g^{-1})(f\chi_{\g^{-1}F})||^2_
{ L^2(F,\nu _r)}=$$
$$=\sum\limits_{\g}||f\chi_{\g^{-1}F}||^2_
{ L^2(\h,\nu _r)}=||f||^2_
{ L^2(\h,\nu _r)}.$$
We have so far checked that $V_r$ is unitary from
 $l^2(\Gamma) \otimes L^2(F,\nu_r)$ onto
 $L^2(\Bbb H,\nu_r)$ and that its inverse is given
by the relation (3.1).

It remains to check that
$$U_r\pt(\gamma)V_r=R_{\gamma}
\overline{\otimes}Id_{B(L^2(F,\nu_r))},
 \text{for \ all} \gamma \in \Gamma.$$

Indeed for any element
 $\sum\limits_{\g}e_{\g}\otimes f_{\g}$
 in $l^2(\Gamma) \otimes L^2(F,\nu_r)$
and any $\sigma \in \Gamma$, we have
$$\lbrace U_r\pt(\sigma)V_r\rbrace
(\sum\limits_{\g}e_{\g}\otimes f_{\g})
=
U_r\pt(\sigma)
(\sum\limits_{\gamma}\pt(\gamma^{-1})f_{\gamma})=$$
$$=\sum\limits_{\gamma}U_r(\pt(\sigma\gamma^{-1})f_{\gamma})=
 U_r(\sum\limits_{\gamma}\pt(\gamma^{-1})
(f_{\gamma\sigma}\chi_F))=$$

$$\sum\limits_{\delta} e_{\delta}\otimes (\pt(\delta)
(\sum\limits_{\gamma}\pt(\gamma^{-1})
(f_{\gamma\sigma}\chi_F)\chi_{\delta^{-1}F}).$$

Note that in general for any $\alpha$ in $\G$
and $g$ in $L^2(F,\nu_r)$ the support of the
function $\pt(\alpha^{-1})(g\chi_F)$ is in
$\chi_{\alpha^{-1}F}.$
Hence, in the above summ, the only nonzero terms that may
occur, are those for which $\delta=\gamma$. Consequently
$$(U_r\pt(\sigma)V_r)
(\sum\limits_{\g}e_{\g}\otimes f_{\g})=$$
$$=\sum\limits_{\delta} e_{\delta}\otimes (\pt(\delta)
\lbrack\pt(\delta^{-1})(f_{\delta\sigma})\rbrack=$$
$$=\sum\limits_{\delta} e_{\delta}\otimes
f_{\delta\sigma}=
\sum\limits_{\delta} e_{\delta\sigma^{-1}}
\otimes f_{\delta}=$$
$$=R_{\sigma}\overline{\otimes} \text{Id}
(\sum\limits_{\delta} e_{\delta}\otimes
f_{\delta}).$$
This completes the proof of the proposition.

\proclaim{Corollary 3.6} Assume the hypothesis iii) from
Theorem 3.2. Then $\lbrace \tilde{\pi}_r(\G)\rbrace'$ is isomorphic
to $e(\Cal L(\G) \otimes B(H))e$ for a projection
$e$ in $\Cal L(\G) \otimes B(H)$.
\endproclaim

Proof. We use the notations from Proposition 3.5.
Note that $P_r$, the projection from $L^2(\Bbb H,\nu_r)$
on $H^2(\Bbb H,\nu_r)$ commutes with
$\tilde{\pi}_r(g)$, for all $g \in G$ (as $\tilde{\pi}_r(g)$
maps analytic functions into analytic functions).

Hence, by [St., Zs.], $\lbrace\pi_r(\G)\rbrace'$ is
isomorphic to $P_r\lbrace\tilde{\pi}_r(\G)\rbrace'P_r$.
The condition iii) in Theorem 3.2 shows
that proposition 3.5 applies and hence
$\lbrace \tilde{\pi}_r(\G)\rbrace'$ is isomorphic to
$\Cal L(\G)\overline{\otimes} B(L^2(F,\nu_r))$ (since in
general $\lbrace\Cal R(\Gamma)\rbrace' $ is
isomorphic to $ \cong \Cal L(\G)$).

We use the notations in Proposition 3.5 and Corollary 3.6
and also assume the condition iii) in Theorem 3.2. Let
$Q_r$= $V_r^{\ast}P_r V_r$ be the projection in
$B(l^2(\G)\otimes L^2(F,\nu_r)$ which corresponds to $P_r$ through the
identification of $L^2(\Bbb H,\nu_r)$ with
$l^2(\G)\otimes L^2(F,\nu_r)$, by the unitary
$U_r$.

We have just proved that $Q_r$ commutes with
$R_{\gamma} \overline{\otimes} \text{Id}_{B(L^2(F,\nu_r)}$
for all $\gamma \in \G$ (as $P_r$ commutes
with $\tilde{\pi}_r(\gamma)$ for all $\gamma \in \G$). Moreover
we have just proved that
\proclaim{Proposition 3.7} With the notations
in proposition 3.5, we have that
$\lbrace \tilde{\pi}_r(\G)\rbrace'$ is isomorphic
to $$Q_r(\Cal L(\G)\overline{\otimes} B(L^2(F,\nu_r)))Q_r
\subseteq B(l^2(\G)\overline{\otimes} B(L^2(F,\nu_r)).$$
Moreover $Q_r$ belongs to the algebra
$\Cal L(\G)\overline{\otimes} B(L^2(F,\nu_r)),$
which is the commutant of $$\lbrace R_{\gamma}
\overline{\otimes}\text{Id}_{B(L^2(F,\nu_r))},
\gamma \in \G\rbrace'.$$
\endproclaim

To determine the isomorphism class of the algebra
$\lbrace \tilde{\pi}_r(\G)\rbrace'$ we will choose a trace
$\tau_1$ on $\Cal L(\G)\overline{\otimes} B(L^2(F,\nu_r))$.
The trace $\tau_1$ will be normalized by the condition that
$\tau_1 $
takes the value 1 on minimal
 projections in   $B(L^2(F,\nu_r))$.
We will use the method explained in [GHJ].

\proclaim{Proposition 3.8} We use the notations
in Proposition 3.5. Let $\tau_1$ be the trace
$\tau_1$ on $\Cal L(\G)\overline{\otimes} B(L^2(F,\nu_r))$,
normalized by the condition that $\tau_1 $
takes the value 1 on minimal
 projections in   $B(L^2(F,\nu_r))$. Let
$\lbrace  \epst_n\rbrace_{n \in \Bbb N}$ be an
orthonormal basis in
$L^2(F,\nu_r)$. Let $\lbrace\delta_{\gamma}\rbrace_{\gamma
\in\G}$ be the canonical orthonormal  basis
for $l^2(\G)$ and let $e$ be the neutral
element in $\G$. Then
$$\tau_1(Q_r)= \sum\limits_{n \in \Bbb N}\vert\vert
Q_r(\delta_{e} \otimes \epst_n)\vert
\vert^2_{l^2(\G)\overline{\otimes} L^2(F,\nu_r)}.$$
\endproclaim
Proof.
Let $\tau_{\calg}$ be the canonical normalized trace
on $\calg$ and let
$$\text{tr}=\text{tr}_{B(L^2(F,\nu_r))}$$
be the semifinite trace on $B(L^2(F,\nu_r))$ taking value
1 on the one dimensional projections.
Then $\tau_1= \tau_{\calg}\overline{\otimes} \text{\tr}$.
Moreover for any element $x\ovt y$
in $\lbg$, with $y$ of trace class, we have:
$$\tau_1(x\ovt y)=
 \tau_{\calg}(x)\text{tr}(y)=$$
$$=\tau_{\calg}(x)\lbrack\sum\limits_{n \in \Bbb N}
\langle y\epst_n,\epst_n\rangle_{\lf}\rbrack=$$
$$=\sum\limits_{n \in \Bbb N}\langle
(x\ovt y) (\delta_e \ovt\epst_n),\delta_e \ovt\epst_n
\rangle_{\lfg}.$$
Since $Q_r$ is a weak limit of linear combinations
of elements $x\ovt y$ as above, this
concludes the proof of the lemma.

\proclaim{Remark} Let $\lbrace \eps_n\rbrace_{n \in \Bbb N}
\subseteq   \lh$ be the image under
$V_r$ of the orthonormal system $\lbrace
 \delta_e \ovt\epst_n\ \vert\  n \in \Bbb N\rbrace
\subseteq \lfg$.
As $$\lbrace
(R_{\gamma}\ovt\text{Id}_{\lbf})
 (\delta_e \ovt\epst_n)\ \vert\ \gamma \in \G,
 n \in \Bbb N\rbrace$$
is an orthonormal basis for $\lfg$ it follows that
$$\lbrace \tilde{\pi}_r(\gamma)\eps_n\
\vert\ \gamma \in \G,
 n \in \Bbb N\rbrace$$
is an orthonormal basis for $\lh$.
\endproclaim

The properties of the orthonormal
system $\lbrace \eps_n\rbrace_{n \in \Bbb N}$
in $\lh$ are summarized in the next proposition
(whose proof has allready been completed).

\proclaim{Proposition 3.9} We use the above notations.
Let $Q_r$ be the projection
$U_rP_rV_r$ in $\lbg$.
Then $\lbrace\pi_r(\G)\rbrace'$ is isomorphic to
$Q_r(\lbg)Q_r$.

 Let $\tau_1$ be the trace
$\tau_1$ on $\Cal L(\G)\overline{\otimes} B(L^2(F,\nu_r))$,
normalized by the condition that $\tau_1 $
takes the value 1 on minimal
 projections in   $B(L^2(F,\nu_r))$ and so that $\tau_1=
\tau_{\calg}$ if
restricted to $\calg\cong\calg\ovt \text{Id}_{\lbf}$.

Then there exists an orthonormal system
$\lbrace \eps_n\rbrace_{n \in \Bbb N}$ in $\lh$
so that

 i).$\ \lbrace \tilde{\pi}_r(\gamma)\eps_n\
\vert\ \gamma \in \G,
 n \in \Bbb N\rbrace$
is an orthonormal basis for $\lh$.

ii).\ $\tau_1(Q_r)$= $\sum\limits_{n \in \Bbb N}
\vert\vert P_r\eps_n\vert\vert^2.$

\endproclaim

We conclude now the prof of iii). in Theorem 3.2.

\proclaim{Proposition 3.10} Asssuming the hypothesis
from Proposition 3.9, it follows that
$$\sum
\limits_{n \in \Bbb N}\vert\vert P_r\eps_n\vert\vert^2
= d_{\pi_r} \text{covol}(\G) =\frac{r-1} {\pi} \text{covol}
(\G).$$
\endproclaim

Proof. The proof is now exactly as in [GHJ]. Since
our context is a bit different we will recall it anyway.

Let $\eta$ be any unit vector in
$$H_r =\hr =P_r(\lh).$$
Then $\pr(g) \eta= \eta, g \in G$ and
$$\vert\vert\tilde{\pi}_r(g)\eta\vert\vert_{\lh}=
\vert\vert\eta\vert\vert_{\lh}=1, \text{\ for\ all\ }
 g \in G.$$

Let $\Cal F$ be a fundamental domain for
$\G$ acting on the right on $G$.
Then covol($\G)$ is $\int\limits_{\Cal F} 1 \text{d} g$ and this
may be infinite. We have
$$\text{covol}(\G)=\int\limits_{\Cal F} 1 \text{d} g=
\int\limits_{\Cal F}
\vert\vert\tilde{\pi}_r(g)\eta\vert\vert_{\lh}\text{d} g=$$
$$=\sum\limits_{n\in \Bbb N, \gamma \in \G}
\int\limits_{\Cal F}\vert\langle \tilde{\pi}_r(g)\eta,
\pr(\gamma)\eps_n\rangle\vert^2\text{d} g=$$
$$=\sum\limits_{n\in \Bbb N}\int\limits_{G}\vert
\langle \tilde{\pi}_r(g)\eta,
\eps_n\rangle\vert^2\text{d} g=$$
$$=\sum\limits_{n\in \Bbb N}\int\limits_{G}\vert
\langle \pi_r(g)\eta,
P_r\eps_n\rangle\vert^2\text{d} g=$$
$$=\sum\limits_{n\in \Bbb N}(d_{\pi_r})^{-1}
\vert\vert \eta\vert\vert^2 \cdot
\vert\vert P_r \eps_n\vert\vert^2=$$
$$=\sum\limits_{n\in \Bbb N}(d_{\pi_r})^{-1}
\vert\vert P_r \eps_n\vert\vert^2=(d_{\pi_r})^{-1}
\tau_1(Q_r). $$
Here we used Proposition 3.1 (the genralized form of
Schurr orthogonality relations for projective representations).
This and Lemma 3.9 concludes the proof for Theorem
3.2, if one assumes that the orthogonality relations
hold.

Note that  we  proved in fact a slightly more precise
statement than the determination of the isomorphism
class of $\lbrace\pi_r(\G)\rbrace'$ in the case
of $\calg$ being a factor (compare [JHG]).

\proclaim{Proposition} Let $\G$ be a discrete
fuchsian subgroup of $G=\psl$ and assume that
$\calg$ is a factor.

Then the von Neumann algebra $\lbrace\pi_r(\G)\rbrace''$
 generated in $\b$
by $\pi_r(\G)$ is isomorphic to $\calg$.
Denote (following [J]) by $dim_M H$ the coupling constant
 ([MvN]) for a type $II_1$ factor $M$ acting on a Hilbert space
$H$. Then

$$dim_{\lbrace\pi_r(\G)\rbrace''}\hr=
dim_{\calg} \hr= \text{covol}(\G)d_{\pi_r}, r >1.$$
\endproclaim

To complete the proof of our theorem we need to
show that the projective unitary representations
$(\pi_r)_{r >1}$ have all square integrable coefficients in
the sense of Definition 3.1. We have
\proclaim{Proposition 3.11} The projective
unitary representations
$(\pi_r)_{r >1}$have all finite  formal dimension
$d_{\pi_r}=\frac{r-1} {\pi}$. In particular

$$\int\limits_{\psl}\vert\langle
\pi_r(g)\zeta,\eta\rangle_{H_r}\vert^2\text{d}g=
\frac{\pi} {r-1} \vert\vert \zeta\vert\vert^2
\cdot \vert\vert \eta\vert\vert^2, \zeta,\eta \in \hr.$$
\endproclaim

Proof. As in [HGJ]it is sufficient to prove this
statement for the case of projective representations
$\pi_r:$  $G=$ SU(1,1) $\rightarrow B(L^2(\Bbb D,\mu _r))$.

Recall that $\mu_r$ is the measure on the unit disk
$\Bbb D$ with density
$$z \rightarrow (1-\vert z\vert^2)^{r-2}$$
with respect to the Lesbegue measure on the disk.
The representations
$\pi_r$ act on $H^2(\Bbb D,\mu_r)$ accordingto the
formula
$$(\pi_r(g) f) (z)=(j(g,z))^{-r} f(g^{-1}z), f \in
H^2(\Bbb D,\mu_r), g \in G,z \in \Bbb D.$$

The modular factor $j(g,z)$ is now given by the formula
$$j(g,z)= (\overline {b}z +\overline {a}),z \in \Bbb D$$
 for $$g=\pmatrix a&b\\{\overline b}&{\overline a}
\endpmatrix \vert \  a,b \in \Bbb C,
|a|^2 - |b|^2 =1 $$
an arbitrary element in SU(1,1). To define $j(g,z)^{-r}$
ane chooses a normal branch for
$t=\text{arg}(\overline {b}z +\overline {a})$ with values in
the interval $\-\pi<t\leq \pi$.

 Because of definition 3.1 it is sufficient to prove
the statement for a single non zero
vector $\zeta=\eta$ and we choose this vector to be the
constant function 1 on $\Bbb D$.

We use the method in ([Ro],chapter 20). Let
$$g=\pmatrix a&b\\{\overline b}&{\overline a}
\endpmatrix,  a,b \in \Bbb C,
|a|^2 - |b|^2 =1,$$ be  an arbitrary element in SU(1,1).
Then
$$\vert \langle\pi_r(g)1,1\rangle\vert=
\vert\int\limits_{\Bbb D} \frac {1} {(\overline {b}z
 +\overline {a})^{-r} }
 (1-\vert z\vert^2)^{r-2} \text{d}z\text{d} \overline z
\vert =$$
$$=
\vert\int\limits_{0}^{1} (1-t^2)^{r-2}\lbrack\int\limits_{0}^{2\pi}
(\overline a +\overline b t e^{i\theta} )^{-r}
 \text{d}\theta\rbrack \text{d} t\vert=$$
$$=
\vert\int\limits_{0}^{1} (1-t^2)^{r-2}\lbrack\int
\limits_{\vert w \vert=1}(iw)^{-1}
(\overline a+\overline b w)^{-r} \text {d} w\rbrack \text{d} t\vert.$$
Since $\vert a\vert>\vert b\vert$,
 the function $ w\rightarrow (iw)^{-1}
(\overline a+\overline b w)^{-r}$ has its only pole in
$\Bbb D$ at $w=0$ and thus the above integral
is
$$\vert\int\limits_{0}^{1} (1-t^2)^{r-2}
\frac{2\pi} {(\overline a)^r}\text{d} t\vert=
\frac{\pi}{ (r-1) \vert a\vert^r}.$$
The proof now follows line by line the one in
[Ro] and we end up with the equality
$$\int\limits_{G} \vert\langle\pi_r(g)1, 1\rangle_{H_r}
\vert^2 \text{d} g = ((r-1)/\pi)^{-1} \vert\vert 1
\vert\vert^2_{H_r}.$$
(See also [GHJ] for a discution on the difference in constants
that arises by working with SU(1,1) instead of $\psl$).
This finishes the remainig part of the proof of Theorem
3.2.
\bigskip
We end this paragraph by explaining why the kernels
$k=k(z,\zeta)=k(\gamma z,\gamma \zeta)$,$\gamma \in \G$ on
$\Bbb D^2$, analytic in $z$ and antianalytic in $\zeta$
 should be considered as a  generalizition of automorphic forms.
We will discuss only the case of automorphic forms of
integral, even weight for $\G=PSL(2,\Bbb Z)$.

 Recall that
an automorphic form of even, integral weight
$2k$ for $\G=PSL(2,\Bbb Z)$ is an analytic function $f$
on $\Bbb H$ so that
$$f(\gamma z) = j(\gamma, z)^{-2k} f(z), z \in \Bbb H,\gamma \in \G$$
and
$$\vert f(z)\vert \leq \text{const}(\text{Im}\ z)^{-k}, z \in
\Bbb H.$$

It was shown in [GHJ] (see also Jones manuscript notes)
that if $f,g$ are automorphic forms of same
integral weight $2k$ then the linear
multiplication operators  $M^n_f$ and respectively $M^n_g$, on $H_n$,
 with the functions $f$ and
$g$ respectively are bounded, with values in
$H_{n+2k}$. Moreover both $M^n_f$ and $M^n_g$ are
intertwining operators for the representations
$\pi_n,\pi_{n+2k}$ resticted to $\G$, that is
$$M^n_f \pi_n(\gamma)=\pi_{n+2k}(\gamma)M^n_f$$
and similarly for $g$.

Hence $(M^n_f)^{\ast} M^n_g$ belongs to
$\Cal A_n$ and $M^n_g(M^n_f)^{\ast}$ belongs to
$\Cal A_{n+2k}$. Moreover the value of the (unique)
traces on $\Cal A_n$ and  $\Cal A_{n+2k}$ on
this elements is computed in [GHJ] and it is
equal (modulo constants depending on $n$ and $k$) to
the Petterson scalar product $\langle f,g\rangle_{\text{Pet}}$
 (see [Ma]) which is defined by
$$\langle f,g\rangle_{\text{Pet}}=\text{const}\int\limits_{F}
f(z) \overline{g(z)} \text{Im}^{2k-2} \text {d}\nu_0(z)=
\int\limits_{F}
f(z) \overline{g(z)} \text {d}\nu_{2k}(z).$$

We note that this computation is now generalized by the trace
formula in  Theorem
3.2. This follows from the following:

\proclaim{Proposition 3.12} Let $f,g$ be automorphic
 forms of
integral weight $2k$. Let $M^n_f$ and respectively $M^n_g$
be   the linear, continuous
multiplication operators  with
 the functions $f$ and
$g$ respectively, on $H_n$
 with values in
$H_{n+2k}$.  Then the Berezin symbol
$k=k(z,\zeta)=k(\gamma z,\gamma \zeta)$,$\gamma \in \G$,
for the operator $M^n_g(M^n_f)^{\ast}$ in
$\Cal A_{n+2k}\subseteq B(H_{2k+n})$ is given by the
formula:
$$k(z,\zeta)= \frac{c_n} {c_{n+2k}}\overline{f(\zeta)}g(z)
((z-\zeta)/2i)^{2k}, z,\zeta \in \Bbb H.$$
Note that the factor $((z-\zeta)/2i)^{2k}$ makes this
symbol $\G$ invariant. In particular (modulo a scalar), by Theorem 3.2,
 the trace of $M^n_g(M^n_f)^{\ast}$ in
$\Cal A_{n+2k}$ is $\int\limits_{F}
f(z) \overline{g(z)} \text{Im}^{2k-2} \text {d}\nu_0(z).$

\endproclaim

Proof. It is obvious that for all $z \in \Bbb H$ one
has:
$$(M^n_f)^{\ast}e^{n+2k}_z= \overline{f(z)}e^n_z.$$
Hence the symbol for $M^n_g(M^n_f)^{\ast}$ is
$$k(z,\zeta)=
\frac{\langle(M^n_f)^{\ast}e^{n+2k}_{\zeta},
 (M^n_g)^{\ast}e^{n+2k}_{z}\rangle}
{\langle e^{n+2k}_{\zeta},
 e^{n+2k}_{z}\rangle}=$$
$$\overline{f(\zeta)}g(z)
\frac{\langle e^n_{\zeta},e^n_z\rangle}
{\langle e^{n+2k}_{\zeta},
 e^{n+2k}_{z}\rangle},z,\zeta \in \Bbb H.$$

\proclaim {Remark 3.13}. This shows that the union of
all symbols in $\Cal A_r$ when $r$ tends to infinity
exhausts all possible pairs of automorphic
functions.
\endproclaim

We have to take $r\rightarrow \infty$ above, since
for fixed $r$, only the automorphic forms of weight
$2k< r-1$ may occur in symbols of bounded
operators in $\Cal A_r$by the above method.

\def\pr{\tilde{\pi}_r}
\def\epst{\tilde{\eps}}
\def\hr{H^2(\Bbb H,\nu_r)}
\def\tr{\text{tr}}
\def\lh{L^2(\Bbb H,\nu_r)}
\def\lbh{B(\lh)}
\def\eps{\epsilon}
\def\ovt{\overline{\otimes}}
\def\calg{\Cal L(\G)}
\def\lbg{\Cal L(\G)\overline{\otimes} B(L^2(F,\nu_r))}
\def\lfg{l^2(\G)\overline{\otimes} L^2(F,\nu_r)}
\def\lbf{B(L^2(F,\nu_r))}
\def\lf{L^2(F,\nu_r)}
\def \pt {\tilde{\pi}^0 _r}

\def\b{B(H_r)}
\def\nk{||k||_{\lambda ,r}}
\def\G{\Bbb {\Gamma}}
\def\h{\Bbb H}

\def\i{\int_ G}
\def\c{c_r}
\def\l{L^2 (G)}
\def\s{\sup\limits_{\zeta \in \Bbb H}}
\def\a{\hat A(z,\zeta)}
\def\ra{\rightarrow}

\def\dz{|d(z,\zeta)|^r}
\def\dn{d\nu _0 (\zeta)}
\def\k{k(z,\overline \zeta)}
\def\bh{\hat B(H_r)}
\def\g{\gamma}
\def\hg{\Bbb H /\Bbb {\Gamma}}
\def \e{\epsilon}
\def\n{||\ ||_{\lambda ,r}}
\def\na{||A||_{\infty ,r}}
\def\an{\langle e^{\frac{s-r}{2}}_z , e^{\frac{s-r}{2}}_z \rangle}
\def\z{\zeta}
\def\st{*_s}
\def\fm{|\langle\pi (g)\zeta,\eta \rangle _H|^2}

\bigskip

\def \pslr{$PSL(2,\Bbb R)$}
\def\pr{\tilde{\pi}_r}
\def\epst{\tilde{\eps}}
\def\hr{H^2(\Bbb H,\nu_r)}
\def\tr{\text{tr}}
\def\lh{L^2(\Bbb H,\nu_r)}
\def\lbh{B(\lh)}
\def\eps{\epsilon}
\def\ovt{\overline{\otimes}}
\def\calg{\Cal L(\G)}
\def\lbg{\Cal L(\G)\overline{\otimes} B(L^2(F,\nu_r))}
\def\lfg{l^2(\G)\overline{\otimes} L^2(F,\nu_r)}
\def\lbf{B(L^2(F,\nu_r))}
\def\lf{L^2(F,\nu_r)}
\def \pt {\tilde{\pi}^0 _r}

\def\b{B(H_r)}
\def\nk{||k||_{\lambda ,r}}
\def\G{\Bbb {\Gamma}}
\def\h{\Bbb H}

\def\i{\int_ G}
\def\c{c_r}
\def\l{L^2 (G)}
\def\s{\sup\limits_{\zeta \in \Bbb H}}
\def\a{\hat A(z,\zeta)}
\def\ra{\rightarrow}

\def\dz{|d(z,\zeta)|^r}
\def\dn{d\nu _0 (\zeta)}
\def\k{k(z,\overline \zeta)}
\def\bh{\hat B(H_r)}
\def\g{\gamma}
\def\hg{\Bbb H /\Bbb {\Gamma}}
\def \e{\epsilon}
\def\n{||\ ||_{\lambda ,r}}
\def\na{||A||_{\infty ,r}}
\def\an{\langle e^{\frac{s-r}{2}}_z , e^{\frac{s-r}{2}}_z \rangle}
\def\z{\zeta}
\def\st{*_s}
\def\fm{|\langle\pi (g)\zeta,\eta \rangle _H|^2}
\def\d{\text{d}}
\def\Im {\text{Im}}
\def\ovl{\overline}
\def\car{\Cal A_r}
\def\cas{\Cal A_s}
\def\cat{\Cal A_t}
\def\lar{L^2(\car)}
\def\las{L^2(\cas)}
\def\lat{L^2(\cat)}
\def\vv{\vert\vert}
\def\tbrf{T^r_{B^{-1/2}_rf}}
\def\jsr{j_{s,r}}
\def\brs{B^{-1/2}_sB^{1/2}_r}
\def\lr{l_{s,r}}
\centerline{ 4.The covariant symbol
 in invariant Berezin quantization}
\bigskip

In this paragraph we will analyse the
deformation quantization for $\Bbb H/\Gamma$
from the viewpoint of covariant
symbols.
We use the notation $M^r_f$ for the multiplication
operator  on $H^2(\Bbb H,\nu_r)$ with the
function $f$.
Let  $P_r$ be the orthogonal projection
from $L^2(\Bbb H,\nu_r)$ onto
$H^2(\Bbb H,\nu_r)$.
 Recall that
 $A\in B(H^2(\Bbb H,\nu_r))$
admits a contravariant Berezin's
 symbol $f$  in
$L^{\infty}(\Bbb H)$ if $A$ is a Toeplitz operator
on $H_r=H^2(\Bbb H,\nu_r)$ with symbol $f$,
that is
$A= T^r_f=P_rM^r_fP_r.$

The relation between the Berezin's contravariant
and covariant symbols
is  a duality relation involving the trace on $B(H_r)$.
 We use the notation
$\overset\text{$\circ$}\to{A}(z,\overline{z})=f(z)$ for the
 contravariant symbol of A. Let  $B$ be any element
in
$B(H_r)$ of covariant
symbol $\hat B (z,\overline z)$. Assume that
 $AB$ is a trace class operator. The duality relation
between the two type of  symbols is given by the
 following equality:
$$\leqno (4.0)\ \ \  \text{tr}_{B(H_r)}(AB)=
\int _{\Bbb H}
\hat B (z,\overline z)
 \overset\text{$\circ$}\to{A} (z,\overline z)
 d \nu _0 (z).$$

In this paragraph we will extend this
relation to the case of $\G$- invariant symbols. These
symbols correspond to
linear operators in $B(H_r)$ that commute with
$\pi_r(\G)$.Recall
that we used in Proposition 3.5 the notation
$\pr :G=PSL(2,\Bbb R)\ra \lbh$ for
 the projective, unitary
representation of $G$    on $\lh$
that is defined by the same algebraic formula
as $\pi_r$.

If $f$ is any bounded,
 measurable function on $\Bbb H$ then
$M^r_f$ commutes with $\pr(\G)$.
Moreover, $P_r$ commutes with $\pr(\G)$
 (and in fact with $\pr(G)$ as $\hr$
is invariated by $\pr(G)$).  Hence
 $T^r_f=P_rM^r_fP_r$ commutes with
$\pi_r(\G)$.
We have thus proved that if $A=T^r_f$ and
if $f$ is a $\G$-invariant and bounded function
on $\Bbb H$ then $A$ belongs to
$\Cal A_r= \lbrace \pi_r(\G)\rbrace'$.

The duality relation (4.0) between
the covariant and
contravariant symbol will now
be replaced with a new relation in which the
trace tr$_{B(H_r)}$ on $\b$ is replaced
by a trace on the semifinite von Neumann
algebra $\Cal A_r$. This will correspond
 to the fact that
in the formula (4.0) we will be rather integrating
over $F$, a fundamental domain for the action
of $\G$ on $\Bbb H$ rather then on $\Bbb H$ as in the
classical setting.

In the last part of this paragraph
we will introduce a third type
of a symbol for operators on
$H_r$. We will call this
type of symbol an intermadiate symbol
for operators in $\b$ because
 it inherits  properties from both
 the covariant and the contravariant
symbol.

Recall that $B_r(\Delta)$ was
the positive operator which assigns to
a function $f$ on $\Bbb H$ the restriction
to the
diagonal of the contravariant symbol of the
 associated Toeplitz
operator $T^r_f$.
The intermediate symbol for an operator
$A$ in $\b$ is the
 operator function
$\lbrack B_r(\Delta)\rbrack^{1/2}$ applied to the
the covariant symbol of $A$.

We start with a rigurous definition of the operator
$B_r(\Delta)$ in the $\G$-invariant case.

\proclaim{Proposition 4.1.} Let
$r>1$and
let $F$ be a fundamental domain for
the action of $\G$ on the upper half plane.
  Let f be any bounded function on $\Bbb H$ that is
$\G-$ invariant and let
$A=T^r_f=P_rM^r_fP_r$ be the Toeplitz
operator on $\hr$ with symbol $f$.

Then $A$ commutes with $\lbrace \pi_r(\G)\rbrace$
(with the notations in the third paragraph, this is
$A\in \Cal A_r$). Let
$\hat A=\hat A(z,\zeta)$ be the contravariant
symbol of $A=T^r_f$. Let
$K_r(z,\zeta)$ be the kernel function on $L^2(F)$
defined by
$$K_r(z,\eta)=\c \sum\limits_{\g \in \G}
\vert d(z,\g\eta)\vert^{2r},$$
where $\vert d(z,\eta)\vert^2=(\text{Im}\ z)
(\text{Im}\ \eta)
\vert z-\overline{\eta}\vert^2,z,\eta \in\Bbb H$
is a function on the hyperbolic distance between
$z$ and $\eta$ in $\Bbb H$.
Then
$$\hat A(z,\overline z)= (B_rf)(z)=$$
$$=\c \int\limits_{\Bbb H} f(\eta)
\vert d(z,\eta)\vert^{2r}\text{d}\nu_r(z)=
\int\limits_F K_r(z,\eta) f(\eta)\text{d} \nu_0(\eta).$$

Moreover the linear operator
$B_r$ defined above on $L^{\infty}(F)$ with
values in $L^{\infty}(F)$ extends to
a bounded, positive, contractive operator
on $L^2(F,\nu_r)$.  The operator $B_r$ is injective and the
operators $(B_r)_{r>1}$ are pairwise commuting.
Moreover $B_r$ tends strongly to 1 as $r$ tends
to infinity.
\endproclaim

Proof. The kernel $K_r(z,w), z,w \in F$ is symmetric
with positive values. Moreover
$$\int\limits_F K_r(z,w)\d \nu_0(w)=
\c\int\limits_{\Bbb H}
\vert d(z,w)\vert^{2r}\d \nu_0(w)=$$
$$=\c (\Im z)^r\int\limits_{\Bbb H}\frac {1}
{\vert(z-\ovl w)/2i\vert^{2r}} \d (w)=
\c^{-1} (\Im z)^r\int\limits_{\Bbb H}\frac
{\c^2} {\vert(z-\ovl w)/2i\vert^{2r}} \d (w)=$$
$$= \c^{-1} (\Im z)^r\langle e^r_z,e^r_z\rangle_{H_r}=
 \c^{-1} (\Im z)^r e^r_z(z)=$$
$$=\c^{-1} (\Im z)^r
\frac{\c} {((z-\ovl z)/2i)^r}=1,\text{\ for\ all\ }
z \in F.$$

Hence the interpolation arguments in Theorem 2.4,
page 1131, [Be] show that $B_r$
extends to a contractive operator
from $L^2(F,\nu_r)$ into $L^2(F,\nu_r)$
and it also extends to a bounded
operator from  $L^1(F,\nu_r)$ into
$L^1(F,\nu_r)$.

That $B_r$ is a positive injective operator follows from
the Corollary on page 66 in Patterson paper ([Pa]) in
Math. Proc. Cambr, (81). We  could have proved the
positivity of $B_r$ by using the corresponding
property of the similar operator
on $\lh$.

The pairwise commutativity of the
operators $B_r$ follows, for example, from the
above quoted paper ([Pa]). In fact the operators
$B_r$ are functions of the
invariant laplacian on
$\hg$. This type of analysis was first considered by Selberg
[Se]. This completes the proof.

\bigskip

We now prove the duality  relation
between the two type of Berezin's symbols
in the $\G-$ equivariant case. Since obviously
(at least in the infinite covolume case)
 the trace
formula can not hold true for all elements
in $\Cal A_r$, we will restrict our
 consideration  to elements
in $\hat {\Cal A}_r=\bh\cap \Cal A_r$.

\proclaim{Proposition 4.2}
Let $\Cal A_r$ be the commutant of
$\pi_r(\G)$ in $\b$ and let
$\hat {\Cal A}_r=\bh\cap \Cal A_r$.
  Let $r>1$ and let $A$
be any operator in $\hat {\Cal A}_r$.
Let $\hat A(z,\ovl{\zeta})$,
$z,\zeta \in \Bbb H$, be the contravariant symbol of the
operator $A$.
We choose a fundamental domain   $F$ for $\G$ in
$\Bbb H$ and let $f$ be $\G-$ equivariant
function on $\Bbb H$.
 Assume that $f$ is in
$L^1(F,\nu_0)$.
 Let $T^r_f=P_rM^r_fP_r$ be the Toeplitz
operator on $\hr$ with symbol $f$. Denote
by $\tau_{\Cal A_r}$
 the (semifinte) faithful
trace on $\Cal A_r$ that was constructed in
Theorem 3.2. Then
$$\tau_{\Cal A_r}(AT^r_f)=
(const)\int\limits_F \hat A(z.\ovl z) f(z)\d \nu_o(z).$$
The value of the constant in front
    of the integral
is (area\ $F)^{-1}$ in the finite covolume case
and 1 otherwise.
Moreover
$$\vert \tau_{\Cal A_r}(AT^r_f)\vert \leq
\c \vert\vert A\vert\vert_{\lambda,r}\vert\vert
f\vert\vert_{L^1(F,\nu_0)}.$$
\endproclaim

Proof.
The symbol for
the operator $AT^r_f$ is
the  iterarated integral:
$$k(z,\ovl{\zeta})=\c^2
((z-\ovl\zeta)/2i)^r\int\limits_{\Bbb H}
\int\limits_{\Bbb H}
\frac{\hat A(z,\ovl\eta)f(a)}
{\lbrack(z-\ovl\eta)/2i\rbrack^r
\lbrack(\eta-\ovl a)/2i\rbrack^r
\lbrack(a-\ovl \zeta)/2i\rbrack^r}\d\nu_r(a,\eta),$$
were we first integrate $a$ and then $\eta$.

Consequently
$$\tau_{\Cal A_r}(AT^r_f)=
(\text{const})\c^2\int\limits_{F,z}
\int\limits_{\Bbb H,\eta}\int\limits_{\Bbb H,a}
\frac{\hat A(z,\ovl\eta)f(a)}
{\lbrack(z-\ovl\eta)/2i\rbrack^r
\lbrack(\eta-\ovl a)/2i\rbrack^r
\lbrack(a-\ovl z)/2i\rbrack^r}\d\nu_r(z,a,\eta),$$
were we first integrate $a$, then $\eta$ and then $z$.
We indicated  above for each
 variable the domain of integration.
 By letting
the measure be $\d \nu_0(z,a,\eta$ and
 by collecting
all the densities into the integrand,
 the integrand itself
becomes a $\G-$invariant function in the variables
$z,a,\eta$ on $\Bbb H^3$. If the integral
 were to be absolute
convergent then  we could integrate
on any fundamental domain of
$\G$ acting on $\Bbb H^3$, e.g.we could
integrate on
 $\Bbb H_z \times F_a \times \Bbb H_{\eta}$.
In this case the integral would be (modulo a
constant):
$$\int\limits_F f(a)
\lbrace \c^2 \int\int\limits_{\Bbb H}
\frac{\hat A(z,\ovl\eta)}
{\lbrack(z-\ovl\eta)/2i\rbrack^r
\lbrack(\eta-\ovl a)/2i\rbrack^r
\lbrack(a-\ovl z)/2i\rbrack^r}\d\nu_r(z,\eta)
\rbrace\d\nu_r(a)=$$
$$\int\limits_F f(a)
\frac {\hat A(a,\ovl a)} {\lbrack (a-\ovl a)/2i\rbrack^r}
 \d \nu_r (a)= \int\limits_F f(a)
\hat A(a,\ovl a) \d \nu_0 (a).$$

To prove the absolute convergence of the
integrand it is sufficient to check this
 on any fundamental domain.
Thus it is sufficient to
estimate:
$$\leqno (4.1)\ \
\int\limits_F \vert f(a)\vert
 \lbrack\c^2 \int\int\limits_{\Bbb H}
\frac{\vert\hat A(z,\ovl\eta)\vert }
{\vert\lbrack(z-\ovl\eta)/2i\rbrack\vert^r
\vert\lbrack(\eta-\ovl a)/2i\rbrack\vert^r
\vert\lbrack(a-\ovl z)/2i\rbrack\vert^r}
\d\nu_r(z,\eta)
\rbrack\d\nu_r(a).$$

By Proposition 2.7, if $A$ belongs to $\bh$,
 then the operator  on
$\lh$,
with integral
kernel $\frac{ \vert\hat A(z,\ovl\eta)\vert}
{\vert\lbrack(z-\ovl\eta)/2i\rbrack\vert^r}$
 is bounded of uniform norm less
than
$\vert\vert A\vert\vert_{\lambda,r} $.
Moreover the function on $\Bbb H$ defined by
$z \ra \vert\lbrack(a-\ovl z)/2i\rbrack\vert^{-r}$,
belongs to $\lh$ and has
norm less than
$$\c\int\limits_{\Bbb H}  \vert\lbrack(a-\ovl z)/2i\rbrack\vert^{-2r}
\d\nu_r(z)=\c^{-1} \langle e^r_a,e^r_a\rangle_{H_r}=
(\text{Im}\ a)^{-r}.$$

Hence the inner integral in the formula
(4.1) is estimated by
$\vert\vert A\vert\vert_{\lambda,r}
 (\text{Im}\ a)^{-r}$
and thus the integral itself
is bounded by
$$(\text{const})\int\limits_F \vert f(a)\vert
\d \nu_0(a) =(\text{const})
\vert\vert A\vert\vert_{\lambda,r}\cdot
\vert\vert f\vert\vert_1.$$
This ends the proof.
\bigskip
In  the next proposition we  show
that  Toeplitz "operators" may be defined even if
 the
symbol function is not bounded, but only in $L^2$.
In this case   the corresponding Toeplitz
 operator will
be an element in $L^2(\Cal A_r)$.

\proclaim{Lemma 4.3} We identify $L^{\infty}(F)$
 with a subspace of all
bounded functions on $\Bbb H$, by extending them
outside $F$ by $\G- invariance$. Let $r>2$.
Let $S_r$ be the bounded linear operator on
$L^{\infty}(F)$ with values in $\Cal A_r$ defined
by $S_rf= T^r_f=P_r M^r_fP_r$.
Then $S_r$ extends to a contractive linear map from
$L^2(F)$ into $L^2(\Cal A_r)$.
Moreover $$\langle S^r_f, S^r_g\rangle _{L^2(\Cal A_r)}
=\langle B_rf,g\rangle_{L^2(F)},\ \text{for\ all\ }
 f,g \in L^2(F).$$

In particular,
$\{S^r_f\vert f \in L^2(F)\cap L^{\infty}(F) \}\subseteq
L^2(\Cal A_r)\cap \Cal A_r.$
\endproclaim
\bigskip

Proof.
We  compute
the term $\langle S^r_f,
S^r_g\rangle _{L^2(\Cal A_r)}$. We assume first
 that $f$ is in
$L^1(F)\cap L^{\infty}(F)$ and g is in $L^{\infty}
(F)$.
In this case  the
iterated integral defining
$\langle S^r_f,
S^r_g\rangle _{L^2(\Cal A_r)}$ is
(modulo a constant which is (area $F)^{-1}$ in the
finite covolume case and 1 otherwise)
$$
\int\limits_{F_z}((z-\ovl z/)2i)^r
\int\limits_{\Bbb H_{\eta}}
\int\int\limits_{\Bbb H^2_{a,b}}
\frac {f(a) \ovl{g(b)}}
{\lbrack(z-\ovl a)/2i\rbrack^r
\lbrack( a-\ovl \eta)/2i\rbrack^r
\lbrack (\eta  -\ovl b)/2i\rbrack^r}
\d \nu_r(a,b)\d \nu_r(\eta)\nu_0(z).$$

The above integral is an iterated integral: first we integrate
$a,b$ and then $\eta$ and $z$.

If the integral is absolutely convergent, then we may
integrate in any order the variables
$a,b,\eta,z$. In this  case, the
integrand with respect to  the measure
$\d \nu_0 (a,b,\eta,z)$ is $\G-$ invariant.
 To evaluate the integral, we may henceforth
change the domain of integration into
 $F_a\times \Bbb H_b \times \Bbb H _{\eta,z}$.
Under the absolute convergence
 assumption the integral
 will be thus equal to
$$(4.2)\int
\limits_{F_a}\int\limits_{\Bbb H_b}f(a)g(b)\lbrack
\int\int\limits_{\Bbb H_{\eta,z}}\frac
{\text{d}\nu_r(\eta)}
{\lbrack(a-\overline{\eta})/2i\rbrack^r
 \lbrack(\eta-\overline b)/2i\rbrack ^r}
\frac{\text {d}\nu_r(z)}
{\lbrack(a-\overline z)/2i\rbrack^r
\lbrack (z-\overline b)/2i\rbrack^r}
\rbrack \text{d}\nu_r(a,b)=$$
$$\int
\limits_{F_a}f(a)(\int\limits_{\Bbb H_b}g(b)
\vert\lbrack(a-\ovl b)/2i\rbrack\vert^{2r}\text{d}
\nu_r(b))\text {d}\nu_r(a)=
\langle f,B_rg\rangle_{L^2(F)}.$$
 By Fubini's theorem, to show
 the absolute convergence of the integral,
 it is sufficient to check
absolute convergence for the first integral in formula
(4.2). The integral of the
absolute value of the integrand is
bounded by
$$(4.3)\vert\vert g \vert\vert_{\infty}
 \int\limits_F\vert f(a)\vert\
M_r(a)\text{d}\nu_0(a),$$
where
$$M_r(a)=
\int\int\int\limits_{\Bbb H^3}
\vert d(a,\eta) d(\eta,b)d(b,z)d(z,a)\vert^r
\text{d}\nu_0(z,\eta,b).$$
Recall that we use the notation
$d(z,\ovl \zeta)=(\text{Im}\ z)^{1/2}
(\text{Im}\ \zeta)^{1/2}\lbrack (z-\ovl \zeta)/2i
\rbrack^{-1}$
and recall that $\vert d(z,\ovl \zeta)\vert$
 is a function of the hyperbolic distance
between $z$ and $\zeta$, $z,\zeta \in \Bbb H$.

It is easy to conclude that $M_r(a)$ is
a $PSL(2,\Bbb R)-$invariant function on $\Bbb H$. Since
$PSL(2,\Bbb R)$ acts transitively on $\Bbb H$
it is thus sufficient to check that the integral
defining $M_r(a)$ is convergent for a single value of $a$.
Also to check the finiteness of thie integral defining
$M_r(a)$ we may use the unit disk
$\Bbb D$ formalism instead
of the corresponding formalism for
 the upper half plane  $\Bbb H$.
We let $a=0 \in \Bbb D$ and hence we have to estimate:
$$\int\int\int\limits_{\Bbb D^3_{z,\eta,b}}
(1-\vert z\vert^2)^{r/2} (1-\vert\eta\vert)^{r/2}
(d(z,b))^r (d(\eta,b))^r \text{d}\nu_0(b)\text{d}
\nu_0(z,b).$$
By integrating first in the parameter
$b$ and then in the parameters
$z,\eta$ and since
the quantities
$$\sup\limits_z\int\limits_{\Bbb D} ((d(z,b))^{2r}
\text {d}\nu_0(b),\ \sup\limits_{\eta}\int\limits_{\Bbb D}
(d(\eta,b))^{2r}\text{d}\nu_0(b)$$
are finite, by using the Cauchy-Buniakowsky-Schwarz
inequality, we obtain that the above integral is dominated by
$$(\text{const})\int\int\limits_{\Bbb D^2}
(1-\vert z\vert^2)^{r/2} (1-\vert\eta\vert)^{r/2}
\text{d}
\nu_0(z,\eta),$$
which is finite if $r>2$.

Thus we have proved in particular that for
$f$ in $L^1(F)\cap L^{\infty}(F)$ the following equality
holds:
$$\langle S^rf,S^rf\rangle_{L^2(\Cal A_r)}=
\langle B_rf,f\rangle_{L^2(F)}.$$
As $B_r$ is bounded and contractive, this shows
that $$\vert\vert S^rf\vert\vert_{L^2(\Cal A_r)}
\leq
\vert\vert f \vert\vert_{L^2(F)}$$
Hence  the above equality also
extends to all $f \in L^2(F)$.

This concludes the proof of the lemma. Incidentaly
we have also proved the following:
\proclaim{Proposition 4.4} Let $r>2$.
Let $S^r:L^2(F)\rightarrow L^2(\Cal A_r)$ be
defined as above by
$$S^rf=T^r_f=P_r M^r_fP_r.$$
Then there exists
a constant $C_r>0$ so that for all
$g$ in $L^{\infty}(F)$ and for
all $f$ in $L^1(F,\nu_0)$ one has the
inequality
$$\vert\tau_{\Cal A_r}((S^rf)(S^rg))\vert
\leq C_r\vert\vert f \vert\vert_{L^1(F,\nu_0)}
\vert\vert g \vert\vert_{\infty}.$$
\endproclaim

In fact it is easier to understand the map $S^r$
by looking at its adjoint. In the next proposition we will
identify $(S^r)^{\ast}$ with the restiction map
$R_r$ on $\Cal A_r$, which associates to any kernel
$k(z,\ovl\zeta)$ on $\Bbb H^2$ its
restriction to the diagonal $z=\zeta$.

\proclaim{Proposition 4.5} Let $r>2$,
let $F$ be a fundamental domain for the action of
$\Gamma $ on $\Bbb H$. Let
$R_r:\Cal A_r  \rightarrow L^{\infty}(F)$ be the
map associating to an operator
$A$ in $\Cal A_r$ with kernel
$\hat A$, the function on $F$ defined by:
$z \rightarrow \hat A(z,\ovl z).$

Then $R_r$ extends to  a continuous linear
map $R_r:L^2(\Cal A_r,\tau)\rightarrow L^2(F)$.
Moreover $R_r$ is the adjoint of
$S^r$ and
$R^{\ast}_rR_r=B_r$.
\endproclaim
Proof.
Indeed we have allready checked that
$$\langle T^r_f,T^r_g\rangle_{\tau}=
\langle B_rf,g\rangle_{L^2(F)}
=\langle R_r(T^r_f),g\rangle_{L^2(F)},$$
for all $f,g \in L^2(F)$ or for all
$f \in L^1(F,\nu_0)$ and  all $g \in L^{\infty}(F)$.

This implies that
$$\langle T^r_f,S_r(g)\rangle_{L^2(\Cal A_r}=
\langle R_r(T^r_f),g
\rangle_{L^2(F)}$$
for all $g$ in $L^2(F)$ and all
$f$ in $L^{\infty}(F)$. Hence
$R_r$ is graph-contained (as an eventually
unbounded operator) in the adjoint of
$S_r$. Since $S_r$ is a bounded linear operator
and $R_r$ is closable it follows that $R_r$
extends to a bounded operator.
This completes the proof.
\bigskip

We will now introduce a third type of symbol which is intermadiate
betweeb Berezin's contravariant and covariant symbols. For
$A \in B(H_r)$ we let $V_rA$ to be the
the value of the operatorial inverse square root
of $B_r$ applied to the function obtained
by restriction to
the diagonal of the contravariant symbol
$\hat A$ of $A$.

When
using the intermadiate
symbol $V_rA$,
the main simplification   occurs   in computations
 involving the
scalar product: the scalar product
$\langle A,B\rangle_{L^2(\Cal A_r)}$
 of two elements
$A,B$ is equal to  the canonical scalar product in $L^2(F)$
of their intermadiate symbols $V_r A, V_r B$.

We will use this property to define
a different representation of the
Berezin deformation. In this representation the
trace of a product of two elements (viewed as
symbols functions on $\Bbb H$) doesn't depend on
the deformation parameter.
The main properties for the symbol $V_r$ are
outlined in the following proposition.

\proclaim{Proposition 4.6} Let
$U_r: L^2(F) \rightarrow L^2(\Cal A_r)$ be
the unitary operator defined by
$$U_rf=T^r_{B^{-1/2}_r}f, f\in L^2(F).$$
Note that $U_r$ is  first defined  on a dense set
and then extended by continuity.

For $A$ in $\Cal A_r$ with symbol
 $\hat A=\hat A(z,\ovl \zeta), z,\zeta \in \Bbb H$,
 let $R_rA$ be the
function on $\Bbb H$ defined by
$(R_r A)(z)= \hat A(z,\ovl z)$. The inverse
of $U_r$ is the unitary $V_r: L^2(\Cal A_r)
\rightarrow L^2(F)$ defined by the formula:
$$V_r A= B^{-1/2}_r (R_r A), A\in L^2(\Cal A_r).$$
Let $j_{s,r}: \Cal A_r \rightarrow \Cal A_s$
be the map that associates to
$A$ in $\Cal A_r$ the element $j_{s,r}(A)$
in $\Cal A_s$ having the same contravariant symbol
as $A$.

Then the following diagram is commutative:
\def\mapright#1{\smash{
  \mathop{\longrightarrow}\limits^{#1}}}

\def\mapup#1{\Big\uparrow
  \rlap{$\vcenter{\hbox{$\scriptstyle#1$}}$}}
$$\matrix L^2(\Cal A_r)&\mapright {j_{s,r}} &L^2(\Cal A_s)
\cr\mapup {U_r}&\
&\mapup {U_s}\cr L^2(F)&\mapright{B^{1/2}_rB^{-1/2}_s} &
L^2(F) \endmatrix$$
\endproclaim
Proof. We first define $U_r$ on the set $\Cal S=
B_r^{1/2}(L^2(F)\cap L^{\infty}(F))$. By
the injectivity
and continuity of $B_r$, $\Cal S$
 is a dense subspace
of $L^2(F)$.
For every vector $f$ in $B_r^{1/2}(L^2(F)\cap L^{\infty}(F))$
we have that
$$\vert\vert U_rf\vert\vert^2_{L^2(\Cal A_r)}=
\langle T^r_{B^{-1/2}_r}f,T^r_{B^{-1/2}_r}f\rangle_
{L^2(\Cal A_r)}.$$
By Proposition 4.4 this is
$$\langle B_rB_r^{-1/2}f,B_r^{-1/2}f\rangle_{L^2(F)}=
\vert\vert f\vert\vert_{L^2(F)}.$$
Thus $U_r$ extends by continuity
to an isometry on $L^2(F)$. To
prove that $U_r$ is in addition
an unitary, it will be sufficient to check
that $V_r$ is a left inverse for $U_r$.

First we note that $V_r$ is also a well defined isometry.
Observe that if $A$ is of the form
$T^r_f$ with $f$ a $\G-$equivariant function
on $\Bbb H$, then $V_r$ is well defined as
$V_r(T^r_f)=B_r^{1/2}f.$
Hence $V_r$ is  well defined on the
following dense subset of
$L^2(\Cal A_r)$:
$$\lbrace T^r_f\vert\ f \in L^2(F)\cap L^{\infty}(F)
\rbrace.$$
Morever
$$\vert\vert V_r A \vert\vert ^2_{L^2(F)}
=\langle V_r A,V_r A\rangle_{L^2(F)}=$$
$$
\langle  B^{-1/2}_r R_r A, B^{-1/2}_r R_r A\rangle
_{L^2(F)}=
\langle B^{-1}_r B_r f, B_r f\rangle_{L^2(F)}$$
and by Proposition 4.4 this is
$\vert\vert T^r_f\vv_{\lar}.$
Thus $V_r$ also extends to an isometry on $\lar$.

The restriction to the
diagonal of the contravariant symbol of the
operator $U_rf$ is equal to  the restriction to
the diagonal of the
symbol of $\tbrf$ which is
$$B_r(B^{-1/2}_r)f=B^{1/2}_rf, f
\in L^2(F).$$

Hence the restiction to the diagonal of
$U_r(V_rA)$ is equal to
$$B^{1/2}_r(V_r A)=B^{1/2}_r B^{-1/2}_r R_r A =
 R_r A$$
for all $A \in \car\cap\lar$. Thus for such  $A'$s,
the restriction  to the diagonal of the symbol of
$U_r(V_rA)$ and the restriction to the diagonal
of the symbol of $A$ are equal. Hence
$U_r(V_rA)=A$ for $A$ in a dense set and
hence $$U_rV_r=\text{Id}_{\lar}.$$

Similarly, for $f \in L^2(F)\cap L^{\infty}(F)$
$$V_rU_r f = V_r(\tbrf)=
B^{-1/2}_r R_r \tbrf=B^{-1/2}_rB_rB^{1/2}_rf=f.$$
Thus $U_r,V_r$ are unitaries, inverse one to
the other.

To complete the proof of the proposition it remains
to check the commutativity for the diagram.
 It is
obvious that $j_{s,r}$ is bounded with respect to the
$L^2$ norms on the corresponding spaces
(since the absolute value of the function
 $d$ entering the formulae
for these norms takes only subunitary values).
Also it  will be proved bellow  that the
operator $B^{1/2}_rB^{-1/2}_s$ is  bounded and
contractive. Hence it will be sufficient to check
 the commutativity of the diagram, for vectors
 in a dense set.

For $f$ in $L^2(F)\cap L^{\infty}(F)$,
 the restiction to the diagonal of
$j_{s,r}(U_rf)$ is
$$R_s(\tbrf)=R_r(\tbrf)=B^{1/2}_rf.$$
On the other hand
$$R_s(U_s(B^{-1/2}_sB^{1/2}_r f))=
R_sT^s_{B^{-1/2}_sB^{1/2}_rB^{-1/2}_s}=$$
$$
R_sT^s_{B^{-1}_sB^{1/2}_r}=
B_sB^{-1}_sB^{1/2}_rf=B^{1/2}_rf.$$
 Since any two
elements in $\cas$ whose symbols  coincide
on diagonal, are equal it follows
that the diagramm is commutative.

\proclaim{Corollary 4.7} For $s\geq r >1$, the
vector space $j_{s,r}(\car)$ is weakly
dense in $\cas$.
Hence, for $r >3$, the algebra $\hat {\car}$
is weakly dense in $\car$.
\endproclaim
Proof. The statement is equivalent to showing that
$j_{s,r}(\lar)$ is dense in $\las$.
By the commutativity of the  diagramm
in the previous proposition, this is
equivalent to proving that
$B^{-1/2}_sB^{1/2}_r$ has a dense image which
is also equivalent (as $B^{-1/2}_sB^{1/2}_r$
is selfadjoint) to showing that this bounded operator
has trivial kernel which will be proved by the
arguments at the end of this paragraph.

To prove that the algebra $\hat {\car}$
is weakly dense in $\car$ we use the fact
(allready proved in the paragarph 2 that for
all strictly positive $\eps$, the algebra
$\Cal A_{r-2-\eps}$ is contained in $\hat {\car}$.
This completes the proof of the corollary.

The diagramm in Proposition 4.7 also shows that
the operatorial absolute value $\vert j_{s,r}\vert$
in the polar decomposition of the inclusion
map $j_{s,r}:\car \ra\cas$ is unitary
equivalent (by the unitary $V_r$) to the
operator $B^{-1/2}_sB^{1/2}_r$.

This will be usefull in understanding the
differentiation,
 with respect the deformation parameter,
of the Berezin quantization.

\proclaim{Corolllary 4.8} Let $s\geq r >1$
and let $\vert j_{s,r}\vert$ in
$B(\lar)$ be the operatorial absolute value
in the polar decomposition of the inclusion
map $j_{s,r}:\car \ra\cas$
(i.e. $\vert j_{s,r}\vert=((\jsr)^{\ast}\jsr)^{1/2}$).

Then the following diagramm is commutative:

$$\matrix L^2(\Cal A_r)&
\mapright {\vert j_{s,r}\vert} &L^2(\Cal A_r)
\cr\mapup {U_r}&\
&\mapup {U_r}\cr L^2(F)&\mapright{B^{1/2}_rB^{-1/2}_s} &
L^2(F) \endmatrix$$
\endproclaim

Proof.
By the previous corollary we have that
$$(\jsr)^{\ast}\jsr=
U_r(\brs)U^{\ast}_sU_s(\brs)U^{\ast}_r=
U_r(B^{-1}_sB_r)U^{\ast}_r.$$
The assertion now follows from the fact
that $B^{-1}_sB_r$ is a positive
operator.

\bigskip

We will now show that the
"coefficient" function
$s\ra \langle U_s,\ U_t g\rangle_t$
is differentiable at any
point $s$ for $f,g$ in a dense subspace.
Moreover its derivative will be equal to an
 expression involving the bilinear functional giving the
derivative (in the deformation parameter) of
the scalar product.
This will allow us to simplify the epression of
the derivative of the Berezin product when using
intermadiate symbols.

\proclaim{Corrolary 4.10}
Let $t>1$ and let $f,g$ be two arbitrary
vectors in $L^2(F)$ so that the
function $\alpha$ given  by the formula
$$\alpha(s)=\langle U_sf,\ U_t g\rangle_{\las}$$
is defined in a neighbourhood of $t$
and is differentiable at $s=t$.
We also assume that the function
$\psi(s)$ defined by
$$\psi(s)=\langle U_sf,U_tg\rangle_{\lat}$$
is also defined in a neighbourhood of $t$.

Then the function $\psi(s)$ is also differentiable at
$t$ and
$$\psi'(t)=\frac {\text{d}}{\text{d}s}
 \langle U_sf,U_tg\rangle_{\lat}=
-1/2\alpha'(t).$$
\endproclaim
Proof. The proof is essentialy contained in the following
observation:
\proclaim{Observation}
Let $D$ be a selfadjoint (eventually) unbounded
operator acting on a Hilbert space $H$. Let
$a_p=a_p(D)$ be positive, selfadjoint,injective
(eventually unbounded) operators, that are
functions on $D$ (in the sense of the borelian
functional calculus), for $p$ in an
interval $(t-\eps,t+\eps)$. We assume that
$a_t(D) =\text{Id}_H$ and that the map
$p\ra a_p(D)\eta$ is
continuously strongly differentaiable in $p$
in a neighbourhood of $p$, for $\eta$ in a dense
subset of $H$.

Let $f,g$ be two vectors in $H$ so that the
function $\alpha(p)= \langle a_p(D)f,g\rangle_H$
is defined in a neighbourhood of $t$ and it is
differentiable at $p=t$. We assume also that the
function $\beta(p)=\langle a^{-1/2}_pf,g\rangle$ is
defined in a neighbourhood of $t$.

Then $\beta$
is differentiable $p=t$ and
$\beta'(t)=-1/2\alpha'(t).$
\endproclaim

Proof. With no loss of generality we may assume
that $D$ has multiplicity $1$. The setting is then
the following: we are given a measure $\mu$
on a subset $\sigma=\sigma(D)$  and functions
$a_p=a_p(x), x \in \sigma$, positive on $\sigma$
and non vanishing at any atom of $\mu$. The operator
$D$ is then
the operator of multiplication with the independent
 variable
 on $H=L^2(d,\mu)$.

The vectors $f,g$ are now simply two
functions in $L^2(d,\mu)$. By hypothesis,
for $\mu-$ almost all $x\in \sigma$, the
application $p\ra a_p(x)$ is differentiable
(and finite) in a neighbourhood of $t$.

 Also by hypothesis,  the functions on $\sigma$,
 in the variable
$p$, are defined in a neighbourhood of $t$ by
$$\alpha(p)=\int\limits_\sigma a_p(x)f(x)\ovl{g(x)}
\text{d}\mu(x)$$
and
$$\beta(p)=\int\limits_\sigma (a_p(x))^{-1/2}
f(x)\ovl{g(x)}
\text{d}\mu(x)$$
are the integrals are absolutely convergent.
Finnaly, the hypothesis also gives that
$\alpha(p)$ is differentiable at $p=t$ and that
$a_t$ is the constant function 1.

The conclusion of the observation then simply
 follows from
Lesbegue theorem of differentiability under the
integral sign.

Proof of the Corollary. Proof of the Corollary.

We have $\langle U_s f,U_t g\rangle _t =
\langle a_{t,s}f,g\rangle _{L^2(F)}=
\langle a^{-1}_{s,t}f,g\rangle _{L^2(F)}$.

 Moreover $\langle a^{2}_{s,t}f,g\rangle
=\langle |j_{s,t}|^2 U_t f,U_t g\rangle =
\langle U_t f,U_t g\rangle _s$
(we use here the following property of $j_{s,t}$:
$$\langle j_{s,t}^* j_{s,t} A,B \rangle _t =
\langle A,B \rangle _s \
\text{ for all } A,B \text{ in }L^2(\car)).$$
Thus if the derivative
 $\frac{d}{ds} \langle U_t f,U_t g\rangle _s \vert_{s=t}$
exists then this will imply that the derivative
$$\frac{d}{ds} \langle a^2_{st} f, g\rangle _s \vert_{s=t}
\ \text{ exists}.$$
The proof now follows from the observation.

Finnaly to describe the multiplication of Berezin
symbols of operators in the expression given
by the intermediate symbols and to be able
to define the differentiation of this multiplication
with respect to the deformation parameter
we will have to find a dense set of vectors that
are well behaved under this type of operations.
This is realized in the following lemma.

\proclaim{Lemma 4.10}
Let $r>1$ be fixed and let $(a,b)$ be an interval
with $r<a<b$. Then there exists a weakly
dense set $\Cal E$ in $L^{\infty}(F)$
so that $\Cal E\cap L^2(F)$ is dense in
$L^2(F)$ and so that for all
$s\in (a,b)$ and all $f \in \Cal E$,
$U_sf=T^s_{B^{-1/2}_s}$ belongs to $\cas\cap\las$
and the function on $(a,b)$ defined by
$$s\ra \vv U_sf\vv_{\infty,r}$$
is locally bounded.
\endproclaim
Proof. We will let $\Cal E$ be the vector space
$B^{3/2}_r(L^{\infty}(F)).$ It is clear that
$\Cal E\cap L^2(F)$ is dense in
$L^2(F)$ as this follows from the fact that
$$B^{3/2}_r(L^{\infty}(F)\cap L^2(F))
\subseteq B^{3/2}_r(L^{\infty}(F))\cap L^2(F)
\subseteq \Cal E\cap L^2(F).$$
The first inclusion follows from the fact that
$B^{3/2}_r$ maps $L^2(F)$ into $L^2(F)$ and thus
$B^{3/2}_r(L^{\infty}(F)\cap L^2(F))
\subseteq L^2(F)$. On the other hand
$L^{\infty}(F)\cap L^2(F)$ is dense in
$L^2(F)$. Since $B^{3/2}_r$ is continuous it
follows that $B^{3/2}_r(L^{\infty}(F)\cap L^2(F))$
is dense in $B^{3/2}_r(L^2(F))$. Since
$B^{3/2}_r$ has a dense range (as it is a selfadjoint
operator with zero kernel) it follows that
also $B^{3/2}_r(L^{\infty}(F)\cap L^2(F))$
is dense in $L^2(F).$

We want to prove that for any $s$ and any
$f\in \Cal E\cap L^2(F)$ there exists
$g$ (depending on $s$ and $f$) so that
$T^r_g$ is bounded ans so that
$$U_sf=T^s_{B^{-1/2}_s f}=
\jsr(T^r_g)=\jsr(U_r(B^{1/2}_rg)).$$
Since $(U_s)^{\ast}\jsr U_r=\brs$
(by Corollary 4.7) the above equality is equivalent
to
$$f=B^{1/2}_rB^{-1/2}_sB^{1/2}_rg=B^{-1/2}_sB_r g$$
and hence this is equivalent to
$$g=B^{1/2}_sB^{-1}_rf.$$
We need to prove that with this $g$,
$T^r_g$ is bounded. As $f$ is an element
in $B^{3/2}_r(L^{\infty}(F))$ and thus
$f=B^{3/2}_r(\theta)$ for some
$\theta$ in $L^{\infty}(F)$.
Hence
$$g=B^{1/2}_sB^{-1}_rf=B^{1/2}_sB^{1/2}_r\theta,
\theta \in L^{\infty}(F).$$
Thus, to complete the proof of the statement
it will be sufficient the bounded
operator $B^{1/2}_sB^{1/2}_r$ maps
$L^{\infty}(F)$ into $L^{\infty}(F)$ and that
the function $s\ra \vv U_s\theta\vv_{\infty,r}$
is locally bounded for all $\theta$ in
$L^{\infty}(F)$.

To do this it will be sufficient
to show that the operator
$B^{1/2}_sB^{1/2}_r$ in $B(L^2(F))$, is given,
analoguous to $B_r$, by an integral,
$PSL(2,\Bbb R)-$ invariant, kernel function
on $\Bbb H$ ([Ku],[Sel]).

More precisely it is thus sufficient to
find a kernel function $L=
L_{s,r}(z,\ovl w)$ on $\Bbb H^2$ so that
$$B^{1/2}_sB^{1/2}_rf(z)=\int\limits_F
L_{s,r}(z,\ovl w)f(w)\text{d}\nu_0(w), z\in \Bbb H.$$
Moreover it is required that there exists
a  kernel function
$l_{s,r}=l_{s,r}(z,\ovl w),z,w\in \Bbb H$
on $\Bbb H^2.$ The following properties
should hold true for $l_{s,r}$:

a).$l_{s,r}$ is $PSL(2,\Bbb R)-$, diagonally invariant,
that is $$\lr(gz,\ovl{gw})=\lr(z,\ovl{w}), z,w \in \Bbb H,
g\in  PSL(2,\Bbb R).$$

b).$\lr(z,\ovl{w})=\sum\limits_{\gamma}
\lr(\gamma z,\ovl{w}), z,w \in \Bbb H.$

c). The expression
$$M(s)=\sup\limits_{\z \in \Bbb H}
\int\limits_{\Bbb H}\vert \lr(z,\ovl w)\vert
\text {d}\nu_0$$
is finite and moreover the function
$M(s)$ is locally finite.

Assume we have find such an $\lr$ having the
properties a), b), c). Then for any
$\theta \in B^{3/2}_r(L^{\infty}(F)\cap L^2(F)$
we will have that
$$B^{1/2}_sB^{1/2}_rf(z)=
\int\limits_{\Bbb H}
l_{s,r}(z,\ovl w)\theta(w)
\text{d}\nu_0(w), z\in \Bbb H.$$
Hence $s\ra \text{sup} \vv B^{1/2}_r B^{1/2}_s \theta \vv_{\infty}$
is locally bounded.

To prove that the operator $B^{1/2}_r B^{1/2}_s$ has the required
property we will use the technique of the Selberg transform
(which in fact as mentioned in [Ve] is a particular case of a
more general transform).

Assume that $B_r, B_s$ are given by the functions $h_r,h_s$ as a
function of the invariant laplacian. Then $(h_r h_s)^{1/2}$ is
the function corresponding to the operator $B^{1/2}_r B^{1/2}_s$.

Let $\tilde{h}_r (t),\tilde{h}_s (t)$ be the functions
$h_r (t^2 +1/4), \  h_s (t^2 +1/4)$. Let $\phi _{r,s}(t)=
(h_rh_s)^{1/2}(t^2+1/4)$.

By Selberg Theorem (see also [Za]), $(h_rh_s)^{1/2}$ will be
represented by a kernel $k_{r,s}$ as we wanted if $\phi _{rs}$
has an holomorphic continuation in the strip $|\text{Im} t|<1/2$
and $\phi_{rs}(t)$ is of rapid decay in this strip.

On the other hand, Berezin formula for $B_r$ as function of the
laplacian shows that
$$ \tilde{h}_r (t)=\prod\limits^{\infty}_{n=1}[1+(t^2+1/4)[(1/r+n)
(1/r+(n-1)]^{-1}]^{-1}.$$
This expression shows that $\phi_{rs}$ has the required property.

\def\pr{\tilde{\pi}_r}
\def\epst{\tilde{\eps}}
\def\hr{H^2(\Bbb H,\nu_r)}
\def\tr{\text{tr}}
\def\lh{L^2(\Bbb H,\nu_r)}
\def\lbh{B(\lh)}
\def\eps{\epsilon}
\def\ovt{\overline{\otimes}}
\def\calg{\Cal L(\G)}
\def\lbg{\Cal L(\G)\overline{\otimes} B(L^2(F,\nu_r))}
\def\lfg{l^2(\G)\overline{\otimes} L^2(F,\nu_r)}
\def\lbf{B(L^2(F,\nu_r))}
\def\lf{L^2(F,\nu_r)}
\def \pt {\tilde{\pi}^0 _r}

\def\b{B(H_r)}
\def\nk{||k||_{\lambda ,r}}
\def\G{\Bbb {\Gamma}}
\def\h{\Bbb H}

\def\i{\int_ G}
\def\c{c_r}
\def\l{L^2 (G)}
\def\s{\sup\limits_{\zeta \in \Bbb H}}
\def\a{\hat A(z,\zeta)}
\def\ra{\rightarrow}

\def\dz{|d(z,\zeta)|^r}
\def\dn{d\nu _0 (\zeta)}
\def\k{k(z,\overline \zeta)}
\def\bh{\hat B(H_r)}
\def\g{\gamma}
\def\hg{\Bbb H /\Bbb {\Gamma}}
\def \e{\epsilon}
\def\n{||\ ||_{\lambda ,r}}
\def\na{||A||_{\infty ,r}}
\def\an{\langle e^{\frac{s-r}{2}}_z , e^{\frac{s-r}{2}}_z \rangle}
\def\z{\zeta}
\def\st{*_s}
\def\fm{|\langle\pi (g)\zeta,\eta \rangle _H|^2}

\centerline{5. A cyclic 2-cocycle associated to
 a deformation quantization}

In the paragraph 3 we have  constructed a family
 of semifinite
von Neumann algebras $\Cal A_r \subseteq B(H_r)$
 which are a deformation
quantization, in the sense of Berezin, for $\hg$.
In this section we are constructing a cyclic
 2-cocycle which is defined on a weakly dense
subalgebra of $\Cal A_r$, for each $r$.
This 2-cocycle  is associated in a canonical way
with the deformation.
It is likely that an abstract setting for
 this construction should be
found in the general machinery developed in [R.N,B.T].

In the particular case of the deformation for $\hg$, the cyclic
 2-cocycle we obtain this way  is very similar in  form
 with  the cyclic
cocycles that are constructed in the paper by
 Connes and Moscovici ([CM]).
 We will prove that that the cyclic 2-cocycle of the equivariant
deformation may be
obtained from a $\G-$invariant Alexander-Spanier two cocycle on
$\Bbb H$ by a procedure very similar to the constructions in
 chapter 4 in the above
mentioned paper.

The relation between the cyclic 2-cocycle and the deformation
becomes more transparent if one uses the intermediate
 type of Berezin's symbols $A\ra U_rA$ that
we have  introduced in the preceding paragraph.
This happens due to the fact that
 this type of symbols have the advantage that
 the trace of a product of two symbols is constant in the
deformation parameter.
 We first introduce  a rather formal abstract setting in
which this type of cocycles may be constructed. In the
 last part
of this paragraph  we will
perform more precise computations for
 the Berezin's, $\G-$ invariant
quantization.

\proclaim{Definition 5.1}
 Let $(\Cal B_s)_{s \in (a,b)}$ be a family of
semifinite von Neumann algebras.
 Denote by $(\ast_s)_{s \in (a,b)}$
the corresponding products operations
 on these algebras.
 We will call
$(\Cal B_s,\ast_s)_{s \in (a,b)}$  a "nice deformation"
 if the
following properties hold true:
\item {i.}\ For $r\o s$, $\Cal B_r$ is contained as a vector subspace
in $\Cal B_s$. Let $\jsr$ be the corresponding
inclusion map. We assume that $\Cal B_r$ is weakly
dense in $\Cal B_s$ and that $\jsr$ is weakly and normic
continuous. Moreover,  for all $s\geq r$,
$\jsr$ preserves the involution and the unit.
\item{ii.}There exists a linear functional $\tau$ on a
subalgebra of the union
 $\mathop\cup\limits_{s \in (a,b)}\Cal B_s$, so that
for each $s$, $\tau$ is a semifinite, normal faithfull trace
on $\Cal B_s$. In addition the maps $\jsr$ are trace
preserving.
\item {iii.} Let $L^2(\cbs,\tau)$ be the Hilbert
of  the Gelfand-Naimark-Segal construction for
$\tau$ on $\cbs$. We assume that
$\jsr$ extends to a contractive linear map
from $\lbr$ into $\lbs$. Moreover,
$\jsr$ maps the positive part of $\cbr$
into the positive part of $\cbs$.(In fact we are not
going to make use of this last property for
$\jsr$.)\item {iv.} There exist a selfadjoint
vector subspace $\Cal D
\subseteq \mathop\cap\limits_{s\in (a,b)}\cbs$
which is closed under all the multiplication operations
$(\ast_s)_{s \in (a,b)}$.  Also, the space
$\Cal D$ be dense in  the Hilbert spaces
$\lbs$, for all $s$.
\endproclaim
\proclaim{Remark 5.2} All of the above  properties
hold true for the family of algebras
$(\cas)_{s \in (a,b)}$ in the (equivariant) Berezin
 deformation quantization
of $\hg$ that was constructed in Theorem 3.2.
In this case, one could let $\Cal D$  be
 any of the algebras $\car$ for
$r <a-2$, or one might take $\Cal D=\hat\Cal A_a$.
\endproclaim
For a "nice" deformation" as above, by requiring some
additional properties, we construct  a cyclic
2-cocycle, which measures, to a certain extent,
the obstruction on the algerbras in
the deformation to be  mapped isomorphically one onto the
other
by  a family of isomorphisms depending smoothly on the deformation
parameter.
\proclaim{Definition 5.3}
Let $(\Cal B_s)_{s \in (a,b)}$
 be a "nice deformation"
as in the  Definition 5.1. We will call
$(\Cal B_s)_{s \in (a,b)}$  a
"nice differentiable deformation"
if in addition there exists
weakly dense, selfadjoint subalgebras
$\hat\cbs\subseteq\cbs,$ for all $s$,
with the following properties:
\item{i).} The algebras $(\hat\cbs)_{s \in (a,b)}$
are all  unital with the same unit
as $\cbs$. Moreover $\lbs \cap \hat\cbs$ is
weakly dense in $\cbs$.
Also, we let $\ns$ be a Banach algebra norm
 on the algebras $\hat\cbs$ for all $s$. The unit balls
for the norms $\ns$ are
 are weakly compact.
Moreover $||B^{\ast} ||_{\lambda ,s}=
||B ||_{\lambda ,s}$ for all $s$.

\item{ii).}  For
$s \geq r$, the inclusions $\jsr$,
map $\hat\cbr$ continuously
into $\hat\cbs$ with respect to the
norms $\n$ and $\ns$, correspondingly.
Moreover $\hat\cbr$ is closed under  the products
$(\ast)_s$ for all $s \geq r$. There exist positive
 constants
$c_{r,s}$ so that the function $s \ra c_{r,s} $ is locally
bounded for all $r$ and so that
$\vv A\ast_sB\vv_{\lambda ,r}
\leq c_{r,s}\vv A\vv_{\lambda ,r}\vv B\vv_{\lambda ,r},$
for all $A,B$ in $\hat \Cal B_r$.
\item {iii).} The space $\Cal D$ is contained in
$\mathop\cap\limits_s(\hat\cbs\cap\lbs)$.
\item {iv).}For all $s\in (a,b)$ The following functionals are well
defined on $\Cal D$.
\item {a).}$\mu_t(c,(a,b))= \frac {\text{d}}
 {{\text d}s}\tau(c\ast_t(a \ast_s b))\vert_{s=t}$.
\item{b).}$\phi_t(a,b)=\dds \tau(a\ast_sb)\vert_{s=t}.$
\item{c).}$\theta_t(a,b,c)=
\dds \tau(a \ast_sb\ast_sc)\vert_{s=t}.$
\endproclaim
Again we note that all of the above conditions will
hold true for the deformation quantization from
Theorem 3.2, with $\cbs=\cas$ and
$\hat\cbs=\hat\cas=\cas \cap \hat B(H_s)$, for all
$s \in (a,b)$. We also let
 $\Cal D$ be
$\car\cap\lar$ for some $r<a-2 $ or we let $\Cal D=
\hat\Cal A_a \cap L^2(\Cal A_a)$.

The above cocycles have the  formal properties listed in the
following proposition.
\proclaim{Proposition 5.4}
Let $\Cal D\subseteq\hat\cbs\subseteq\cbs $, $s\in (a,b)$
be a "nice differentiable deformation" as
in Definition 5.3. Let $t$ be fixed in $(a,b)$.
 Define on $\Cal D$
the following additional cocycles:
$$\alpha_t(a,b,c)=
\phi_t(a\ast_tb,c)+\phi_t(b\ast_tc,a)+
\phi_t(c\ast_t a,b),$$
$$\psi_t(a,b,c)=
\theta_t(a,b,c)-(1/2)\alpha_t(a,b,c),$$
for  $a,b,c \in \Cal D$. Then the following properties
hold true:
\item{i).} The linear
 functionals $\psi_t,\theta_t,\alpha_t$ are
cyclic,  that is
$\psi_t(a,b,c)=\psi_t(b,c,a)$  and  similarly
for $\theta_t,\alpha_t$. Moreover
$\phi_t$ is antisymmetric, that
is $\phi_t(a,b)=-\phi_t(b,a)$, for all
$a,b,c\in \Cal D$.
\item {ii).} For all $a,b,c \in \Cal D$ one has
$\theta_t(a,b,c)=
\mu_t(c,(a,b))+\phi_t(a\ast_tb,c).$
\item {iii.}$\psi_t$ belongs to
$Z^2_{\lambda}(\Cal D,\Bbb C)$, that is $\psi_t$
is a cyclic two cocycle in the sense of
Connes' cyclic cohomology (Co]):
$$\psi_t(a,b,c)=\psi_t(b,c,a),$$
$$\psi_t(a\st b,c,d)-
\psi_t(a, b\st c,d)+
\psi_t(a, b,c\st d)-
\psi_t(d \st a, b,c)=0,a,b,c,d \in \Cal D.$$
\item{iv.} $\mu_t$ is a Hochschild 2-cocycle,
that is $\mu$ verifies the second property
listed for $\psi_t$ above.
\item{v.}The following equality holds true
$\ovl{\psi_t(a,b,c)}=
\psi_t(b^{\ast},a^{\ast},c^{\ast})$
for all $a,b,c \in \Cal D$. Morever
if $1 \in \Cal D$, (which corresponds to the case
when all the algebras $\cbs$ are finite) then
$\psi_t(1,b,c)=0$ for all $b,c$ in $\Cal D$.
\endproclaim
Proof.
The statement in i). follows from the fact that
$\tau$ is a trace. The assertion
ii). is a consequence of the
product rule for differentiation:
$$\dds\tau(a\ss b\ss c)\vert_{s=t}=
\dds\tau(c\ss(a\st b))\vert_{s=t}+
\dds\tau(c \st(a \ss b))\vert_{s=t}=$$
$$\phi_t(c,a \st b)+\mu_t(c,(a,b)),
 a,b,c \in \Cal D.$$
In particular this implies that the following
formula,  relating $\psi_t,\theta_t,\phi_t$, holds true:
$$\leqno {(5.1)}\ \psi_t(a,b,c)=
\theta(a,b,c)-1/2\phi_t(a,b,c)=$$
$$ \mu_t(c,(a,b))+1/2\lbrack
-\phi_t(c \st a,b) +\phi_t(c,a\st b)
-\phi_t(b\st c,a)\rbrack.$$
Note that formula (5.1) shows that i). and
iv). imply iii).

The property iv). follows from the identity
$$\dds \tau(d \st ((a \ss b)\ss c))\vert_{s=t}=
\dds\tau(d \st (a \ss (b \ss c))\vert_{s=t},$$
by using the product rule for differentiation.

Finally property v) is a
consequence  of the following two
equalities:
$$\theta_t(a,b,c)=
\ovl{\theta_t(b^{\ast},a^{\ast},c^{\ast})}$$ and
$\phi_t(a,b)=\ovl{\phi_t(b^{\ast},a^{\ast})}$, which hold for
all
$a,b,c \in \Cal D.$
These  properties follow both from corresponding properties
of the trace. This ends the proof.

\bigskip
We will check that the model
described in Theorem 3.2 has the properties
in Definition 5.3. We will obtain bounds for the cocycles
 by estimating the absolute for  the integrals
representing the cocycles.

\proclaim { Proposition 5.5}
Let $a>1$ and let $\Cal D$ be the vector space
$L^2(\Cal A_a) \cap\hat\Cal A_a$. Then $\Cal D$
is weakly dense in
$\Cal A_s,$ for all $s$ in $(a,b)$ and  the
conditions in Definition 5.3 hold true for $\Cal D$.
In particular $\Cal D$ is closed under all the multiplication
the operations $(\ss)_{s \in (a,b)}$. For any $t$ in $(a,b)$
and
$r<a$, the
cocycles $\mu_t,\theta_t,\phi_t$ are defined on
$\Cal D$. Moreover there exists a constant $c$ depending
on $r$ so that
 for all $A,B,C \in \hat\Cal A_r\cap L^2(\Cal A_t)$ we have:
$$\vert\phi_t(A,B)\vert
\leq c\vv A\vv_{2,t}\vv B\vv_{2,r},$$
$$\vert\phi_t(A\st B,C)\vert
\leq c \vv C\vv_{\lambda,r}
\vv A\vv_{2,t}\vv B\vv_{2,t},$$
$$\vert\mu_t(C,(A,B))\vert\leq
c(\vv C\vv_{\lambda,r}
\vv A\vv_{2,t}\vv B\vv_{2,t}+
\vv A\vv_{\lambda,r}
\vv B\vv_{2,t}\vv C\vv_{2,t}+
\vv B\vv_{\lambda,r}
\vv C\vv_{2,t}\vv A\vv_{2,t}).$$
\endproclaim

Proof. We only have the check the inequalities. By
Lesbegue theorem on differentiation
under the integral sign, the
derivatives involved in $\phi_t,\mu_t$ will
exist as soon as the absolute value of the
derivatives of the integrands have finite integral.

For $A,B,C$ in $\hat\Cal A_r\cap L^2(\Cal A_t)$,
we let
$A_t(z,\ovl\eta)=A(z,\ovl\eta)
\lbrack(z-\ovl\eta)/2i\rbrack^{-t},z,\eta \in \Bbb H$
and we use  a similar notation for $B$ and $C$.
We deduce the following
expressions
for $\mu_t(C,(A,B))$ and
$\phi_t(A\st B,C)$ (were by const we denote
$(\text{area\ }F)^{-1}$
or 1 according to the case when $\G$ has finite or
infinite covolume in \pslr):
$$\leqno(5.2) \ \mu_t(C,(A,B))=
\dds \tau(C\st(A\ss B))\vert_{s=t}=
\frac{c'_t} { c_t}\tau(A\st B \st C)+$$
$$c^2_t(\text{const})
\int\limits_{F_z}
\int\limits_{\Bbb H_{\eta}}
\int\limits_{\Bbb H_{\zeta}}
A_t(z,\ovl\eta)B_t(\eta,\ovl\zeta)
C_t(\zeta,\ovl z)m(z,\eta,\zeta)\text{d}\nu_t
(z,\eta,\zeta),$$
with
$$\leqno(5.3)\ m(z,\eta,\zeta)=
\ln\lbrack(\eta-\ovl\eta)/2i\rbrack+
\ln\lbrack(z-\ovl\zeta)/2i\rbrack-
\ln\lbrack(z-\ovl\eta)/2i\rbrack-
\ln\lbrack(\eta-\ovl\zeta)/2i\rbrack, z,\eta,
\zeta\in\Bbb H.$$
Similarly,
$$\leqno(5.4)\ \phi_t(A\st B, C)=
\frac{c'_t} { c_t}\tau(A\st B \st C)+$$
$$c^2_t(\text{const})
\int\limits_{F_z}
\int\limits_{\Bbb H_{\eta}}
\int\limits_{\Bbb H_{\zeta}}
A_t(z,\ovl\eta)B_t(\eta,\ovl\zeta)
C_t(\zeta,\ovl z)
\ln\vert d(z,\ovl\zeta)\vert^2
\text{d}\nu_t
(z,\eta,\zeta).$$
We will cary over only the estimate
for $\mu_t$, since the other one
is
similar. We recall that we used the notation
$d(z,\ovl\zeta)=(\text{Im}\ z)^{1/2}
(\text{Im}\ \ovl\zeta)^{1/2}
\lbrack(z-\ovl\zeta)/2i\rbrack^{-1},z,\zeta \in
\Bbb H.$
 Also recall that the absolute value
$\vert d(z,\ovl\zeta)\vert^2$ depends only on the
hyperbolic distance between
$z,\zeta $ in $\Bbb H$.
Let $l(z,\zeta)$ be the function on $\Bbb H^2$
defined by
$$\leqno \ (5.5)
l(z,\zeta)=\ln (\text{Im}\ z)^{1/2}+
\ln (\text{Im}\ \ovl\zeta)^{1/2}-
\ln \lbrack(z-\ovl\zeta)/2i\rbrack^{-1},z,\zeta \in
\Bbb H.$$
We  obviously have:
$$ \leqno\ (5.6) \
m(z,\eta,\zeta)=
l(z,\eta) +l(\eta,\zeta)-l(z,\zeta),
z, \eta, \zeta \in
\Bbb H.$$

Consequently to show that the absolute value of the integral
in the formula (5.2) is convergent, it is
sufficient to estimate the following
integral (and two other similar ones).
$$ \int\limits_{F_z}\int\limits_{\Bbb D_{\eta}}
\int\limits_{\Bbb D_{\zeta}}
\vert A_t(z,\ovl\eta)\vert
\vert B_t(\eta,\ovl\zeta)\vert
\vert C_t(\zeta,\ovl z)\vert
\vert l(\eta,\zeta)\vert\text{d}\nu_t
(z,\eta,\zeta).$$
For fixed $z$ in $F$, denote
$f_z(\eta)=\vert A_t(z,\ovl\eta)\vert $,
$g_z(\zeta)= \vert B_t(\eta,\ovl z)\vert $
for $\eta \in \Bbb H$.
By Proposition 1.5.a, the functions
$f_z,g_z$ belong to $L^2(\Bbb H,\nu_r)$ for
all $z$ in $F$. Moreover
$$\int\limits_F\vv f_z\vv^2_{L^2(\Bbb H,\nu_r)}=
\int\limits_F\int\limits_{\Bbb H}
\vert A_t(z,\ovl\eta)\vert ^2\text{d}\nu_t(z,\eta)=
\vv A\vv^2_{L^2(\Cal A_t)}.$$
On the other hand if $B$ belongs to
$\Cal A_r$, then
$$\sup\limits_{z \in \Bbb H}
c_r \int\limits_{\Bbb H}\vert B(z,\ovl\zeta)\vert
\vert d(z,\ovl\zeta)\vert^r \text{d}\nu_0(\zeta)
\leq \vv B \vv_{\lambda,r}.$$
Since, if $r<t$, there exists a constant
$c(r)$ such that
$$x^t\ln x\leq c(r)x^r, 0\leq x\leq 1$$
 and since
$\vert l(z,\zeta)\vert\leq \vert
\ln\vert d(z,\ovl\zeta)\vert\vert$ for all
$  z,\zeta \in \Bbb H$, it follows that there
exists a constant $c(r)$ so that:
$$\sup\limits_{z \in \Bbb H}
\int\limits_{\Bbb H}
\vert B(\eta,\ovl\zeta)\vert
\vert l(\eta,\ovl\zeta)\vert^t
\text{d}\nu_0(\zeta)
\leq c(r) \vv B \vv_{\lambda,r}.$$
By Proposition 2.7, it follows that the
kernel $$K_B(\eta,\ovl\zeta)=
\vert \frac {B(\eta,\ovl\zeta)l(\eta,\ovl\zeta)}
{\lbrack (\eta-\ovl\zeta)/2i\rbrack^t}\vert,
\eta,\zeta \in \Bbb H,$$
defines a bounded operator on
$L^2(\Bbb H,\nu_r)$ of operatorial uniform
norm less than $c(r) \vv B \vv_{\lambda,r}.$

 The integral in formula (5.6) is clearly
equal to
 $$\int\limits_{F_z}\int\limits_{\Bbb D_{\eta}}
\int\limits_{\Bbb D_{\zeta}}
K_B(\eta,\ovl\zeta) f_z(\eta), g_z(\zeta)
\text{d}\nu_t(\eta,\zeta) \text{d}\nu_t(z).$$
By the above arguments this is less than
$$\int\limits_{F_z}
\vv K_B\vv_{B(L^2(\Bbb H,\nu_r))}
\vv f_z\vv_{H_r} \vv g_z\vv_{H_r}\text{d}\nu_t(z)\leq$$
$$c(r)\vv B \vv_{\lambda,r}
\lbrack\int\limits_F\vv f_z\vv_{H_r}^2\text{d}
\nu_t(z)\rbrack^{1/2}
\lbrack\int\limits_F\vv g_z\vv_{H_r}^2\text{d}
\nu_t(z)\rbrack^{1/2}=$$
$$c(r)\vv B \vv_{\lambda,r}\vv A\vv_{2,t}
\vv B\vv_{2,t}.$$
This is an upper bound for one of the three terms
that apear in the expression
of $\mu_t$. The other two terms, listed
in the statement of the proposition, are,
by similar arguments, upper bounds for the
other two integrals. This completes the
proof.
\bigskip

We will now determine the expression for the
cyclic 2-cocycle $\psi_t$,
$t \in (a,b)$ that is associated to the
deformation $(\cas)\sinab$.
The expression for $\psi_t$ will be very similar
to the one that appears in the construction
in the paper by Connes and Moscovici.
The formula for
$\psi_t(A,B,C)$ is obtained by  superposition in the integral
formula   for the trace $\tau(A\st B\st C)$ of a $\G-$
invariant Alexander-Spanier cocycle $\theta$.

Recall that by using the notation
$A_t(z,\ovl\eta)=A(z,\ovl\eta)
\lbrack(z-\ovl\eta)/2i\rbrack^{-t},
z,\eta \in \Bbb H$ and similarly for
$B,C$, the formula for $\tau(A\st B\st C)$ is
$$\tau(A\st B\st C)=
c^2_t(\text{const})
\int\limits_{F_z}
\int\limits_{\Bbb H_{\eta}}
\int\limits_{\Bbb H_{\zeta}}
A_t(z,\ovl\eta)B_t(\eta,\ovl\zeta)
C_t(\zeta,\ovl z)
\text{d}\nu_t
(z,\eta,\zeta).$$
The cocycle $\theta$ is a bounded
measurable
function on $\Bbb H^3$
and this may be  used to get
better estimates for $\psi_t$. We find such an estimate
 in the next statement.
 \proclaim{Proposition 5.7}
Let $1<r<t$ and let
$A,B,C$ be arbitrary elements
 in $\hat\car \cap \lar$.
Let $\phi(z,\eta)=
i\text{arg}\lbrack(z-\ovl\eta)/2i\rbrack=$
$\ln \lbrack(z-\ovl\eta)/2i\rbrack-
\ln \lbrack(\eta-\ovl z)/2i\rbrack$ for $ z,\eta
\in \Bbb H.$
Let
$$\theta(z,\eta,\zeta)=
1/2\lbrack \phi(z,\zeta)+
\phi(\zeta,z)+\phi(\eta,\z)\rbrack, z,\eta,\zeta
\in \Bbb H.$$
Clearly $\theta$ is a $\G-$ invariant,
bounded function on $\Bbb H$.
We  use the notation \break $A_t(z,\ovl\eta)=A(z,\ovl\eta)
\lbrack(z-\ovl\eta)/2i\rbrack^{-t},
z,\eta \in \Bbb H$ and similarly for
$B,C$. Then
$$\psi_t(A,B,C)=
1/2(\frac{c'_t}{c_t})\tau(A\st B\st C)+$$
$$
c^2_t(\text{const})
\int\limits_{F_z}
\int\limits_{\Bbb H_{\eta}}
\int\limits_{\Bbb H_{\zeta}}\theta(z,\eta,\zeta)
A_t(z,\ovl\eta)B_t(\eta,\ovl\zeta)
C_t(\zeta,\ovl z)
\text{d}\nu_t
(z,\eta,\zeta).$$
Moreover the following estimates holds for $\psi_t$.
$$\vert \psi_t(A,B,C)\leq\text{const}
\lbrack\vv A\vv_{\lambda,t}\vv B\vv_{2,t}
\vv C \vv_{2,t}\rbrack,\text{\ for\ all\ }
A,B,C \in \hat\car \cap \lar.$$
\endproclaim

Proof. To deduce the expression for
$\psi_t(A,B,C)$ we use the formulae 5.2 and 5.4.
Because
$$\psi_t(A,B,C)=
\mu_t(C,(A,B))+
1/2\lbrack -\phi_t(C\st A,B)+
\phi_t(C,A\st B)-\phi_t(B\st C,A)\rbrack$$
we obtain that
$$\psi_t(A,B,C)=1/2(\frac{c'_t}{c_t})\tau(A\st B\st C)+$$
$$
c^2_t(\text{const})
\int\limits_{F_z}
\int\limits_{\Bbb H_{\eta}}
\int\limits_{\Bbb H_{\zeta}}\gamma(z,\eta,\zeta)
A_t(z,\ovl\eta)B_t(\eta,\ovl\zeta)
C_t(\zeta,\ovl z)
\text{d}\nu_t
(z,\eta,\zeta),$$
where $\gamma $ has the following expression:
$$\gamma(z,\eta,\zeta)=
m(z,\eta,\zeta)+
1/2\lbrack-\ln \vert d(\eta,\ovl\zeta)\vert^2+
\ln \vert d(z,\ovl\zeta)\vert^2-
\ln \vert d(\eta,\ovl z)\vert^2\rbrack.$$
Since $$\ln \vert d(z,\ovl\zeta)\vert^2=
\ln (\text{Im\ } z)+\ln (\text{Im\ } \zeta)-
\ln\lbrack(z-\ovl\zeta)/2i\rbrack -
\ln\lbrack(\zeta-\ovl z)/2i\rbrack,$$
and
$$m(z,\eta,\zeta)=
\ln (\text{Im\ } \eta)
+\ln\lbrack(z-\ovl\zeta)/2i\rbrack-
\ln\lbrack(z-\ovl\zeta)/2i\rbrack
-\ln\lbrack(\eta-\ovl\zeta)/2i\rbrack,
\text{for\ all \ }z,\eta,\zeta \in \Bbb H,$$
one obtains that $\gamma=\theta$.

The estimate for $\psi_t$ is now obtained by the same
procedure as the one used for $\eta_t$ in
the preceding paragraph, with the only difference
that computations are now easier
by the fact that the function is bounded.

Indeed we have to estimate the following integral:
$$ \int\limits_{F_z}\int\limits_{\Bbb D_{\eta}}
\int\limits_{\Bbb D_{\zeta}}
\vert A_t(z,\ovl\eta)\vert
\vert B_t(\eta,\ovl\zeta)\vert
\vert C_t(\zeta,\ovl z)\vert
\text{d}\nu_t
(z,\eta,\zeta)$$ for $A,B,C \in \hat\car \cap \lar$.
 One denotes
for a  fixed $z$ in $F$,
$f_z(\eta)=\vert A_t(z,\ovl\eta)\vert $,
$g_z(\zeta)= \vert B_t(\eta,\ovl\zeta)\vert $
for $\eta,\zeta \in \Bbb H$.
Since $B$ is in $\hat\cat$ it follows that
the kernel on $\Bbb H^2$ defined  by $\zeta,\eta\ra
\vert B_t(\eta,\ovl\zeta)\vert$ defines a bounded operator
$L^2(\Bbb H,\nu_r)$ of norm less than
$\vv B \vv _{\lambda,t}$. But then the above integral is
(modulo a universal constant) less
than
$$\vv B \vv _{\lambda,t}\int\limits_{F_z}
\vv f_z\vv_t\vv g_z\vv_t\text{d}\nu_t(z)=
\vv B \vv _{\lambda,t}\vv A\vv_{2,t}
\vv B\vv_{2,t}.$$
This completes the proof.
\bigskip
In the last part of this paragraph
we will  the use
 the intermadiate symbols $U_r$ on the
algebras $\car$  for yet another aproach
to the construction of the cocycle $\psi_t$.
We will prove that there exists  a dense domain
$\Cal E\subseteq L^2(F)$, which is closed under all the
multiplication operations in all the algebras
$\cas$. Also, we will prove that on
$\Cal E$ the following formula holds true:
$$\dds\tau(U_sf\ss U_sg\ss U_sh)\vert_{s=t}=
\psi_t(U_tf,U_tG,U_th).$$
This formula explains more clearly the reason for
which $\psi_t$ is a cyclic cocycle.

This is because if we transfer the product
operation on $\Cal D$ by
$$f \circ_t g= U^{\ast}_t(U_tF\st U_tg), f,g \in \Cal E$$
and define a trace by $$\tau(f)=
\text{const}\int\limits_{F} f
\text{d} \nu_0,f \in \Cal E$$
and $\tilde\psi(f,g,h)=\psi_t(U_tf,U_tg,U_th)$,
then the above formula will show that:
$$\dds \tau (f\circ_s g\circ_s h)=
\tilde \psi(f,g,h), f,g \in \Cal D.$$

The reason for which $\tilde \psi_t$
 is a cyclic two cocycle is now easy to deduce
 because
 $\tau(f \circ_tg)$ is a  constant
(depending on $f,g$ only). We will use for the next
statement the formalism
 introduced in Definition 5.1 and 5.3.

\proclaim{Proposition 5.8} With the formalism
in Definitions 5.1 and
5.3
let $\Cal D\subseteq\hat B_s\cbs, s\in (a,b)$ be
a "nice differentiable deformation". Let $t$
be an arbitrary in $(a,b)$. We assume
in addition that there exists a Hilbert space
$H$ and unitaries $U_s \ra \lbs$ with the following
additional property:

Whenever  $f,g$  be vectors in $H$ so that if
$U_sf, U_sg$ belong to $\Cal D$ for $s$ in a small
neighbourhood of $t$, then
$$\dds \lbrack U_sf,U_sg\rbrack_{\lbt}=
(-1/2)\phi (U_tf,(U_tg)^{\ast}).$$

Then if $f,g,h$ in $H$ have the
property that $U_s,U_sg,U_sh$ belong to $\Cal D$
for $s$ in a small neighbourhood of $t$, then
$$\dds \tau(U_sf\ss U_sg\ss U_Sh)\vert_{s=t}$$
exists and is equal to
$\psi_t(U_tf,U_tg,U_th)$.
\endproclaim

\centerline
{6. Bounded cohomology and the cyclic 2-cocycle
  of the Berezin's deformation
quantization}
\bigskip
\bigskip

In this section we prove some
 facts about the cyclic 2-cocycle that was
constructed in the previous section for
 a deformation quantization of algebras.
Recall that $\psi_t$ was a cyclic 2-cocycle
 defined on a dense $*$-subalgebra $\hcat$
of the deformation quantization $\cat$ for $\hg$
 constructed in paragraph 3.
 We will
show that the cohomology class of $\psi_t$
in the second cyclic cohomology group $H^2(\hat\cat,\Bbb C)$ ([Co])
is closely related to a canonical element in the bounded
cohomology of the group $\G$.

 In the last part of this paragraph
we will show that a deformation in which the 2-cyclic cocycle
is bounded with respect to the uniform norms from the
von Neumann algebras will have the property that the algebras
in the deformation are isomorphic. Indeed in this
case, by the next paragraph, the cyclic 2-cocycle
vanishes in the cyclic cohomology group of  the von
Neumann algebra.
 We then prove that
 there exists a linear, nonautonomuous differential equation,
with bounded linear operators, whose evolution operator implements the
isomorphism.

 We first recall the integral formulae and the estimates
 that we found in the
last paragraph for the cyclic  2-cocycle $\psi_t$ associated to the
deformation quantization of $\hg$ that we introduced in paragraph 3.
  Let
$ \phi(z,\ovl{\zeta})=
i\text{arg}((z-\ovl{\zeta})/2i)= \text{ln}
((z-\ovl{\zeta})/2i)-\ovl{\text{ln}((z-\ovl{\zeta})/2i)}$,
 for $ z,\zeta\text{ in }
\h$ and let
$$ \ \theta(z,\eta,\zeta)=\phi(z,\ovl{\zeta})+\phi(\zeta, \ovl{\eta})
+\phi(\eta,\ovl{z}),\ \ z,\eta,\zeta \text{ in } \h.$$
Then $\theta$ is a $\G$-invariant continuous function on $\h^3$
 which is an
Alexander-Spanier cocycle. Let $t>1$ and  let $A,B,C$ be elements
 in $\hat\cat$.
Let
 $\hat A=\hat A(z,\eta)$ be the Berezin's, contravariant
symbol of $A$.
 We use the notation $A_t(z,\ovl{\eta})=
\hat A(z,\eta)((z-\ovl{\eta})/2i)^{-t},
\ z,\eta \in  \h$,  and similarly
 for $B$ and $C$. Then the formula for
$\psi_t$ is
$$(6.1)\ \ \psi_t(A,B,C)=
1/2(\frac{c'_r}{c_r})\tau (A\ast_tB\ast_tC)+$$
$$\int\limits_{F_z}\int\limits_{\times}
\int\limits_{\Bbb H^2_{\eta ,\zeta}}
\theta(z,\eta,\zeta)A_t(z,\ovl{\eta})
B_t(\eta, \ovl{\zeta})C_t(\zeta,\ovl{z})
d\nu_t(z,\eta,\zeta).$$
The formula for $\psi_t(A,B,C)$ should be compared with the similar
 formula for the
trace $\tau (A*_tB*_tC)$ which is
$$\tau (A*_tB*_tC)=
\int\limits_{F_z}\int\limits_{\times}\int\limits_{\Bbb H^2_{\eta ,\zeta}}
A_t(z,\ovl{\eta})B_t(\eta, \ovl{\zeta})C_t(\zeta,\ovl{z})
d\nu_t(z,\eta,\zeta).$$
Note that the integral formula for $\psi_t$ like the formula
for $\tau (A*_tB*_tC)$  is an iterated integral. This integral
  converges absolutely
 if $A$ belongs to $\hcat$. In fact for the absolute convergence of the
integral it is sufficient that
$$\sup\limits_{z}\int\limits_{\h}|A(z,\ovl{\zeta})|
\ |d(z,\ovl{\zeta})|^r d\nu_0(\zeta))
\leq ||A||_{\lambda,t}<\infty.$$
Also recall that we found the following estimate:
$$(6.2)\ \ |\psi_t(A,B,C)|\leq \text{const} ||A||_{\lambda,t}
 ||B||_{L^2(\cat)} ||C||_{L^2(\cat)}$$
The relation between the 2-cocycle $\psi_t$ with the deformation
 quantization is
more transparent  from the viewpoint of the
the intermediate symbols $V_rA=B^{-1/2}_r(A(z,\ovl{z}))$, $A$
in $\car$ introduced
at the end of paragraph 4. Recall that
$V_r$ maps $L^2(\car)$ isometrycally onto $L^2(F)$. If $U_r$ is
 the inverse
for $V_r$ then we found that
$$ \psi_t(U_t f,U_t g,U_t h)=
\frac{d}{ds}\tau (U_sf*_sU_sg*_sU_sh)|_{s=t}$$
if $f,g,h$ run through a dense subset of $L^2(F)$.

 We will single out some obstructions
for the cocycle $\psi_t$ to be trivial
 in the cyclic cohomology group
$H^2_{\lambda}(\hat\cat,\Bbb C)$. The  condition that $\psi_t$
 vanishes
in this cohomology group  is equivalent to
the existence of a bilinear form on $\cat$ so that
$$\leqno (6.3)\ \chi_(A,B)=-\chi_t(B,A),$$
$$ \leqno (6.4)\
\psi_t(A,B,C)=
\chi_t(B\st A,C)-\chi_t(B,A\st C) +\chi_t(C\st B, A),$$
for all  $A,B,C$ in $\hat\cat$.
Because of the antisymmetry for $\chi_t$, the
relation (6.4) is equivalent to:
$$ \leqno (6.5)\
\psi_t(A,B,C)=
\chi_t(B\st A,C)+\chi_t(A\st C,B) +\chi_t(C\st B, A).$$
A  natural candidate
 for  a bilinear form $\chi_t$
is given by the formula
$$\chi_t(A,B)=
\int\limits_{F_z}\int\limits_{\Bbb H}
A_t(z,\ovl\eta)B_t(\eta,\ovl z)d(z,\ovl\eta,\eta,\ovl  z)
\text{d}\nu_t(z,\eta),$$
 for a suitable
$\G-$equivariant function $d$ on $\Bbb H^2$.

In fact we wil rather use  this formula to construct  a solution for a
 perturbed problem
with respect to the equation in (6.5).
This is contained in the following statement.
\proclaim{Proposition 6.1}
Let $d=d(z,\ovl\eta,\eta,\ovl  z)$ be  a
$\G-$equivariant function  on $\Bbb H^2$, having purely imaginary
values and with the following properties :
$$
d(z,\ovl\eta,\eta,\ovl  z)+
d(\eta ,\ovl\zeta,\zeta,\ovl\eta)+
d(\zeta,\ovl z,z,\ovl\zeta)=\theta(z,\eta,\zeta),$$
$$
d(z,\ovl\eta,\eta,\ovl  z)=
-d(\eta ,\ovl z,z,\ovl \eta ),$$
$\text {\ for \ all\ }
z,\eta,\zeta \in \Bbb H.$  Let $A,B$ be
 in $\lar$, with Berezin's contravariant symbols
$\hat A,\hat B$.  We use the  notation $A_t(z,\ovl \eta )$ for
$\hat A(z,\ovl \eta)\lbrack (z-\ovl\eta)/2i\rbrack^{-r} $ and
similarly for $B$.
Let $\chi_t(A,B)$ be defined  by
$$\chi_t(A,B)=
\int\limits_{F_z}\int\limits_{\Bbb H}
A_t(z,\ovl\eta)B_t(\eta,\ovl z)d(z,\ovl\eta,\eta,\ovl  z)
\text{d}\nu_t(z,\eta).$$
Then   $\chi_t$ is an antisymmetric
bilinear form ($\chi_t(A,B)$=-$\chi_t(B,A)$).
 The form domain for $\chi_t$ is the linear
space of all $A,B$ for which the integrals
in the definition of $\chi_t$ are absolutely convergent.  Moreover the
following equality holds true for all operators $A,B,C$ in $\lat$ whose
symbols are so that the  integrals  involved in the
formula are absolute convergent:
$$ \leqno (6.6)\
\psi_t(A,B,C)-1/2(\frac{c'_r}{c_r})\tau (A\ast_tB\ast_tC)=
\chi_t(B\st A,C)+\chi_t(A\st C,B) +\chi_t(C\st B, A).$$
\endproclaim

Proof.
Indeed, if the integrals are absolutely
convergent, then we have that:
$$\chi_t(A\st B,C)=
\int\limits_{F_z}\int\limits_{\Bbb H_{\eta}}
\int\limits_{\Bbb H_{\zeta}}
d(z,\ovl\zeta,\zeta,\ovl  z)
A_t(z,\ovl\eta)B_t(\eta,\ovl\zeta)
C_t(\zeta,\ovl z)\text{d}\nu_t(z,\eta,\zeta),$$
and simillarly for the other two terms.

Because the absolute value of the integrands has finite integral and
the integrands are $\G-$ invariant functions on $\Bbb H^3$,
by Fubini's theorem, we may  choose any domain of integration,
as long it is a fundamental domain for $\G$ in $\Bbb H^3$.
 Then $(6.6)$ reduces
 to the first
property for the function $d$  (by  using also formula (6.1)).
This completes the proof.
\bigskip

In general it is difficult to check that the domain were the
above identity holds is sufficiently rich. In fact, in the
case of groups $\G$  of finite covolume, the identity
vector $1\in \cat\subseteq \lat$  makes the integrals involved in the
formulae divergent.

 We will   find  a condition on the group
$\G$ for which a function $d$ with the properties in Proposition 6.1
 exists. To do this we need to  recall
 the constuction of canonical a group 2-cocycle
in the second  cohomology group $H^2(PSL(2,\Bbb R),\Bbb Z)$.

\proclaim {Definition 6.2} Let $N(\gamma_1,\gamma_2),$
$\gamma_1,\gamma_2\in \G$ be the group cocycle in the
second  cohomology group $H^2(PSL(2,\Bbb R),\Bbb Z)$,
defined by the formula:
$$(2\pi)N(g_1,g_2)=
\text{arg}\ j(g_1g_2,z)-\text{arg}\ j(g_1,g_2 z)-
\text{arg}\ j (g_2,z), $$
 for all  $g_1,g_2 \in
PSL(2,\Bbb R)$ and for all $z \in \Bbb H.$

Then $N$ is a non-trivial element in
$H^2(PSL(2,\Bbb Z)$.  The  only possible   values
 for $N$ are  -1 , 0 or 1.
(see e.g the book of Maas ([Ma], pp. 113).
Denote by $N_{\Gamma} $ the restriction of $N$  to
$\G\times \G$. Then $N_{\G} $ vanishes
in $H^2(\G,\Bbb Z)$ if $\G$ is not cocompact (see e.g [Pat]).
\endproclaim

The reasons for which the 2-cocycle $N_{\G}$ is a coboundary in
$H^2(\G,\Bbb Z)$
when $\G$ has finite covolume are more transparent in the
case when $\G =PSL(2,\Bbb Z)$. In this case (because the
commutator subgroup of $PSL(2,\Bbb Z)$ is cyclic of finite order)
there exists at most one $\Bbb Z-$ valued
 cocycle $c=c(\gamma),\gamma\in\G$
so that
$$N_\G(\gamma_1,\gamma_2)=
c(\gamma_1\gamma_2)-c(\gamma_1)-c(\gamma_2),\gamma_1,\gamma_2\in \G.$$
Also, when  $\G =PSL(2,\Bbb Z)$ it is easy to determine the cycle $c$.
  A possible
formula for $c$ is
$$c(\gamma)=\ln(\Delta(\gamma z))-\ln(\Delta(z)),
\gamma \in \G,z \in \Bbb H.$$
The explicit formula for $c$ in terms of the generators for
$\G$ has been determined allready by Rademacher in [Ra].
(Recall that $\Delta$ is the unique modular form for $\G$
of order 12 and that $\ln\Delta$ is defined in $\Bbb H$.)

In the next proposition we find  a sufficient criteria on a discrete,
fuchsian subgroup $\G$ of \pslr\  so that there exists a bounded function
$d$ for $\G$ with the properties in Proposition 6.1. The fact that
$d$ is bounded implies that the bilinear form $\chi_t$ constructed in that
proposition is bounded. As we will see later the criteria on the
group $\G$ can not hold true unless $\G$ has infinite covolume.

\proclaim {Proposition 6.3} Let $\G$ be a fuchsian group
such that $N_\G$ vanishes not only in $H^2(\G,\Bbb Z)$ but also
in the bounded cohomology group $H^2_{\text{bound}}(\G,\Bbb Z)$
(that is there exist a bounded cochain $c:\G \ra \Bbb Z$ so that
$N_\G(\gamma_1,\gamma_2)=
c(\gamma_1\gamma_2)-c(\gamma_1)-c(\gamma_2),\gamma_1,\gamma_2\in \G.$)

Then there exists
a bounded measurable function $d$ on $\Bbb H$
so that the function on $\Bbb H^2$ defined by
$z,\zeta\ra \text{arg\ }\lbrack(z-\ovl\zeta)/2i\rbrack+d(z)-d(\zeta),$
is diagonally  $\G-$ invariant.
\endproclaim

Proof. Define for each $\gamma \in \G$
$$J(\gamma,z)=\text{arg\ }(j(\gamma,z))- (2\pi)c(\gamma),z \in \Bbb H.$$

We clearly have then that
$$\leqno (6.8)\ J(\gamma_1\gamma_2,z)=
J(\gamma_1,\gamma_2 z)+J(\gamma_2,z) \text{\ for\ all\ }
\gamma_1,\gamma_2 \in \G,z \in \Bbb H.$$
Moreover, the quantity
$\sup\limits_{z \in \Bbb H,\gamma \in \G}\vert J(\gamma,z)\vert$
is finite.

We let $d$ to be any bounded measurable function on $F$ and then
  we define
$d$ outside $F$ by the relation
$$d(\gamma z)=d(z)+i J(\gamma,z), z\in F,\gamma \in \G/\lbrace e\rbrace.
$$
The conditon (6.8) shows that in this case
the relation $d(\gamma z)=d(z)+J(\gamma,z),$
will hold true for all $z\in\Bbb H$ and all $\gamma \in \Bbb H$.
 We clearly obtain now
that the following equality
$$ \text{arg\ }\lbrack(\gamma z-\ovl{\gamma \zeta})/2i\rbrack-
\text{arg\ }\lbrack( z-\ovl{ \zeta})/2i\rbrack=
\text{arg\ }(j(\gamma ,z)-\text{arg\ }(j(\gamma,\zeta))=
J(\gamma, z)-J(\gamma,\zeta),$$
 holds true for all $z,\zeta$ in $\Bbb H$ and $\gamma \in \G$. Hence,
with the above choice for the function $d$, the function
$\text{arg\ }\lbrack(z-\ovl\zeta)/2i\rbrack+d(z)-d(\zeta)$ is
diagonally $\G-$ invariant on $\Bbb H$. This completes the
proof.
\bigskip

As we mentioned before the statement of the preceding proposition,
the vanishing of the cocycle $n_{\G}$ in the bounded cohomology
amounts to the fact that the bilinear form in Proposition 6.1 may
be chosen  to be bounded. This is more precisely stated in the
 following
corollary.
\proclaim{Corollary 6.4} If the 2-cocycle $n_\G$ defined in
6.3 vanishes in $H^2_{\text{bound}}(\G,\Bbb Z),$ then
there exists a bounded, antisymmetric
operator $X$ on
$\lat$ so that the bilinear functional $\chi_t$ defined on
$\lat\times \lat$ by
$\chi_t(A,B)=
\langle X(A),B\rangle_{\lat},$ $A,B \in \lat,$
is a solution to the equation (6.6).
\endproclaim

We will now prove that this can only happen if
$\G$ has infinite covolume.

\proclaim{Corollary 6.5}
Let $n_\G$ be the integer valued, two cocycle on the group $\G$,
 defined by the
following relation, in which the choice of $z$ in $\Bbb H$ is
 irelevant:
$$(2\pi )n_\G(\g_1,\g_2)=
\text{arg\ }(j(\g_1\g_2,z)-\text{arg\ } j(\g_1,\g_2 z)-
\text{arg\ } j (\g_2,z), \g_1,\g_2 \in
\G, z\in \Bbb H. $$
If $\G$ has finite covolume then $n_\G$ is a
nonzero element in $H^2_{\text{bound}}(\G,\Bbb Z)$.
\endproclaim

Proof.  Assume that we have a group $\G$ of finite
covolume so that $n_\G$ vanishes in $H^2_{\text{bound}}(\G,\Bbb Z)$.
Then, by the preceding corollary, there exists a
bounded, antisymmetric operator
$X_t$ on $\lat$ so that the equation
(6.7) holds true with
$\chi_t(A,B)=
\langle X_t(A),B\rangle_{\lat}.$ In this case $1$ belongs to
$\lat$, so we may take
in the equation (6.6), $A=B=C=1$.
As $\psi_t(1,1,1)=0$ we obtain that
$$\langle X_t 1,1\rangle_{\lat}=
-\frac{c'_r}{c_r}\langle 1,1\rangle_{\lat}.$$
This contradicts the fact that $X_t$ is antisymmetric.
This completes the proof of the corollary.

\bigskip

Assume  that    $\G$ is such that the cocycle $n_\G$ is zero
in $H^2_{\text{bound}}(\G,\Bbb Z)$. Then the equation (6.6)
shows a better estimate for $\psi_t$.

\proclaim{Proposition 6.6} Assume that $\G$ is a fuchsian
subgroup of \pslr\ (necessary of infinite covolume)
so that the group cocycle   $n_\G$  introduced in
Definition 6.3 vanishes
in $H^2_{\text{bound}}(\G,\Bbb Z)$.
Then, for all $A,B,C \in \lat \cap \cat$, we have
$$\vert\psi_t(A,B,C)\vert\leq\text{(const)}
\lbrack\vv A\vv_2\vv B \st C\vv_2+
\vv B\vv_2\vv C \st A\vv_2 +
\vv C\vv_2\vv A \st B\vv_2\rbrack.$$
In particular $\psi_t$ extends to a 3-linear functional
on $\lat \cap \cat$.
\endproclaim

The discussion shows that  when $\G$ is
of finite covolume, we can not expect to be able
to solve the equation (6.6) with  $\chi_t$ of the form in
Proposition 6.1 and so that simultaneously 1 be contained in the domain of
the quadratic form $\chi_t$.

Finally we show that if don't require for $d$ to be bounded then a
 function
$d$ with the properties in Proposition 6.1 may be easy to construct
 for discrete
groups like $PSL(2,\Bbb Z)$. This corresponds to the fact that
 the Alexander-Spanier cocycle $\theta$ on $\Bbb H^3$ defining
  the cyclic 2-cocycle$\psi_t$
 is a coboundary even in
the $\G-$ equivariant form of the Alexander-Spanier cohomology.

\proclaim{Remark 6.7} Assume that $\Gamma$ is a fuchsian
subgroup of \pslr\ so that  there exists an automorphic
form $\nu$ of order $k$, $k \in 2\Bbb N$, which is nowhere zero
in $\Bbb H$. For example $\G$ could be $PSL(2,\Bbb Z)$ and
$\nu$ could be the unique automorphic form  $\Delta$ of weight
12 for $PSL(2,\Bbb Z)$. Let
$\alpha$ be the function on $\Bbb H^2$ defined by
$$\alpha(z,\ovl\zeta)=
(1/k)\lbrace\ln \nu(z)+\ovl{\ln\nu(\zeta)}+
k\ln\lbrack(z-\ovl \zeta)/2i\rbrack\rbrace,z,\zeta \in \Bbb H.$$

Then $\alpha$ is $\G$-invariant and the hypothesis
of (6.7) are fullfilled with
$$d(z,\ovl\zeta,\zeta,\ovl z)=
\alpha(z,\ovl \zeta)-\alpha(\zeta,\ovl z),z,\zeta \in \Bbb H.$$
Moreover the (unbounded) quadratic form $\chi_t$ associated to
$d$ has the following form
$$\chi_t(A,B^{\ast})=
\langle(\alpha\cdot A),B\rangle_{\lat}-
\langle (A,(\alpha\cdot B)\rangle_{\lat}.$$
We suppose that $A,B$ run through a subspace $D=D(\chi_t)$  of
  $\Cal A_t$ and assume
that $D$ is so that for $A,B$ in $D$
their contravariant symbols
 led to absolutely convergent integrals in the formula
for  $\chi_t$.
\endproclaim
By using the Berezin intermediate symbols $U_t$, that were
introduced in paragraph 4, it is interesting to observe
 the  epression
for $\chi_t(U_tf, U_tg)$, when $chi_t$ is defined
by a function $d$ constructed as above.
We state this separately
\proclaim{Remark} We use the notations from the preceding
remark. Let  $f,g$ in $L^2(F)$ so that $U_tf,U_tg$
belong to the domain $D(\chi_t)$. Let $M_{\alpha}$ be  the
(unbounded) multiplication
operator on $L^2(F)$ with the restriction of the function
of $\alpha$ to the diagonal $z=\zeta$. Then
$$\chi_t(U_tf, U_tg)=\langle B^{-1/2}_t\lbrack B_t,M_{\alpha}\rbrack
B^{-1/2}_tf,g\rangle_{L^2(F)},$$
where $f,g$ belong to the domain of the (unbounded) operator on the right
hand side of the equality.
\endproclaim

Proof. We   only have to  check  the
last formula in the statement of the remark.
Let  $A,B$ be in the domain of $\chi_t$  and assume
 $B=T^r_g$, $A=T^r_f$ for some $f,g$ in $L^2(F)$. Then, by Proposition 4.6,
$$\chi_t(A,B)=
\int\limits_F(\alpha\cdot A)(z,\ovl z)\ovl{g(z)}-
f(z)(\alpha\cdot B^\ast)(z,\ovl Z)\text{d}\nu_0(z)=$$
$$\int\limits_F\alpha(z,\ovl z)(A(z,\ovl z)\ovl {g(z)}-
f(z)\ovl{B(z,\ovl z)}\text{d}\nu_0(z).$$
Let $k,l$ be the intermediate  Berezin symbols (see paragraph 4)
for $A,B$, that is $A=T^t_{B^{-1/2}_tk}$,
$B=T^t_{B^{-1/2}_tl}$.
Then
$f=B^{-1/2}_tk$,$g=B^{-1/2}_tl$ and
$A(z,\ovl z)=B_tf(z)=B^{1/2}_tk(z)$
and
$B(z,\ovl z)=B_tg(z)=B^{1/2}_tl(z)$.
Hence
$$\chi_t(A,B^{\ast})=
\langle M_{\alpha}B^{1/2}_tk,B^{-1/2}_tl\rangle_{L^2(F)}-
\langle M_{\alpha}B^{-1/2}_tk,B^{1/2}_tl\rangle_{L^2(F)}=$$
$$\langle B^{-1/2}_t\lbrack B_t,M_{\alpha}\rbrack
B^{-1/2}_tf,g\rangle_{L^2(F)}.$$
This completes the proof.
\bigskip
The problem with the solutions we have constructed so far
is that we rather solved  the perturbed equation (6.6) instead of
(6.5), with antisymmetric cocycles $\chi_t$. This corresponds
to a scaling factor, which is then cancelled out by numerical
factors in the formulae for the
traces on the algebras $\cas$.

In the remaining part of this paragraph we prove that if one could find
a  solution to the equation (6.5) with bounded antisymmetric cycles
would imply that the  algebras $(\cas)_{s\in (a,b)}$ are isomorphic.
We will state the procedure of constructing such an isomorphism
 in an abstract setting that that formally uses the  intermediate
Berezin symbols introduced in paragraph 4.

\proclaim{Definition 6.8}
Let $H$ be a Hilbert space and let $\Cal E\subseteq H$
be a dense subspace with an involution  denoted
by $\ast$ and a lenght 1 vector denoted by $1$, $1\in \Cal E$.
For $t\in (a,b)$, let $\Cal D_t\subseteq B(H)$ be  type $II_1$ factors
 with  unit $1_{B(H)}$. Moreover,
assume that the trace $\tau_{\Cal D_t}$ on $\Cal D_t$ is
computed by the formula:
$$\tau_{\Cal D_t}(x)=\langle x(1),1\rangle_H,\text{\ for\ all\ }
x \in \Cal D_t\subseteq B(H).$$
When no confusion arrises we  denote the trace
$\tau_{\Cal D_t}$ simply by $\tau$.
In particular $H$ is canonically identified with the
Hilbert  space $L^2(\Cal D_t,\tau)$
 of the Gelfand-Naimark-Segal construction for the
 trace $\tau$ on $\Cal D_t$.
We require that the invololution on $H$ is exactly the one
corresponding to the canonical involution on $L^2(\Cal D_t)$ for
all $t$. For $x,y$ in $L^2(\Cal D_t)\cap \Cal D_t$ we
denote their  product in $\Cal D_t$ by
$x \circ_t y$.
In addition we assume that the subspace $\Cal E$ is
contained in the intersection of all
$L^2(\Cal D_t)\cap \Cal D_t$ for all $t$.

 Let
$\vv\cdot\vv_{\infty,t}$ be the norm defined on a dense subspace
of $H$ which corresponds to the uniform norm
on $\Cal D_t$.
 We assume that the function on $(a,b)$ defined by
$ s\ra\vv f\vv_{\infty,s},$
is locally bounded for all  $f$ in $\Cal E$.
Also we require that
the derivative
$$\tilde \psi_t(f,g,h)=\dds \tau(f\circ_s g\circ_s)\vert_{s=t},$$
exists for all $f,g,h \in \Cal E.$

We will call  a family $(1,\Cal E,H,(\Cal D_t)_{t \in (a,b)}$),
with the  above properties, a
"nice intermediate deformation".
\endproclaim

The reason for this terminology ("nice intermediate deformation")
is that this type of deformation corresponds to the
Berezin quantization deformation, when we use
the intermediate symbols corresponding to the operators $U_r$ acting
on $\lar$. This is explained in the following remark

\proclaim{Remark 6.9}
Assume that $\G\subseteq$ \pslr\ is a fuchsian group of
finite covolume. Let $(\car)_{r \in (a,b)}$ be the family of
von Neumann  algebras that are associated with  the $\G$-invariant form
of the Berezin quantization (see Theorem 3.2).

Let $V_r:\lar\ra L^2(F)$ be the unitary corresponding
to the intermediate symbols defined
Proposition 4.9 and let $\Cal E$ be the dense
subspace constructed at the end of the
paragraph 4.

Let $(\Cal D_r,\circ_r)
=V_r(\car)V^{\ast}_r$ be the type
$II_1$ factor  represented on $L^2(F)$, obtained by
transporting the multiplication structure from $\car$:
$$f\circ_rg= V_r(V^{\ast}_rf\ast_rV^\ast_rg), f,g\in \Cal D_r.$$
Then  $(1,H=L^2(F),\Cal E,(\Cal A_t)_{t \in (a,b)}$), is
a "nice intermediate deformation" in the sense of the
preceding definition.
\endproclaim
\proclaim{Definition 6.10}We use
the notations from Definition 6.8. For $s \in (a,b)$,
 assume that  the cocycles
$\tilde\psi_t$ in Definition 6.8, have in addition  the property that
for all $f,g,h \in \Cal E$,
$$\tilde \psi_t(f,g,h) \vert
\leq(\text{const})\vv f\vv_{\infty,t}
\vv g\vv_{L^2(F)}\vv h\vv_{L^2(F)}.$$
We  then  call the deformation
$(1,\Cal E,H,\ast,(\Cal D_t)_{t \in (a,b)}$)
a "rigid nice intermediate deformation".
\endproclaim
\proclaim{Remark 6.11} The property that the function
$s\ra \vv f\vv_{\infty,s}$ is locally bounded for $f$ in
$\Cal E$ shows that with the additional property
in Definition (6.10), in a "nice intermediate deformation",
the cocycle $\tilde\psi_t$ may be extended by
continuity to $\Cal D_t\times L^2(\Cal D_t)\times L^2(\Cal D_t)$.
Moreover we have that
$$\dds\tau(f\circ_sg\circ_sh)\vert_{s=t}$$ exists and is equal
to
$\tilde\psi_t(f,g,h)$ for all $f$ in $\Cal E$ and
$g,h \in L^2(F)$.
\endproclaim

We will show now that for  a "rigid deformation" as in
Definition 6.10, the algebras $\Cal D_t$ in the deformation
are all isomorphic. The proof will consist into two steps:
one is to show that the boundedness conditions in the
Definition 6.10 imply that $\tilde\psi_t$ is a (bounded) coboundary
in  Connes's cyclic 2-cohomology group. The
other step will be to show that the evolution operator
for  a differential equation associated to the deformation
realizes precisely this the isomorphism.
The precise statement is the following:
\proclaim{Proposition 6.12}
Let $(1,\Cal E,H,\ast,(\Cal D_t)_{t \in (a,b)}$) be a
"rigid nice intermediate deformation" in the sense of
the Definitions 6.8 and 6.10. Assume in addition that there
exists a bounded operator $A(t)$ on $H$ for each $t$ in
$(a,b)$ with the following properties:
\item {(a).} $t\ra A(t)$ is a (norm) bounded measurable
function with values in $B(H)$.
\item{(b)} $A(t)$ is antisymmetric and if we define $\phi_t$ by
$\phi_t(x,y)=\langle A(t)x,y^{\ast}\rangle_{L^2(F)}$ for $ x,y \in L^2(F),$
then
$$\tilde\psi_t(f,g,h)=
\phi_t(f\circ_tg,h)-\phi_t(f,g\circ_th)+\phi_t(h\circ_tf,g),$$
for all $f\in \Cal D_t$ and $g,h \in L^2(\Cal D_t)$.
\item {(c)} In addition, $A(t)$ maps
$L^1(\Cal D_t,\tau)$ and $\Cal D_t$ continuously into
$L^1(\Cal D_t,\tau)$ and respectively $\Cal D_t$.
\item {(d)} $A(t)$ preserves the involution on $H$ and
$A(t)1=0$.

 For $t,s \in (a,b)$, let $U(t,s)$ be the evolution operator
([Sim]) corresponding to the linear, nonautonomuous, diferential
equation :
$$\dot y(t)=A(t) y(t).$$
Then $U(t,s)$ is a unitary for all $t,s$.  By definition $U(t,s)$  has
the property:
$$\dds U(s,t)= A(s)U(s,t).$$
Moreover $U(t,s)$ maps $\Cal D_t$ into
$\Cal D_t$ and $U(t,s)$ is an algebra
isomorphism from the algebra $\Cal D_t$ into $\Cal D_s$.
\endproclaim

We will prove into the next paragraph that the existence
of a bounded measurable function $t\ra A(t)$
for which the properties a), b), c), d) hold true follows
automatically from the boundedness property
for $\tilde\psi_t$ in Definition 6.10
(that is  $\tilde \psi_t(f,g,h) \vert
\leq(\text{const})\vv f\vv_{\infty,t}
\vv g\vv_{L^2(F)}\vv h\vv_{L^2(F)},$
for all $f,g,h \in \Cal E$.)
Proof of Proposition 6.12.
We will divide the proof into a series of lemmas.
\proclaim {Lemma 6.13}
For $f,g$ in $H$, the
derivative
$\dds (f\circ_sg)\vert_{s=t}$ exists in the weak sense
and
$$\dds (f\circ_sg)\vert_{s=t}=
A(t)(f\circ_tg)-A(t)f\circ_tg-f\circ_tA(t)g.$$
Note that the right hand side makes perfectly sense as an
element in $L^1(\Cal D_t,\tau)$.
\endproclaim
Proof. We check the above equality by taking the scalar
product of both terms of the equation with a vector
$h$ in $\Cal E$ and use the fact that
$\dds \tau(f\circ_sg\circ_sh)\vert_{s=t}$ is
equal to $\psi_t(f,g,h)$. We then use   condition (b) from the
hypothesis of Proposition 6.12.
\proclaim {Lemma 6.14}
Let $f$ be any selfadjoint elelment
in $\Cal E$ and fix $t$ in $(a,b)$. Let $\lambda$ be any
complex number with $\text{Im\ } \lambda\ne 0$.
Denote the inverse of an element $a$ in
$\Cal D_s$ (if it exists) by $a^{-1,s}$. Then
$$\dds(f+\lambda)^{-1,s}\vert_{s=t}$$
exists (weakly in $H$) and it is equal to:
$$ A(t)((f+\lambda)^{-1,t})+
(f+\lambda)^{-1,t}\circ_t(A(t)(f+\lambda))\circ_t(f+\lambda)^{-1,t}.$$
\endproclaim
Proof. Formally this follows from
the equality
$$\dds\lbrack(f+\lambda)^{-1,s})\circ_s(f+\lambda)\rbrack\vert_{s=t}=0.$$
This implies that
$$\dds\lbrack(f+\lambda)^{-1,s}\rbrack_{s=t})\rbrack
\circ_t(f+\lambda)=
-\dds \lbrack(f+\lambda)^{-1,t})\circ_s(f+\lambda)\rbrack\vert_{s=t}=$$
$$-A(t)(1)+(A(t)y)\circ_t(f+\lambda)+y\circ_tA(t)(f+\lambda)=$$
$$A(t)y\circ_t(f+\lambda) +y\circ_tA(t)(f+\lambda),$$
which is the required equality (we use the notation
 $y=(f+\lambda)^{-1,t}$).

To obtain a rigirous justification for the above
formal computation we note that the procedure we are
using here is the following:
We  have $F(s,t)$ a function on
$(a,b)^2$ (in our case the function
$(s,t)\ra \rangle(f+\lambda)^{-1,s})\circ_t(f+\lambda),c\rangle_{L^2(F)}$
for a fixed $c$ in $\Cal E$) and we know that
$F(t,t)$ is constant on $(a,b)$.

We want to deduce that
$\dds F(s,t)\vert_{s=t}$ exists
and is equal  to
$-\dds F(t,s)\vert_{s=t},$
if the later term exists.
This  comes from the identity:
$$(s-t)^{-1}\lbrack F(s,t)-F(t,t)\rbrack=
(s-t)^{-1}\lbrack F(s,t)-F(s,s)\rbrack.$$

Thus, the proof would be completed if we can prove that
$(z,y)\ra \frac {\text{d}}{\text{d}z} F(y,z)$ is a continuos
function in $(y,z)$ around $(t,t)$ which in turn will follow if we
knew that
$$(z,y)\ra \tilde\psi_z((f+\lambda)^{-1,y}, f+\lambda,c)$$
is a continuous function $(y,z)$ around $(t,t)$.
This  follows from the following statement
\proclaim {Lemma 6.15} For
$f$ in $\Cal E$ and $g$ in $H$, the map on $(a,b)$ with values in
 $H$ defined by
$$s \ra f\circ_sg$$
is continuous. The same holds for an $n-$ fold product for every
$n\geq 2$ in $\Bbb N$.
\endproclaim
Proof.
Fix $t$ in $(a,b)$ and assume that
$M=\sup\limits_{s \in V}\vv f\vv_{\infty,s}$ is finite in a
neighbourhood $V$ of $t$.

For any $h$ in $H$ and $s$ in $V$ we have that
$$\langle f\circ_sg,h\rangle-
\langle f\circ_tg,h\rangle=
\int\limits^s_t\tilde\psi_p(f,g,h)\text{d}p.$$
Using the estimate
for $\psi_p$ we get that
$$\vert\langle f\circ_sg,h\rangle-
\langle f\circ_tg,h\rangle\vert\leq
M\vert s-t\vert \vv g\vv_H \vv h\vv_H.$$
Since this is valid for arbitrary $h\in H$ the statement follows.

\proclaim {Corollary 6.16}
Let $f$ be any vector in $\Cal E$ and let $\lambda$ be any non real
element  in $\Bbb C$. Then the map
$$s\ra (f+\lambda)^{-1,s}$$ on $(a,b)$ with values in $H$ is
continuous in the norm topology.
\endproclaim
Proof.
We use the
expression:
$$
(f +\lambda)^{-1,s}=
\text{(cst)}\int\limits^{\infty}_0 \exp(-\lambda p)
\exp_{\cas}(itf)\text{d}t.$$
 We also use the preceding corollaary
for the continuous dependence on $s$ of the
function $s \ra \exp_{\cas}(itf)$. (we use here the notation
$\exp_{\cas}$ for the exponential inthe algebra $\cas$.)
This proves the required continuity result and also concludes
the proof of Lemma 6.13.
 From lemma 6.13 we also deduce the following more general result
\proclaim{Lemma 6.17}
Let $f$ be any function on $(a,b)$ with values in $H$ and so that
$f$ is differentiable and $f'(s)$ belongs to $H$ for allmost all
$s$ in $(a,b)$. Let $\lambda$ be any complex number which is not
real.
Then, (in the weak sense)
$\dds(f(s)+\lambda)^{-1,s}\vert_{s=t}$ exists and it is equal to
$$A(t)(f(t)+\lambda)^{-1,t}+(f(t)+\la)\ut\circ_t
\lbrack A(t) f(t)-f'(t)\rbrack\circ_t(f(t)+\lambda)\ut,$$
for almost all $t$ in $(a,b)$.
\endproclaim
Proof. A general fact in Banach algebras shows that:
$$\dds (f(s)+\lambda)\us\vert_{s=t}=
(f(s)+\lambda)\ut \circ_t f'(t) \circ_t (f(s)+\la)\ut, t\in (a,b).$$
This and Lemma 6.13 completes the proof of the Lemma 6.17.

\proclaim{Corollary 6.18} Let $f=f(t)$ be a function
on $(a,b)$ with values in $H$ so that $f$ is
differentiable and so that
$f'(t)=A(t)f(t)$ for almost all $t$. Then, for all
$\lambda \in \Bbb C/\Bbb R$, we have
$$\frac {\text{d}}{\text{d}t}(f(t)+\lambda)^{-1,t}=
A(t)(f(t)+\lambda)^{-1,t},$$
for almost all $t$ in $(a,b)$.
Consequently, the uniqueness in $H$, of the solution for the
equation $\dot y(s)=A(s)y(s)$ shows that
$$U(s,t) (f+\la)\ut =(f+\la)\us,$$
for almost all $s,t$ in $(a,b)$.

Hence, for almost all $s,t$ in $(a,b)$, $U(s,t)$ maps
$\Cal D_t$ into $\Cal D_s$.
\endproclaim

Proof. This follows from Corollary 6.17 and the fact that the
set
$$\lbrace (\la +f)\ut\vert\ f \in H=L^2(\Cal D_t),f=f^{\ast},
\lambda \in \Bbb C/\Bbb R\rbrace$$
is normic dense in $\Cal D_t$ for all $t$.

We have thus proved

\proclaim {Theorem 6.19} Let $(1,\Cal E,H,\ast,(\Cal D_t)_{t \in (a,b)}$),
be a rigid "nice intermediate deformation" as in definitions 6.8 and
6.10.

 Let $(A(t))_{t \in (a,b)}$ be a bounded solution
(measurable in $t$) of the equation
$$\tilde\psi_t(f,g,h)=
\lbrack A(t)(f\circ_tg)-A(t)f\circ_tg-f\circ_tA(t)g,h^{\ast}\rbrack_H,$$
for all $f,g,h$in $\Cal E$ with properties a), b), c), d) in Proposition 6.12
(such a solution exists automatically by the next paragraph).

Let $(U(t,s) )_{t,s \in (a,b)}$ be the evolution operator
associated to the linear differential equation $\dot y(t)=
A(t)y(t)$.
Then $U(t,s) $ maps $\Cal D_t$ onto $\Cal D_s$ and
$U(t,s) $ is an algebra
isomorphism for all $s,t$.
\endproclaim

Proof. Using 6.18, it follows that we only have to check that
$U(t,s) $ is a morphism of algebras. Fix $t$ in $(a,b)$.
For $x,y$ in $\Cal D_t$, we let
$x(s)=U(s,t)x$ and $y(s)=U(s,t)y$. Then we have
$\dot x(s)=A(s)x(s)$ and $\dot y(s)=A(s)y(s)$ and
$x(t)=x,$ $y(t)=y$.

Let $z(s)=x(s)\circ_s y(s),s \in (a,b)$. Then $x(s), y(s)$ belong
to $\Cal D_S$ for all $s$ and
$$\dds z(s)\vert_{s=t}=
\dds (x(s)\circ_sy(s))\vert_{s=t}=
$$
$$\dds (x(s)\circ_t y(s))\vert_{s=t}+
(A(t)x(t))\circ_t+ x(t)\circ_t (A(t)y(t)).$$
By Lemma 6.14 this is
$A(t)(x(t)\circ_ty(t)).$ As $x(t)\circ_ty(t)$ belongs to
$\Cal D_t\subseteq H$ for all $t$ and
$z(t)=x\circ_ty$, the unicity of the solution of the
linear(nonautonomuous) differential equation
shows that
$$z(t)=x(t)\circ_ty(t)=
U(s,t)(x\circ_ty).$$
this completes the proof.

\centerline{7) Vanishing for certains bounded cyclic cohomology}
\centerline
{cocycles
in a finite von Neumann algebra}
\bigskip
The main result of this section is that cyclic cohomology
 2-cocycles $\ps$
on a type $II_1$ factor $M$ with trace $\tau$, that have
 the property that
$$|\ps (a,b,c,)|
\leq (\text{const})\ ||a||_{\infty}\ ||b||_2\ ||c||_2$$
are coboundaries of antisymmetric 1-cococycles on
$M$ defining bounded
operators on $\lt$ and $L^{1} (M,\tau)$.

We will start first by explaining why this result is not
sufficient for our
purposes and then prove the above mentioned result.
The cocycles that  we would like  to be coboundaries in Connes's
cyclic cohomology
 live on dense subsets of the algebras
$\cat$ like $\hat \Cal A_t$.

Let $t$ be any real number in
$(1,\infty)$. Let  $A,B,C$ be  in $L^2(\Cal A_t)$, with contravariant
 symbols $A=\hat A(z,\overline{\zeta}),$ for
$ z,\zeta \text{ in } \Bbb H$ and similarly for $B$ and $C$.
Let
 $ A_t(z,\overline{\zeta})$ be \break $\hat A(z,\overline{\zeta})
\lbrace(z-\overline{\zeta})/2i\rbrace^{-t}$ and similarly
 for $B$ and $C$.
   The 2-cocycle $\psi_t$ associated to the
deformation quantization for $\hg$
 is  defined by the following formula:
(as long as the integrals are absolutely convergent)
$$\psi_t (ABC)=(1/2)(c'_t/c_t)\tau (A*_t B*_t C)+$$
$$+c^2_t\int\limits_{F_z}\int\limits_{\times}
\int\limits_{\Bbb H^2_{\eta ,\zeta}}
i\theta _t (z,\zeta ,\eta)A_t(z,\overline\eta)
B_t(\eta ,\overline\zeta)C_t(\zeta,\overline\z)d\nu_t(z,\zeta , \eta).$$
Recall that $\t$ is a bounded, continuous,
 $\G$ invariant function on
$\h^3$ given by formula (6.3). Note that if in the
last integral above we replace $\t$ by
$1$, then we  get $\tau (A\ast_tB\ast_tC$).  Also recall that
we found a  Banach norm $||\ ||_{\lambda,t}$
 on a weakly dense  subalgebra
$\hcat$ of $\cat$ that behaves nicely with respect to
the algebras in the deformation. Using this norm we found the
estimate
$$|\psi_t(A,B,C)|\leq ||A||_{\lambda,t}||B||_{2,t}||C||_{2,t}.$$
Unfortunately the norm $||A||_{\lambda,t}$ is not equivalent (at least
as long as the Haussdorf dimension
of the limit set $\delta(\G)$ is strictly
greater than $ 1/2$ (see[Pat])) to $||A||_{\infty,t}$.
Moreover,
 even when assuming stronger conditions on the group $\G$, like
vanishing of the canonical cocycle $[n_{\G}]$ (introduced
in Definition (6.2) in the bounded cohomology
 $H^2_{bound}(\G,\Bbb Z)$,
 the best estimate we are able to find for $\psi_t$ is
$$|\psi_t(A,B,C)|\leq \text{\ const\ }||A||_{2,t}||B*_tC||_{2,t}+\text
{\ two other terms by permutation }.$$
This estimate can be improved to a a complete boundedness
condition ([EC,AS]).
\proclaim{Remark 7.1}
Let $\G$ be a discrete, fuchsian subgroup
of \pslr\  such that the canonical cocycle
$[n_{\G}]$, introduced in (6.2),
 vanishes in $H^2_{bound}(\G,\Bbb Z)$.  For $n$ in $\Bbb N$, let
 $\tilde{\psi}_{t,n}$ be the 3-linear functional
on $M_{n}(\cat)\times M_{n}(\cat)\times M_{n}(\cat)$ with values in
$M_n(\Bbb C)$ that is
associated with $\psi_t$. Then  $\tilde\psi_{t,n}$
is defined by requiring that,
 for
$A=(A_{ij}),\ B=(B_{ij}),\ C=(C_{ij})$ in $M_n(\cat)$,
$\tilde{\psi}_{t,n}((A_{ij}), \ (B_{ij}),\ (C_{ij}))$ be
 the matrix with
$i,j$ entries equal to $\sum\limits_{k,l}\psi _{t}(A_{ik},B_{kl},C_{lj})$.
for all $i,j$.
Then
$$\vv\tilde{\psi}_{t,n}(A,B,C)\vv _{M_n (\Bbb C)}\leq \text{const }
(\vv A\vv\  \vv  BC\vv +\vv B\vv\ \vv CA\vv +\vv C\vv \ \vv AB\vv )$$
for all $A, B, C,$ in $M_n(\cat)$, all the norm being uniform norms.
\endproclaim

Proof. Indeed we know that we have in this case a splitting for
$\psi_t(A,B,C)-c'_r/c_r\tau (A*_t B*_t C)$ into a sum of three
other terms:
$$\chi_t(A,B*_tC)+\chi_t(B,C*_tA)+\chi_t(C,A*_tB)$$
with  $\chi_t$  of the form
$$\chi_{t}(A,B)=\int\limits_{F}
\int\limits_{\times\h}(d(z,\zeta,\overline{\zeta},z))
A_t(z,\overline {\zeta })
B_t(\zeta,\overline {z})d\nu_t(\z,z),$$ for a
suitable $\G-$invariant function $d$.
 It is thus  sufficient to prove for $\chi_t$ such a completely
 boundedness  type of estimate.
But $\chi_t (A,B)$ has the following expression
$$\langle \Cal T_d A,B \rangle_{L^2(\cat)}
=
\langle \Cal T_d A *_tB\zeta_0,\zeta_0 \rangle_{L^2(\cat)}, $$
where $\Cal T_d$ is the Toeplitz operator with symbol
 $d(z,\overline{\zeta},\zeta,z)$
on the space of analytic functions on $\h\times\h$,
 that are $\G$ invariant
and square summable and $\zeta_0$ is the unit vector in $L^2(\cat)$.
As $\Cal T_d$ has bounded symbol,the Paulsen and Smith dilation
 lemma for completely
bounded maps applies. (see e.g. [E.K.],[A.S.]). This completes the
proof.
\bigskip
It is conceivable that the techniques in [C.S.] or [Sm.,P.]
 could eventually be
used to show that in this case $\psi_t$
 is a coboundary of a completely bounded cocycle.
To obtain for $\psi_t$ an estimate like the one
in the main theorem of this paragraph, one would need
to have some  more information about the function $z,\zeta\ra
\text{arg\ }[(z-\z)/2i]$ which appears in the expression for $\t$.
\proclaim{Remark 7.2}
If the function $\phi=\phi (z,\zeta)$ on $\h^2$ defined by $(z,\z)\ra
\text{arg\ }[(z-\z)/2i]$ could be shown to belong to the projective
 tensor product
of $L^{\infty}(\h)$ with itself, (or even weaker, if one
could prove that $\phi$  belongs to a weak
limit of some ball in the projective tensor
  product) then it would follow
that $\psi_t$ automatically verifies
 the estimate $$|\psi_t(A,B,C)|\leq
\text{const} ||A||_{\infty,t}||B||_{2,t}||C||_{2,t}.$$
 Thus, if this
would hold true,
 it would
follow that for all lattices $\G$ in \pslr, the algebras
 in the deformation are isomorphic
(i.e. that the fundamental group of $\Cal L(\G)$ is nontrivial).

The same conclusion would also hold if the function $z,\zeta\ra
\text{arg\ }[(z-\z)/2i]$ would be a Schurr multiplier on the
Hilbert spaces
$H^2(\Bbb H,(\text{Im\ }z)^{t-2}\text{d} z\text{d} \ovl z)$
for $t$ in an interval.

\endproclaim
\bigskip
The main result of this paragraph shows that the
estimate\break
$|\psi_t (a,b,c,)|
\leq (\text{const})\ ||a||_{\infty}\ ||b||_2\ ||c||_2,$
implies that the cocycle $\psi_t$ is trivial in cyclic cohomology.

\proclaim{Theorem 7.3}
Let $M$ be a semifinite von Neumann algebra with semifinite, faithful
normal trace  $\tau$. Denote the $L^2$-norm on $M$ by $||\  ||_2$
 and the uniform
norm on $M$ by $||\ ||_{\infty}$.  Let $\psi :(L^2(M)\cap M)^3\ra \Bbb C$
be a 3-linear functional on $M$ with the following properties:

\item{i).} $\psi$ is a cyclic 2-cocycle in the sense of [AC], that is
for all $a,b,c,d$ in $M\cap L^2(M,\tau)$,:
$$\psi(a,b,c)=\psi (b,c,a)$$
$$\psi(ab,c,d)-\psi(a,bc,d)+\psi(a,b,cd)-\psi(da,b,c)=0.$$
\item{ii).} $\ |\psi (a,b,c)|
\leq\text{const}\  ||a||_{\infty}\ ||b||_{2}\ ||c||_{2}\ $
for all $a,b,c$ in $M\cap L^2(M,\tau)$.
 Moreover if $\psi$ is extended
by continuity to $M\times (L^2(M,\tau)\cap M)^2$ then $\psi(1,b,c)=0$,
for all $b,c$ in $M$.
iii) $\psi (a,b,c)=\ovl{\psi (a^*,b^*,c^*)}$ for all
 $a,b,c$ in $M\cap L^2(M,\tau)$.

Then there exists a bilinear
 form $\phi:(M\cap L^2(M,\tau))^2\ra\Bbb C $
so that
for all $a,b,c$
in $M\cap L^2(M)$,
$$\psi(a,b,c)=\phi(ab,c)-\phi(a,bc)+\phi(ca,b).$$
In addition $\phi(a,b)=-\phi(b,a)$ and $\phi$
 may be chosen so that if
$\chi=\chi_{\phi}$ is the linear operator
 on $L^2(M,\tau)$ defined by the equality $
\langle\chi a, b\rangle_{\tau}=\phi(a,b^*), a,b \in
M\cap L^2(M,\tau) $, then $\chi$ is a bounded antisymmetric operator.
Finnally, we may in addition assume that $\chi$ maps  $M$ into
$M$ and that $\chi$ maps $M_{sa}$ into $M_{sa}$.
\endproclaim
Proof. We will consider a convex set $K_{\psi}$ of bounded bilinear
 functionals
on $L^2(M,\tau)$.
By identifying $K_{\psi}$ with a convex compact subset
of the unit ball of $B(L^2(M,\tau))$ we will be able to apply the
fixed point theorem of Ryll-Nardjewski.
This is a standard procedure when solving cohomology problems in
von Neumann algebras (see [Ek], [Ka]).

For each bounded, bilinear functional
 $\phi$ on $(L^2(M,\tau))^2$,
we associate a bounded linear operator $T_{\phi}$ in $B(L^2(M,\tau))$
which is defined
by $\langle T_{\phi}(x),y^*\rangle_{\tau}=\tau(yT_{\phi}(x))
=\phi (x,y)$,
$x,y$ in  $L^2(M,\tau)$.
For $u$ a unitary in $M$ let $\phi_u$ be the bounded linear functional
defined by
$$\phi_u (x,y)=\psi(yu^*,u,x),\ \  x,y\text{\ in\ } (L^2(M,\tau)$$
and let $T_u$ be the associated bounded operator on $(L^2(M,\tau)$ which
is thus defined by
$$\langle T_{u}(x),y\rangle_{\tau}=\phi_u (x,y^*)=
\psi(y^*u^*,u,x),\ \  x,y\text{\ in\ } (L^2(M,\tau).$$
To simplify our setting we will assume that the constant
 in (ii) is 1 so that
$T_u$ belongs to the unit ball $B(L^2(M,\tau)_1$ for all unitaries
$u$.
We  consider now the weakly compact convex set $K$ in $B(L^2(M,\tau)_1$
defined by
$$K=\ovl{co}^w\{T_u\ |\ u\in \Cal U(M)\}.$$
The $w$- topology is the weak operator topology on the unit ball
 $B(L^2(M,\tau)_1$.

We show  that for $T$ in $K$, $T$ extends to a bounded linear operator
on $M$. We first note the identity:
$$(7.3)\  \langle T_{u v}(x),y\rangle_{L^2}=
\langle T_{u}(vx),vy\rangle_{L^2}
- \phi_u (v,xyv^*)+\langle T_{ v}x,y\rangle$$
which is valid for all $x,y$ in $L^2(M)$ and for all
$u,v$ in $\Cal U(M)$.
This is easy to check because it corresponds to
$$\psi(y^*v^*u^*,uv,x)=
\psi(y^*v^*u^*,u,vx)-\psi(xyv^*u^*,u,v)+\psi(yv^*,v,x).$$
This is equivalent to
$$\psi(yv^*,v,x)-\psi(y^*v^*u^*,uv,x)+
\psi(y^*v^*u^*,u,vx)-\psi(xyv^*u^*,u,v)=0,$$
which is exactly the identity for $\psi$ with
$a=y^*v^*u^*,\ b=u,\ c=v,\ d=x.$

 From the relation (7.3),
 by using the continuity for $T_u,\ T_v,\ T_{uv}$ we deduce  that
$$ (7.4)\ |\phi_u (v,xyv^*)|\leq 3||x||_2\ ||y||_2\ \text{ for any }
x \ \text{in }L^2(M,\tau).$$
Take $z$ be arbitrary in $L^1(M,\tau)\cap M,$ and let
 $ z=z_+-z_-$ be the canonical decomposition of $z$ as a difference
of positive elements.

The preceding relation shows that
$|\phi_u (v,z_{\underline{+}})|
\leq 3||z_{\underline{+}}||_{L^1(M,\tau)}$
and hence that
$$|\phi_u (v,z)|\leq 3(||z_{+}||_{L^1(M,\tau)}+||z_{-}||_{L^1(M,\tau)})=
3||z||_{L^1(M,\tau)}.$$
Thus
we
 have shown that for $u$ in $\Cal U(M)$
 the bilinear maps $\phi_u$ on $L^2(M,\tau)$,
 have in addition the property that
$$ |\phi_u (v,z)|\leq 3(||z||_{L^1(M,\tau)})\ \text{ for all }
z \ \text{in }L^1(M,\tau)\cap M.$$
We now use the fact that any $x$ in $M$ is
 a linear combination of four
unitaries (see e.g. 2.24 [S.Z.])
$x={\sum\limits^4_{i=1}}\lambda_i u_i$
$ \ \text{ with }|\lambda_i|\leq 2||x||.$
Consequently $|\phi_u (x,z)|\leq 24||x||_{\infty}\ ||z||_1$ for all $x$
in $L^1(M,\tau)\cap M$, $z$ in $L^1(M,\tau)\cap M.$ Since $\psi$ was
assumed to be weakly continuous this shows that the operator
 $T_u$ associated
to $\phi_u$ maps boundedly $M$ into $M$ and $L^1(M,\tau)$
 into $L^1(M,\tau)$.
The operator norm is bounded in both cases by 24.

But then the same statement holds true for any convex combination
in the operators
 $\{T_u\ |u\in \Cal U(M)\}$.
By taking weak limits one obtains any element $T$ in $K$ has in addition
the property that it extends
 by continuity to $M$ (and $L^1(M,\tau)$) and
$$||T(x)||_{\infty}\leq 24 ||x||_{\infty},\ ||T(x)||_1\leq 24 ||x||_1
\text{  for all } x\ \text{in } L^2(M,\tau)\cap M.$$
We define a family of affine weakly
 continuous maps $(\alpha_v),\ v\in \Cal U(M)$
on $K$ with values in $K$ by
$$\langle\alpha_v(T)x,y\rangle_{L^2}=\langle T(vx),vy\rangle_{L^2}-
\langle T(v),vyx^*\rangle_{L^2}+\langle T_u(x),y\rangle_{L^2},\ x,y
 \ \text{in } L^2(M,\tau)\cap M$$
which, by identifying the elements
 in $T$ with the associated bilinear
$\phi$ functionals $\phi$ in $K$, are
$$\alpha_v(\phi)(x,y)=\phi(vx,yv^*)-\phi(v,xyv^*)+
\phi_u(x,y),\ x,y\in M.$$
Relation (7.3) shows that
$$\alpha_v(T_u)=T_{uv}, \ \text{  for all }
 x\ \text{in } \Cal U(M).$$

By what we have just shown $\alpha_u$ are weakly
 continuous and well defined on $K$.
Assume that there exists a common fixed point $\phi$ in $K$
 for all the maps
$(\alpha_u)_{u\in \Cal U(M)}$.
Then $\alpha_u(\phi)(x,y)=\phi(x,y)$
 for all $u$ in $\Cal U(M)$, $x,y,$
in $L^2(M,\tau)\cap M$ which gives the following relation valid for all
$x,y,$ in $L^2(M,\tau)\cap M$, $v$ in $\Cal U(M)$.
$$\phi(x,y)=\phi(vx,yv^*)-\phi(v,xyv^*)+\psi(yv^*,v,x).$$
This is equivalent to:
$$\psi(yv^*,v,x)=\phi(v,xyv^*)-\phi(vx,yv^*)+\phi(x,y)$$
for all $x,y,$ in $L^2(M,\tau)\cap M$, $v$ in $\Cal U(M)$.

Denote $\phi^{op}(x,y)=\phi(x,y)$. We get
$$\psi(yv^*,v,x)=\phi^{op}(xyv^*,v)-\phi^{op}(yv^*,vx)+\phi^{op}(y,x)$$
$$=\phi(yv^*v,x)-\phi^{op}(yv^*,vx)+\phi^{op}(xyv^*,v).$$
Denoting $a=yv^*, \ b=v,\ c=x$ we get
$$\psi(a,b,c)=\phi^{op}(ab,c)-\phi^{op}(a,bc)+\phi^{op}(ca,b)$$
By continuity we get that $\phi^{op}$ has the property
$$\psi(a,b,c)=\phi^{op}(ab,c)-\phi^{op}(a,bc)+\phi^{op}(ca,b)$$
for all $a,b,c $ in $ L^2(M,\tau)\cap M$ (in fact all $b$ in $M$, $a,c$
in $ L^2(M,\tau)$).

We note that in addition all the elements in $K$ have the
 property that
$T1=0$.  Indeed, because we have a proved that
 $\phi_u$ makes sense on $M\times L^1(M,\tau)$ it follows that
$\psi(uy^*,u,x)$ is defined for all $x$ in $M$ and $y$ in $L^1(M,\tau)$.
By weak continuity since $\psi(a,b,1)=0$ for all $a,b$
 in $ L^2(M,\tau)\cap M$
it follows that the same holds true for $a=uy^*,\ b=u$ and hence
$\psi(uy^*,u,1)=0$ for all $y$ in $M\cap L^1(M,\tau)$.

To show that there exists a fixed point for all
 the affine maps $(\alpha_u),\ u\in \Cal U(M)$
on $K$ we will apply the Ryll-Nardjewski theorem. To do
 that we have to find a
seminorm $p$ which is weakly inferior semicontinuous and
 so that if
$T,S\in K$ then
 $\inf\limits_{u\in \Cal U(M)}p(\alpha_u (T)-\alpha_u (S))>0$.

We choose the seminorm $p$ to be the uniform norm on $M$ (since  the
ball \break
$\{x\  |\ ||x||<c\}$ is allways weakly closed).
We have to show that if $T,S$ belong to $K$ and
$$(7.4)\ \inf\limits_{u\in \Cal U(M)}
||\alpha_u (T)-\alpha_u (S)||_{B(L^2(M,\tau))}=0$$
then $R=S$. Denote $R=T-S$.
 Then $\alpha_u (R)$ is given by the formula
$$\langle\alpha_u(R)x,y\rangle_{L^2(M,\tau)}=
\langle u^* R(ux),y \rangle_{L^2(M,\tau)}
- \langle R(u),uyx^* \rangle$$
for all $x,y $ in $L^2(M,\tau)$.

By (7.4) it follows that for any $\epsilon >0$ there exists $u$ in $u$
in $\Cal U(M)$ with
$$||u^*R(ux)-u^*R(u)x||_{L^2(M,\tau)}
\leq \epsilon\  ||x||_2, \  x\text{ in }
L^2(M,\tau).$$

Equivalently this means that for any $\epsilon$ there exists $u$ with
$$||R(ux)-R(u)x||_{L^2}\leq \epsilon ||x||_{L^2},
 \ \text{ for all }x\text{ in }
M\cap L^2(M,\tau)$$
and consequently
$$||R(x)-(R(u)u^*)x||_{2,L^2(M,\tau)}\leq \epsilon ||x||_2$$
As $R(u)$ belongs to $M$ (since we have shown this
 for all $T_u$), this
shows that $R$ belongs to the uniform norm closure
 in $\Cal L(M,M)$ of the
maps $\{L_a \ | \ a\in M\}$ where $L_a(x)=ax$ for $x$ in $M$.
 But this uniform
norm closure is $\{L_a \ | \ a\in M\}$ itself and hence $R$ is of
 the form
$L_a$ for same $a$ in $M$. As $R(1)=0$ it follows that $R=0$.
 Thus the
Ryll-Nadjewski applies.

The relation
$\psi(a,b,c)=\ovl{\psi(b^*,a^*,c^*)}$  gives that
$$ \psi(a,b,c)=\ovl{\psi(b^*,a^*,c^*)}=\ovl{ \phi (b^*a^*,c^*)}-
\ovl{ \phi (b^*,a^*c^*)}+\ovl{ \phi (c^*b^*,a^*)}.$$
Consequently,
 by using also the properties of  $\psi$ and writting\break
$\psi(a,b,c)=\ovl{\psi(b^*,a^*,c^*)}=0$ we get
$$[\phi(ab,c)-\ovl{ \phi (b^*,a^*c^*)}]
\cdot [\phi(a,bc)-\ovl{ \phi (c^*b^*,a^*)}]
+[\phi(ca,b)-\ovl{ \phi (b^*,a^*c^*)}]=0,$$
for all $a,b,c$ in $M\cap L^2(M,\tau)$.
 By
using the property $\phi (x,y)=-\phi (y,x)$, we get that
$$[\phi(ab,c)+\ovl{ \phi (c^*,b^*a^*)}]-[\phi(a,bc)+
\ovl{ \phi (c^*b^*,a^*)}
+[\phi(ca,b)-\ovl{ \phi (b^*,a^*c^*)}]=0.$$
Now
using the expression $\phi=\frac{\phi_1 +\phi_2}{2}$ where
$\phi_1 (x,y)=\phi(x,y)+\ovl{\phi (y^*,x^*)},$
it follows that there is no loss of generality when
 assuming  the equation
$\psi(a,b,c)=\phi(ab,c)-\phi(a,bc)+\phi(ca,b)$ to suppose that
$$\phi =\phi_2 =\phi (x,y)-\ovl{\phi (y^*,x^*)}.$$
Hence $\phi(x,y)+\ovl{\phi (y^*,x^*)}=0, \ x,y \text{ in }L^2(M,\tau)$
which is the required condition for $\phi$ to be antisymmetric.

The only property that hasn't been checked is that
 $\chi(M_{sa})\subseteq M_{sa}$
which is equivalent to
$$\chi (x^*)= \chi (x^*)\ \text{ for all }x\  \text{ in }M\cap L^2(M,\tau)
\text{ or}$$
$$\langle \chi (x^*)y,z\rangle=\langle y,\chi (x)z\rangle \text{ or}$$
$$\langle \chi (x^*),zy^*\rangle=\langle yz^*,\chi (x)\rangle .$$
Denoting $b=zy^*$, this means that we want that
$$\langle \chi (x^*),b\rangle=\langle b^*,\chi (x^*)\rangle .$$
By the antisymmetry for $\chi$ this is
$$\langle \chi (x^*b),b\rangle=-\langle \chi (b^*),x \rangle$$
which is the same as $\phi(x^*,b^* )=-\phi (b^*,x^*)$.

This completes the proof of our theorem.

\centerline{References}

\item {[Aro]} Aronsjan, Theory of reproducing kernels I,
 Memoirs of the
American Mathematical Society, Providence, 1953.

 \item{[AS]}            M.F. Atyiah, W. Schmidt,
 A geometric construction
of the discrete series for semisimple Lie groups, Invent. Math. 42,
(1977), pp. 1-62.

\item{[Ba]} V. Bargmann, Annals of Mathematics, 1950.

\item {[Be1]} F.A. Berezin, Quantization, Math USSR Izvestija,
Vol 8(1974), No 5.

\item {[Be2]} F.A. Berezin, Quantization in complex
  symmetric spaces,
Math USSR Izvestija, 9(1975),  341-379.

\item{       [Be3]}      F. A. Berezin,  General
 concept of quantization, Comm. Math. Phys. 40
(1975), 153-174.

\item {[Cob]}  C. A. Berger, L. A. Coburn, Heat
 Flow and Berezin-Toeplitz
estimates, American Journal of Math., vol 110, (1994), pp. 563-590
\item{ [Co1]}, A. Connes, Noncommutative diff

\item{[Co2]}       A. Connes, Non-commutative Differential Geometry,
Publ. Math., Inst. Hautes Etud. Sci., 62 (1986), pp. 94-144.

\item{[CM]}      A. Connes, H. Moscovici, Cyclic Cohomology,
the Novikov Conjecture and Hypebolic Groups, Topology 29 (1990),
pp.  345-388.

    \item{[Co3]} A. Connes, Sur la Theorie
 Non Commutative de l'Integration,
Algebres d'Operateurs, Lectue Notes in Math., Springer Verlag,
vol 725, 1979.

\item{[Ch1]}. E. Christensen, F. Pop, A. M. Sinclair, R. R. Smith,
On the cohomology groups of certain finite von Neumann
 algebras, Preprint 1993.

\item{[CE] }      E. Christensen, E. G. Effros, A. M. Sinclair,
 Completely Bounded Multilinear Maps and $C^*$-algebraic cohomology,
Invent. Math., 90 (1987), 279-296.

\item{[Dix]}      Diximier

\item{[Dyk]}       K. Dykema,

\item{[Gh]}    J. Barge,   E. Ghys, Cocycles d'Euler et de Maslov,
Math. Ann., 294, (1992), pp. 235-265.

\item{[Gr]}       M. Gromov, Volume and bounded  cohomology,
Publ. Math., Inst. Hautes Etud. Sci. 56, (1982), pp. 5-100.

\item{[Go]}     R.  Godement

\item{[GHJ]}       F. Goodman, P. de la Harpe, V.F.R. Jones, Coxeter
Graphs and Towers of Algebras, Springer Verlag, New York, Berlin,
Heidelberg, 1989.

\item{[Jo] }     V.F.R. Jones, Manuscript Notes (unpublished)

\item {[Ka]}  R. V. Kadison, Open Problems in Operator Algebras,
Baton Rouge Conference, 1960, mimeographed notes.
\item{[KL]}    S. Klimek, A.Leszniewski, Letters in Math. Phys., 24
(1992), pp. 125-139.

\item{[Ku]}    T. Kubota, Theory of Eisenstein series,

\item{[La]}      S. Lang, $SL(2,\Bbb R)$,

\item{[Leh]}J.  Lehner, Discontinuous Groups and Automorphic Forms,
Amer. Math. Soc., Providence, 1964.

\item{[Ma]} Maass,

\item{[NT]}     Neszt, Tsigan

\item{[Pa1]} S. J. Patterson, Spectral Theory and Fuchsian Groups,
Math. Proc.  Camb. Phil. Soc., (1977), vol 81, pp. 59-75.

\item{[Pa2] }Patterrsom , Cohomology discrete groups,
 Glasgow Math. Journal

\item{[Pau]}     V. I.  Paulsen , R.R. Smith, Multilinear maps and
tensor norms on operator systems, J. Funct Analysis,
 29m (1978), 397-415.

\item{[Pi]}      Pisier

\item{[Pu]}      L. Pukanszki,

\item{[Rade]} H.  Rademacher, Zur Theorie
 der Dedekindschen Summen, Math. Zeit.,
63 (1956), pp.  445-463.

\item{[Ra1]}      F. Radulescu, Inv

\item{[Ra2]}       F, Radulescu, Toeplitz..

\item{[Ri]}      M. A. Rieffel, Deformation Quantization and Operator
Algebras,  Proc. Symp. Pure Math. 51, (1990), p. 411-423.

\item{[Ro]} A. Robert, Introduction to the representation theory
of compact and locally compact groups, Cambridge Univ. Press, 1983.

\item{[Sa]} Sakai, $C^*$ and $W^*$ Algebras, Springer Verlag

\item{[Sa]}      Sally, Memoirs A. M. S.,

\item{[Sch]}  Schatten

\item{[Se]}      A. Selberg, Harmonic Analysis
 and discontinuous groups in weakly
symmetric Riemannian spaces with applications
 to Dirichlet series. J. Indian
Math. Soc. 20 (1956),  47-87.

\item{[Sha]} Shapiro, Shields

\item{[Si]} B. Simon,

\item{[SZ]}       Stratila , Zsido, Lectures on von Neumann Algebras

\item{[Up]} D. Borhwick, A. Lesniewski, H.
 Upmeier, Non-perturbative Deformation
Quantization of Cartan Domains, J. Funct. Analysis, 113 (1993),
pp. 153-176

\item{[Ve]}       Venkov,

\item{[Vo1]}       D. Voiculescu, Free entropy.., Inv. Math.

\item{[Vo2]}       D. Voiculescu,  Random matrices.., Inv Math.

\item{[We]}     A. Weinstein

\end